\documentclass[11pt]{amsart}
\usepackage[margin=1.05in]{geometry}
\usepackage{amsmath,amssymb,amsthm,mathtools}
\usepackage{enumitem}
\usepackage{booktabs}
\usepackage{tikz}
\usepackage[protrusion=true,expansion=false]{microtype}
\usepackage[hidelinks]{hyperref}
\AddToHook{env/remark/begin}{\emergencystretch=8em}
\hypersetup{pdftitle={Linear-exponent bounds for the automatic complexity of rational diagonals modulo prime powers},pdfauthor={Alex Shvets}}

\newtheorem{theorem}{Theorem}[section]
\newtheorem{proposition}[theorem]{Proposition}
\newtheorem{lemma}[theorem]{Lemma}
\newtheorem{corollary}[theorem]{Corollary}

\theoremstyle{definition}
\newtheorem{definition}[theorem]{Definition}
\newtheorem{example}[theorem]{Example}
\newtheorem{remark}[theorem]{Remark}

\numberwithin{equation}{section}

\newcommand{\Z}{\mathbb{Z}}
\newcommand{\QQ}{\mathbb{Q}}

\newcommand{\F}{\mathbb{F}}
\newcommand{\Zp}{\mathbb{Z}_p}
\newcommand{\one}{\mathbf{1}}
\newcommand{\zs}{z_{*}}
\newcommand{\Diag}{\operatorname{Diag}}
\newcommand{\supp}{\operatorname{supp}}
\newcommand{\conv}{\operatorname{conv}}
\newcommand{\CT}{\operatorname{CT}}
\newcommand{\kerp}{\operatorname{ker}_p}
\newcommand{\dg}[3]{\langle #1\rangle_{#2,#3}}
\newcommand{\Hc}{\mathcal{H}}
\newcommand{\Lc}{\mathcal{L}}
\newcommand{\Sig}{\Sigma_p}

\title[Linear exponents for kernels of rational diagonals]{Linear-exponent bounds for the automatic complexity\\ of rational diagonals modulo prime powers}
\author{Alex Shvets}
\address{Independent researcher, Haifa, Israel}
\email{alex@shvets.io}
\thanks{ORCID: \href{https://orcid.org/0009-0005-9802-379X}{0009-0005-9802-379X}}
\date{28 November 2026}

\begin{document}

\begin{abstract}
Let $p$ be a prime and let $G$ be the diagonal of a rational function $N/Q$ over the $p$-adic integers, with $Q(0)$ a unit. We prove that the $p$-kernel of its coefficients modulo $p^\alpha$ has size $K_\alpha(G)\le p^{(|J|+1)\alpha-1}+(p^\alpha-1)/(p-1)$, where the finite set $J$ depends on $p$ and the Newton polytopes of $N$ and $Q$. In particular, $|J|\le\prod_i(h_i+1)$ for $h_i=\max(\deg_{z_i}N,\deg_{z_i}Q)$. The proof counts the coefficient vectors of the polynomial numerators obtained after reading successive digits. Once the last $\alpha-1$ digits are fixed, the preceding vectors determine all but the highest base-$p$ digit of each coefficient in the next vector. This gives an exponent linear in $\alpha$ for fixed input and prime. For the Furstenberg series of an equation of bidegree $(h,d)$ this gives $K_\alpha(F)\le p^{((h+1)d+1)\alpha-1}+(p^\alpha-1)/(p-1)$; at $\alpha=1$ this is the known bound $1+p^{(h+1)d}$, and the asymptotic bounds of Bridy and of Rowland, Stipulanti and Yassawi are sharper there. For the root $F$ with $F(0)=0$ of $bx^2F^2+(1+ax+dx^2)F+cx^2=0$ with odd $2$-adic coefficients $a,b,c,d$, the general bound gives $K_\alpha(F)\le 2^{6\alpha-1}+2^\alpha-1$. We prove the sharper bounds $9\cdot8^\alpha/(64\alpha-32)\le K_\alpha(F)\le 2^{3\alpha-1}+3(2^{\alpha+1}+\alpha+1)(2^\alpha-1)$, so the exact limit of $\log_2 K_\alpha(F)/\alpha$ is $3$. For every choice of odd $a,b,c,d$, the orbit of the initial state of the Rowland--Yassawi automaton under the zero digit has exactly $2^{\alpha+1}+\alpha+1$ elements. When $a=b=c=d=1$, this is the formula conjectured by Rowland and Yassawi.
\end{abstract}

\maketitle

\section{Introduction}\label{sec:intro}

\subsection{The counting problem}\label{subsec:intro-count}
Fix a prime $p$ and an integer $\alpha\ge1$. Think of a machine that receives the base-$p$ digits of an index $n$, starting with the units digit, and returns a coefficient modulo $p^\alpha$. How many different states does it need as the precision $\alpha$ increases?

For a sequence $a=(a(n))_{n\ge0}$, reading $k$ digits with value $v$ leaves the machine with the sequence $(a(p^k n+v))_{n\ge0}$ still to evaluate. The collection of these sequences, over all $k\ge0$ and $0\le v<p^k$, is the \emph{$p$-kernel} of $a$. Two different remaining sequences can be told apart by some continuation of the input. Conversely one can use the remaining sequence itself as a state. Extra leading zeros do not change the output. For a series $G(x)=\sum_{n\ge0}g(n)x^n$, let $K_\alpha(G)$ be the size of the $p$-kernel of $(g(n)\bmod p^\alpha)_{n\ge0}$. This is the minimum number of states.

Our inputs are rational functions $N/Q$ in $m\ge1$ variables $z=(z_1,\ldots,z_m)$. Their coefficients lie in $\Zp$, the ring of $p$-adic integers, and the constant coefficient $Q(0)$ is a unit, written $Q(0)\in\Zp^\times$, meaning that it is nonzero modulo $p$. Thus $N/Q$ has a formal power series expansion at the origin. Write $\one=(1,\ldots,1)$ and $z^u=z_1^{u_1}\cdots z_m^{u_m}$ for an exponent vector $u$. The \emph{diagonal} keeps the coefficients whose exponents are all equal:
\[
  \Diag\Bigl(\sum_{u\in\Z_{\ge0}^m}s_u z^u\Bigr)
    =\sum_{n\ge0}s_{n\one}x^n .
\]
Over $\F_p=\Z/p\Z$, Christol's theorem \cite{Chr79} says that a series is algebraic precisely when its coefficient sequence has a finite $p$-kernel \cite[Theorem~1]{CKMR80}. In the $p$-adic setting, Christol's results show that the orbit of a rational diagonal under the section operators, which select the coefficients in a given residue class modulo $p$, is relatively compact for the uniform $p$-adic norm on coefficients \cite[Corollary~8.3 and Example~9]{Chr86}; hence it is finite modulo each $p^\alpha$. Denef and Lipshitz prove that the coefficients of every multivariate algebraic power series over $\Zp$ are $p$-automatic modulo each $p^\alpha$ \cite[Theorem~3.1(i)]{DL87}. In 2015 Rowland and Yassawi gave explicit digit transitions \cite[Theorem~2.1]{RY15}. We ask how the minimum number of states grows as $\alpha$ increases for fixed $N,Q$ and $p$.

\subsection{A bound at every precision}\label{subsec:intro-bound}
The exponent in our bound depends on a finite set of monomials. For a polynomial, its \emph{support} is the set of exponent vectors with nonzero coefficient, and its \emph{Newton polytope} is their convex hull. Put $\Delta(Q)=\conv(\supp Q)$ and $K(N)=\conv(\supp N)$, with $K(0)=\emptyset$. In the definition below, $(1-\varepsilon)\Delta(Q)+\varepsilon K(N)$ consists of the points $(1-\varepsilon)a+\varepsilon b$ with $a\in\Delta(Q)$ and $b\in K(N)$.

\begin{figure}[htbp]
\centering
\begin{tikzpicture}[x=0.9cm,y=0.65cm]
  \fill[gray!25] (0,0) -- (2,1) -- (2,3) -- cycle;
  \draw[black,line width=0.8pt] (0,1) -- (2,3) -- (2,4) -- cycle;
  \draw[->] (-0.20,0) -- (2.60,0) node[right] {$X$};
  \draw[->] (0,-0.20) -- (0,4.62) node[above] {$Y$};
  \foreach \i in {0,1,2}
    \draw (\i,0) -- (\i,-0.07) node[below,font=\footnotesize] {$\i$};
  \foreach \j in {1,2,3,4}
    \draw (0,\j) -- (-0.05,\j) node[left,font=\footnotesize] {$\j$};
  \foreach \i in {0,1,2}
    \foreach \j in {0,1,2,3,4}
      \fill[gray!55] (\i,\j) circle[radius=1.6pt];
  \foreach \i/\j in {0/0,1/1,1/2,2/2,2/3}
    \fill[black] (\i,\j) circle[radius=2.4pt];
\end{tikzpicture}
\caption{The Newton polygons $\Delta(Q_F)$ (grey) and $K(N_F)$ (outlined) for odd $a,b,c,d$, and the support set $J$ (black dots). The lattice points $(0,1)$, $(2,1)$, $(2,4)$ of the polygons are not in $J$.}\label{fig:J}
\end{figure}
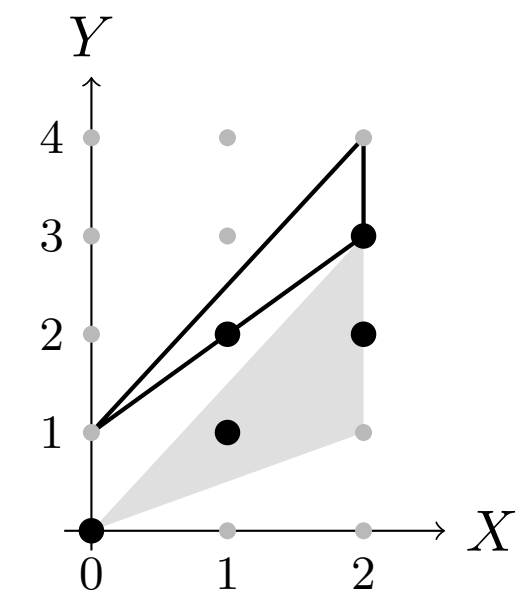

\begin{definition}\label{def:J}
The \emph{support set} of $(N,Q)$ at the prime $p$ is
\[
J=J_p(N,Q)=\bigl\{u\in\Z_{\ge0}^m:\ u+\eta\one\in(1-\varepsilon)\Delta(Q)+\varepsilon K(N)\ \text{for some }\eta\in[0,1),\ \varepsilon\in(0,1/p]\bigr\}.
\]
\end{definition}
The set $J$ is finite and can be computed by linear programming; for the quadratic family of Section~\ref{subsec:intro-family} with odd $a,b,c,d$ it has five points, shown in Figure~\ref{fig:J}.

Write $h_i=\max(\deg_{z_i}N,\deg_{z_i}Q)$ for the largest degree in coordinate $i$, taking the degree of the zero polynomial to be $-\infty$. Every point in $J$ lies in the integer box $0\le u_i\le h_i$, so $|J|\le D=\prod_{i=1}^m(h_i+1)$. We call the first inequality below the \emph{support bound}, since it uses $|J|$, and the resulting bound in terms of $D$ the \emph{degree bound}.

\begin{theorem}\label{thm:main}
Let $p$ be a prime, $m\ge1$, and $N,Q\in\Zp[z_1,\dots,z_m]$ with $Q(0)\in\Zp^\times$. Let $G=\Diag(N/Q)$. Then for every $\alpha\ge1$
\[
K_\alpha(G)\ \le\ p^{(|J|+1)\alpha-1}+\frac{p^\alpha-1}{p-1}\ \le\ p^{(D+1)\alpha-1}+\frac{p^\alpha-1}{p-1}.
\]
\end{theorem}

\emph{Where $J$ comes from.} Divide $N$ and $Q$ by $Q(0)$ so that $Q(0)=1$. After $k$ digits with value $v$, keep the terms of $N/Q$ whose exponents are congruent to $v$ in every coordinate modulo $p^k$, then subtract $v$ and divide by $p^k$ in each coordinate. The resulting series is a sum of fractions with successive denominators $Q(z),Q(z^p),Q(z^{p^2}),\ldots$, where powers of $z$ are taken coordinate by coordinate. Write $m_k$ for the numerator over the first denominator $Q(z)$ in the expansion of Proposition~\ref{prop:levels}. Let $C$ keep the monomials whose exponents are all divisible by $p$ and divide those exponents by $p$. Write $\zs=z_1\cdots z_m$ and $\sigma(A)(z)=A(z^p)$. Define the polynomials $\delta_s$ recursively by
\[
 Q^{p^s-1}=\sum_{t=1}^s\delta_t\,\sigma^t\!\left(Q^{p^{s-t}-1}\right) \qquad(s\ge1).
\]
They give the expansion
\[
 \frac1{Q(z)}=\sum_{t\ge1}\frac{\delta_t(z)}{Q(z^{p^t})},
\]
converging $p$-adically (Lemma~\ref{lem:inverse}). Let $v_{i,t}$ be the value of the last $t$ digits among the first $i$ digits read, with the earliest of those $t$ digits in the units place. Starting with $m_0=N$, we have
\[
 m_i=\sum_{t=1}^{i}T_{t,v_{i,t}}(m_{i-t}),\qquad
 T_{t,v}(M)=C^t\!\left(\zs^{-v}\delta_tM\right)\qquad(1\le i\le k).
\]
For $k\ge1$ the support of $m_k$ lies in
\[
 \frac{(p^k-1)\Delta(Q)+K(N)-v\one}{p^k}.
\]
A recurrence term for a block of $t$ digits adds an exponent from $(p^t-1)\Delta(Q)$, subtracts the value of the block in each coordinate, then divides by $p^t$. The word thus gives $\varepsilon=p^{-k}$ and $\eta=v/p^k$. Letting these parameters range over the intervals of Definition~\ref{def:J} gives one support set $J$ independent of the length $k$ and the value $v$.

To prove the bound we retain the last $\alpha-1$ digits (the \emph{window}) and the last $\alpha$ numerator vectors (the \emph{history}). Write these vectors from oldest to newest, modulo $p,p^2,\ldots,p^\alpha$. Each has $|J|$ coefficients. Counting them independently would give $p^{|J|\alpha(\alpha+1)/2}$ choices. Instead, fix the window and fill in one vector at a time. The first vector contributes one digit per coefficient. For the vector stored modulo $p^j$, the earlier vectors determine its reduction modulo $p^{j-1}$. Indeed, the term $T_{t,v_{i,t}}(m_{i-t})$ reaching back $t$ steps contains Vlasenko's polynomial $\delta_t$, whose coefficients are divisible by $p^{t-1}$ (Lemma~\ref{lem:delta}, \cite[Lemma~3(i)]{Vl18}). Only the highest digit of each coefficient is free. This gives $p^{|J|\alpha}$ compatible histories per window (Lemma~\ref{lem:count}). Multiplying by the $p^{\alpha-1}$ windows and adding the shorter words gives the bound.

For $k\ge\alpha$, let $v=v_{k,k}$ be the value of the entire $k$-digit word and $v_t=v_{k,t}$ the value of its last $t$ digits. Write $[z^u]$ for the coefficient of $z^u$. The output $g(v)$ is the constant coefficient of the section. The term over $Q(z)$ contributes $m_k(0)$, and the terms over $Q(z^{p^j})$ for $j\ge1$ give the sum below (Corollary~\ref{cor:automaton}):
\[
 g(v)\equiv m_k(0)+\sum_{t=1}^{\alpha-1}[z^{v_t\one}]
 \Bigl(m_{k-t}\bigl(\delta_{t+1}+\cdots+\delta_\alpha\bigr)\Bigr)\pmod{p^\alpha}.
\]

\subsection{Algebraic series and the odd quadratic family}\label{subsec:intro-family}
Suppose $P(x,y)$ is a polynomial, $P(0,0)=0$, and the partial derivative $\partial_yP(0,0)$ is a unit. Its \emph{Furstenberg series} is the unique formal series $F(x)$ with $F(0)=0$ and $P(x,F(x))=0$: the coefficient of each new term of $F$ is determined by a linear equation with unit coefficient $\partial_yP(0,0)$. Furstenberg gave a diagonal formula for this root over a field \cite[Proposition~2]{Fur67}. Denef and Lipshitz extended it to integral domains \cite[Lemma~6.3(ii)]{DL87}, so over $\Zp$ we can use the numerator $y(\partial_yP)(xy,y)$ and denominator $P(xy,y)/y$. Section~\ref{sec:cor} proves the next two bounds.

\begin{corollary}\label{cor:furst}
Let $P\in\Zp[x,y]$ with $P(0,0)=0$ and $\partial_yP(0,0)\in\Zp^\times$, $\deg_xP\le h$ and $\deg_yP\le d$, and let $F$ be its Furstenberg series. Then for every $\alpha\ge1$
\[
K_\alpha(F)\le p^{((h+1)d+1)\alpha-1}+\frac{p^\alpha-1}{p-1}.
\]
\end{corollary}

In the next corollary $a,b,c,d$ are coefficients. In particular, $d$ multiplies $x^2y$.

\begin{corollary}\label{cor:family}
Let $a,b,c,d\in\Z_2$, $P=bx^2y^2+(1+ax+dx^2)y+cx^2$, and let $F$ be its Furstenberg series. Then
\[
K_\alpha(F)\le2^{6\alpha-1}+2^\alpha-1\qquad(\alpha\ge1).
\]
If $a,b,c,d$ are odd, then also $K_\alpha(F)\ge2^{\alpha+1}$.
\end{corollary}

When $a,b,c,d$ are all odd (units of $\Z_2$), the support set for this diagonal representation is
\[
 J=\{(0,0),(1,1),(1,2),(2,2),(2,3)\}
\]
(Section~\ref{subsec:proof-family} and Remark~\ref{rem:J-equality}).
The five corresponding monomials are $1,xy,xy^2,x^2y^2,x^2y^3$, as in Section~\ref{subsec:matrix-example}. The support bound therefore gives slope $6$ in the base-$2$ exponent. The bidegree of the defining equation is $(2,2)$, so taking $h=d=2$ in Corollary~\ref{cor:furst} gives the exponent $7\alpha-1$. Theorem~\ref{thm:odd-exponent} proves that the true slope is $3$:
\[
 \frac{9\cdot8^\alpha}{64\alpha-32}\le K_\alpha(F)
 \le 2^{3\alpha-1}+3(2^{\alpha+1}+\alpha+1)(2^\alpha-1),
 \qquad \lim_{\alpha\to\infty}\frac{\log_2 K_\alpha(F)}{\alpha}=3 .
\]

Rowland and Yassawi conjectured in 2024, in the first arXiv version of \cite{RY26} (arXiv:2408.00750v1), that for $a=b=c=d=1$ the zero-digit orbit of their initial polynomial state has exactly $2^{\alpha+1}+\alpha+1$ elements, a count they checked through $\alpha=14$ \cite{RY26}. After the change of coordinates $x\mapsto xy^{-1}$, which we call central-row coordinates, we expand in powers of $x$ with coefficients in the Laurent polynomial ring $\Z_2[y,y^{-1}]$. Taking the coefficient of $y^0$ gives the diagonal. Put $f=P/y$, $q=2^{\alpha-1}$, and choose $\widetilde f$ with the same monomials as $f$ and with $\widetilde f\equiv f\bmod 2$. Rowland and Yassawi build their automaton from such a lift of $f$ modulo $2$ \cite[\S2]{RY26}, so the count should not depend on the choice. A polynomial state is a numerator $S$ over $\widetilde f^q$ modulo $2^\alpha$. The initial numerator is $(y\,\partial_yP)f^{q-1}$. The zero transition is
\[
 S\longmapsto\Lambda_{0,0}(S\widetilde f^{\,q})\bmod 2^\alpha,
\]
where $\Lambda_{0,0}$ keeps terms with both exponents even and halves their exponents. Theorem~\ref{thm:zero} applies to every choice of four odd parameters. Modulo $2^\alpha$, the map $S\mapsto S/\widetilde f^q$ identifies the polynomial orbit with the orbit of $(y\,\partial_yP)/f$ under $\Lambda_{0,0}$ counted in part~(a). Thus the number of distinct polynomial states in this zero orbit, which are states of their automaton rather than of the minimal one, is $2^{\alpha+1}+\alpha+1$ for every permitted lift $\widetilde f$. Taking the coefficient of $y^0$ in each fraction gives a sequence of one-variable series in the kernel. Part~(b) proves that this sequence has exact eventual period $2^{\alpha+1}$, so $K_\alpha(F)\ge2^{\alpha+1}$ throughout the odd family.

To see where the slope $3$ comes from, use central-row coordinates and allow both digits. Once a word leaves the three exceptional paths $0^k$, $1^k$, $01^{k-1}$ ($k\ge1$), each new numerator $m_k$ in central-row coordinates belongs to the rank-two module $W_0=\{c_0+c_1x:c_0,c_1\in\Z_2\}$ (in the coordinates of the diagonal, $W_0$ is spanned by $1$ and $xy$). After another $\alpha-1$ digits the whole history lies there. Two coefficients give $2^{2\alpha}$ compatible histories for each of the $2^{\alpha-1}$ windows, yielding the upper-bound term $2^{3\alpha-1}$. Words on an exceptional path, or with fewer than $\alpha-1$ digits since they left one, account for the second term (Proposition~\ref{prop:odd-upper}).

For the lower bound, fix a window $w$ of $\alpha-1$ digits. The histories with this window that lie in $W_0$ give a plane $\mathcal G_{\alpha,w}$: a $\Z/2^\alpha$-module spanned by two series whose first two coefficients are $(1,0)$ and $(0,1)$. We call these two series the normalized pair of $w$ (Section~\ref{subsec:odd-frames}). The $3\cdot4^{\alpha-1}$ primitive elements of each plane are reachable, so all of them belong to the kernel (Lemma~\ref{lem:odd-reachability}). Planes of different windows overlap little. If $w\ne v$ have a longest common final segment of $\ell$ digits, their normalized pairs can agree modulo $2^s$, with $s\le\alpha$, only if $s\le\ell+3$. We call this inequality the \emph{contact bound} (Lemma~\ref{lem:odd-weak-contact}). Since windows with a long common final segment are rare, a Cauchy--Schwarz count shows that the planes together contain at least $9\cdot8^\alpha/(64\alpha-32)$ kernel elements. The contact bound is proved for $a=b=c=d=1$ by a finite exact computation together with a proof that it suffices at every precision (Appendix~\ref{app:odd-certificate}); a script in the arXiv ancillary files reproduces the computation. Proposition~\ref{prop:odd-universality} shows that the normalized pairs modulo $2^\alpha$ are the same for all odd $a,b,c,d$, so this one computation covers the whole family.

Table~\ref{tab:main} gives finite computations of the entire kernel. For $a=b=c=d=1$, for instance, its sizes modulo $2$ and $4$ are $14$ and $53$. Whether $K_\alpha(F)/8^\alpha$ has a limit remains open.

\subsection{Comparison with earlier work}\label{subsec:intro-prior}
\emph{What is new.} For a fixed input and prime, the general bounds of Rowland and Yassawi modulo $p^\alpha$ have exponents cubic in $\alpha$ for Furstenberg series and of degree $m+1$ in $\alpha$ for diagonals in $m$ variables \cite[Theorems~1, 3 and~56]{RY26}; the count based on Beukers' schemes, discussed below, also has exponent of order $\alpha^{m+1}$. Theorem~\ref{thm:main} gives an exponent linear in $\alpha$ for every rational diagonal with unit denominator at the origin; we do not know an earlier bound with this property for general rational diagonals. The $\delta$-polynomials and the divisibility $\delta_t\in p^{t-1}\Zp[z]$ are due to Vlasenko \cite[Lemmas~3 and~4]{Vl18}. What we add is the observation that, once the last $\alpha-1$ digits are fixed, the earlier numerators determine all but the top base-$p$ digit of each coefficient of the next one. Theorem~\ref{thm:zero} proves the conjecture of Rowland and Yassawi on the zero orbit and extends it to all odd parameters, and Theorem~\ref{thm:odd-exponent} shows that for this family $\log_2K_\alpha(F)=3\alpha+o(\alpha)$, while the general bound gives $6\alpha-1$.

\emph{Reduction modulo a prime.}
Setting $\alpha=1$ in Theorem~\ref{thm:main} gives $K_1(G)\le p^{|J|}+1$. There is also a short direct bound in terms of the degrees alone. Over $\F_p$, one can keep the denominator $Q$ and replace the numerator $S$ on input $r$ by the section of $SQ^{p-1}$ selecting exponent vectors congruent to $(r,\ldots,r)$ modulo $p$. A section divides the selected exponents by $p$ after subtracting that residue. Starting from $S=N$, the inequalities $\deg_{z_i}S\le h_i$ are therefore preserved at every step. Counting the $p$ choices for each of the $D$ coefficients gives the elementary bound $p^D$. In two variables this space of numerators is the one used in the proof of \cite[Theorem~35]{RSY25}.

For diagonals in two variables, Rowland, Stipulanti and Yassawi prove the bound $(1+o(1))p^{h_1h_2}$, where $h_i=\max\{\deg_{z_i}(N\bmod p),\deg_{z_i}(Q\bmod p)\}\ge1$ is the larger of the two degrees in variable $z_i$, for $i=1,2$ \cite[Theorem~35]{RSY25}. The error $o(1)$ tends to zero when one of $p,h_1,h_2$ tends to infinity with the others fixed. If both Newton polygons are the full rectangle $[0,h_1]\times[0,h_2]$, take lifts with that same support. Setting $\eta=0$ puts every lattice point of the rectangle in $J$. Our exponent $(h_1+1)(h_2+1)$ is then always larger than $h_1h_2$. For sparse inputs $|J|$ can be much smaller: the diagonal of $1/(1-x(1+y+y^2))$ has $J=\{0\}$ (Section~\ref{sec:questions}, item~4). For three or more variables, setting $\alpha=1$ in \cite[Theorem~3]{RY26} gives the degree count $\prod_i(h_i+1)$.

For a nonconstant one-variable algebraic series $F$ over $\F_p$, let $h$ and $d$ now be the degrees in $x$ and $y$ of its minimal polynomial equation $A(x,y)=0$, and let $\gamma$ be the genus of the smooth projective model of this curve. Bridy bounds the number of states by $(1+o(1))p^{h+d+\gamma-1}$ \cite[Theorem~1.2]{Bri17}; since $\gamma\le(h-1)(d-1)$, this gives $(1+o(1))p^{hd}$ \cite[Corollary~1.4]{Bri17}. For a nonzero Furstenberg series over $\F_p$ defined by a polynomial of bidegree $(h,d)$, Rowland, Stipulanti and Yassawi give an explicit upper bound and the asymptotic bound $(1+o(1))p^{hd}$ \cite[Theorems~1 and~2]{RSY25}. On the way they prove, with an elementary invariant space, the explicit bound $p^{(h+1)d}+1$ \cite[Propositions~8 and~9 and Corollary~10]{RSY25}. Adamczewski, Bostan and Caruso bound, over a perfect field of characteristic $p$, the dimension of a space that contains an algebraic series and is stable under the section operators, in terms of the Newton polygon of $A$ \cite[Theorem~1.2]{ABC23}; over an arbitrary field of characteristic $p$ they pass to its perfect closure \cite[Remark~1.3]{ABC23}, \cite[Theorem~2.3 and Remark~2.4]{ABC26}. Over $\F_p$ this gives the same explicit bound $1+p^{(h+1)d}$ \cite[Theorem~4.1]{ABC23}. Its exponent $(h+1)d$ exceeds the exponent $hd$ of the asymptotic bounds. When the same degrees apply, Corollary~\ref{cor:furst} gives exactly $p^{(h+1)d}+1$ at $\alpha=1$; its proof in Section~\ref{subsec:proof-furst} shows why. The leading exponent in the bivariate bound of Rowland and Yassawi depends only on the coordinate degrees \cite[Theorem~56]{RY26}.

\emph{Increasing the precision.}
Now fix the prime and the input and let $\alpha$ increase. To compare degrees across precisions, take inputs whose coordinate degrees are preserved under reduction modulo $p$. For a nonzero Furstenberg series of an equation with positive degrees in both variables, the exponent in the asymptotic bound of Rowland and Yassawi is cubic in $\alpha$ \cite[Theorem~1]{RY26}. For a rational diagonal in $m\ge2$ variables with $h_i\ge1$, the exponent in their degree bound is $\sum_{k=1}^\alpha\prod_i(kh_i+1)$, of degree $m+1$ in $\alpha$ \cite[Theorem~3]{RY26}. Their two-variable refinement has leading exponent $h_1h_2\,\alpha(\alpha+1)(2\alpha+1)/6$ \cite[Theorem~56]{RY26}. Theorem~\ref{thm:main} gives a linear exponent for every fixed rational diagonal with unit denominator at the origin. For the odd quadratic family, setting $h=d=2$ in the explicit bound of Rowland and Yassawi \cite[Theorem~2]{RY26} gives $\alpha(\alpha+1)(7\alpha+5)/6$ as the exponent of its main term. Our support exponent $6\alpha-1$ and the exponent $7\alpha-1$ of Corollary~\ref{cor:furst} both fall below it from $\alpha=2$ on ($11$ and $13$ versus $19$). Both inequalities persist at higher precisions (Remark~\ref{rem:numbers}).

Beukers bounds the number of components in a $p$-linear scheme for the coefficients of a rational function with integer coefficients (the components form a vector of sequences $A$, with $A(pn+r)=M_rA(n)\bmod p^\alpha$ for each digit vector $r$ and a fixed matrix $M_r$). Suppose each exponent of $N$ is bounded coordinate by coordinate by some point in the interior of $\Delta(Q)$. His bound is $\rho^m\prod_i\deg_{z_i}Q$ components, where $\rho$ is a positive integer satisfying $\rho-\lceil\rho/p\rceil\ge\alpha-1$ \cite[Corollary~4.2]{Be24}. Rowland and Yassawi \cite{RY26} observe that counting component values modulo $p^\alpha$ gives an automaton with at most $p^{\alpha\prod_i(\rho h_i+1)}$ states for $\deg_{z_i}N\le\deg_{z_i}Q=h_i$ and $h_i\ge1$. Counting the values of Beukers' $\rho^m\prod_i\deg_{z_i}Q$ components directly gives the exponent $\alpha\rho^m\prod_i\deg_{z_i}Q$, which is at most $\alpha\prod_i(\rho h_i+1)$ when $\deg_{z_i}Q=h_i$. For fixed $p$ the least admissible positive $\rho$ grows linearly with $\alpha$, so that exponent has order $\alpha^{m+1}$ at fixed degrees. We obtain a linear exponent by counting the histories that satisfy the congruences between successive numerators.

For constant terms $A(n)=\CT[P^nQ]$ of Laurent polynomials there is another route. Rowland and Zeilberger construct automatic and linear $p$-schemes modulo $p^\alpha$ \cite{RZ14}, and Straub bounds the degrees of the polynomials in these schemes \cite[Theorem~2.4]{Str22}. For the Motzkin numbers this gives an automatic scheme with at most $(p^{\alpha-1}-1)/(p-1)+p^{\alpha(p^{\alpha-1}+2)}$ states, doubly exponential in $\alpha$, and a linear scheme with at most $(p^\alpha-1)/(p-1)+2$ states \cite[Example~2.7]{Str22}. Henningsen and Straub observe that a sequence with $A(0)=1$ satisfies Lucas congruences modulo $p$ exactly when it has a linear $p$-scheme with a single state modulo $p$ \cite[Proposition~1.3]{HS22}. The generating series of the Motzkin numbers is the diagonal of $(1-y^2)/(1-x(1+y+y^2))$, for which $J=\{(0,0),(0,1)\}$, so Theorem~\ref{thm:main} gives $K_\alpha\le p^{3\alpha-1}+(p^\alpha-1)/(p-1)$.

The generating series $G$ of the central trinomial coefficients $[y^0](y^{-1}+1+y)^n$ is the diagonal of $1/(1-x(1+y+y^2))$, for which $J=\{0\}$ and Theorem~\ref{thm:main} gives $K_\alpha(G)\le p^{2\alpha-1}+(p^\alpha-1)/(p-1)$. Kohen constructs, for constant terms $\CT[P^nQ]$ modulo $p^\alpha$, a linear representation in which the first $\alpha-1$ digits select the initial vector \cite[Proposition~23]{Ko26}; it serves his study of uniform recurrence \cite[Theorem~33]{Ko26}, not the size of the kernel. For $q=p^{\alpha-1}$ the terms $y^q$ and $qy^{q-1}$ in $(y^{-1}+1+y)^q$ both survive modulo $p^\alpha$, so no common factor $p$ can be removed from the exponents, and his index interval $[1-q,q-1]$ gives length $L=2p^{\alpha-1}-1$ \cite[Lemma~25 and Definition~26]{Ko26}. Counting vectors crudely, as we do here and not as Kohen does, a representation of length $L$ has at most $p^{\alpha L}$ possible vector states, an exponent exponential in $\alpha$.

A linear digit representation of fixed dimension would immediately give a linear exponent by counting vectors. Proposition~\ref{prop:rank} explains why we count histories: for an irrational D-finite series $F$ with integer coefficients and a fixed prime $p$, the minimum number of generators of a $\Z/p^\alpha\Z$-module containing the $p$-kernel of $F\bmod p^\alpha$ tends to infinity. This is a short consequence of the theorems of Becker and B\'ezivin. A uniform bound would make $F$ $p$-regular over $\QQ$ (its $p$-kernel would span a finite-dimensional space), hence $p$-Mahler by Becker's theorem \cite[Theorem~1]{Bec94}. B\'ezivin proved that a $p$-Mahler series satisfying a nonzero linear differential equation with polynomial coefficients, as a D-finite series does, is rational \cite{Bez94}; see also \cite[Theorem~1]{BCR13}. For a one-variable rational series, a fixed linear recurrence for the coefficients already gives a linear exponent by an elementary count (Section~\ref{sec:rank}).

\section{Preliminaries: kernels, Cartier sections and digit matrices}\label{sec:setting}

\subsection{Kernels}\label{subsec:kernels-prelim}
Let $\Sig=\{0,1,\dots,p-1\}$. A word $w=r_0r_1\cdots r_{k-1}\in\Sig^k$ is read least significant digit first: $v(w)=\sum_{q=0}^{k-1}r_qp^q$ and $|w|=k$. Every integer $0\le v<p^k$ equals $v(w)$ for exactly one $w\in\Sig^k$. For a power series $g=\sum_ng(n)x^n$ and $r\in\Sig$ put $\Lambda_rg=\sum_ng(pn+r)x^n$, and $\Lambda_w=\Lambda_{r_{k-1}}\circ\dots\circ\Lambda_{r_0}$. Induction on $k$ gives
\begin{equation}\label{eq:Lambdaw}
\Lambda_wg=\sum_{n\ge0}g(p^kn+v(w))x^n .
\end{equation}
Indeed, $\Lambda_r\Lambda_wg=\sum_n g(p^k(pn+r)+v(w))x^n$ and $p^kr+v(w)=v(wr)$. Hence
\begin{equation}\label{eq:Kalpha}
K_\alpha(g)=\bigl|\{\Lambda_wg\bmod p^\alpha:\ w\in\Sig^*\}\bigr|,
\end{equation}
and $K_\alpha(g)$ depends only on $g\bmod p^\alpha$.

For completeness, the residual series $\Lambda_wg$ themselves form an automaton: the digit-$r$ transition is $S\mapsto\Lambda_rS$ and its output is $[x^0]S$. Conversely, if two words reach the same state of any automaton generating $g$, then continuing with the digits of any $n$ shows that the two residual sequences have the same coefficient at $n$. Thus distinct residuals require distinct reachable states, and the residual automaton attains this lower bound. This is the classical identification, going back to Eilenberg, of the $p$-kernel with the minimal automaton reading the least significant digit first; see \cite[Proposition~4.2]{ABC23}.

\subsection{Series in several variables}
Write $z^u=\prod z_i^{u_i}$ for $u\in\Z^m$ and $\zs=z^{\one}=z_1\cdots z_m$. Let $\Lc=\bigcup_{c\ge0}\zs^{-c}\Zp[[z]]$ be the ring of formal Laurent series in $z$ whose exponents are bounded below in every coordinate. Define $\Zp$-linear maps $\sigma,C:\Lc\to\Lc$ by
\[
\sigma\Bigl(\sum_us_uz^u\Bigr)=\sum_us_uz^{pu},\qquad C\Bigl(\sum_us_uz^u\Bigr)=\sum_us_{pu}z^u .
\]
Then $\sigma$ is a ring endomorphism fixing $\Zp$. For $r\in\Sig$ and $S\in\Zp[[z]]$ put
\[
C_r(S)=C(\zs^{-r}S),
\]
and for a word $w=r_0\cdots r_{k-1}$ put $C_w=C_{r_{k-1}}\circ\dots\circ C_{r_0}$ (with $C_\emptyset=\mathrm{id}$). We call the $C_r$ \emph{synchronous sections}: they select the same residue $r$ in every variable.

We use the $p$-adic topology on $\Zp[[z]]$: $S_n\to0$ if for every $e$ we have $S_n\in p^e\Zp[[z]]$ for all large $n$. The ring $\Zp[[z]]$ is complete and separated for this topology. All maps used below ($\sigma$, $C$, $C_r$, multiplication by a fixed series) are $\Zp$-linear and map $p^e\Zp[[z]]$ into itself. They are therefore continuous, commute with $p$-adically convergent sums, and are compatible with reduction modulo $p^\alpha$.

\begin{lemma}\label{lem:sections}
Let $A,Y\in\Lc$, $B,S\in\Zp[[z]]$, $r\in\Sig$ and $t\ge0$.
\begin{enumerate}[label=\textup{(\alph*)}]
\item $C(A\,\sigma(B))=C(A)\,B$.
\item $C(\zs^{-r}C^t(Y))=C^{t+1}(\zs^{-rp^t}Y)$.
\item If $t\ge1$ and $0\le v<p^t$, then $C^t(\zs^{-v}S)\in\Zp[[z]]$, and it is a polynomial if $S$ is.
\item $C_r$ maps $\Zp[[z]]$ to itself, and $C_w(S)=C^{k}(\zs^{-v(w)}S)$ for every word $w$ of length $k$.
\end{enumerate}
\end{lemma}

\begin{proof}
(a) We have $A\sigma(B)=\sum_{u,w}a_ub_wz^{u+pw}$, each coefficient being a finite sum because the exponents of $A$ are bounded below and those of $B$ are nonnegative. Moreover $u+pw\in p\Z^m$ if and only if $u\in p\Z^m$. Hence $C(A\sigma(B))=\sum_{u',w}a_{pu'}b_wz^{u'+w}=C(A)B$.

(b) The coefficient of $z^u$ on the left is $[z^{pu+r\one}]C^t(Y)=[z^{p^t(pu+r\one)}]Y=[z^{p^{t+1}u}](\zs^{-rp^t}Y)$, which is the coefficient of $z^u$ on the right.

(c) All exponents of $\zs^{-v}S$ are $\ge-v>-p^t$ in every coordinate. An exponent divisible by $p^t$ is therefore $\ge0$, and $C^t$ divides it by $p^t$. Finitely many monomials remain if $S$ is a polynomial.

(d) The first claim is (c) with $t=1$. The formula follows from (b) by induction on $k$, since $C_r(C^k(\zs^{-v(w)}S))=C^{k+1}(\zs^{-v(w)-rp^k}S)$ and $v(w)+rp^k=v(wr)$.
\end{proof}

\begin{lemma}\label{lem:diag}
For $S\in\Zp[[z]]$ and $r\in\Sig$ we have $\Diag(C_rS)=\Lambda_r\Diag(S)$. Consequently $\Diag(C_wS)=\Lambda_w\Diag(S)$ for every word $w$, and
\[
K_\alpha(\Diag S)\le\bigl|\{C_wS\bmod p^\alpha:\ w\in\Sig^*\}\bigr|.
\]
\end{lemma}

\begin{proof}
$[x^n]\Diag(C_rS)=[z^{n\one}]C(\zs^{-r}S)=[z^{pn\one}](\zs^{-r}S)=[z^{(pn+r)\one}]S=[x^{pn+r}]\Diag S$. Iterating gives the second claim. By \eqref{eq:Kalpha}, $K_\alpha(\Diag S)$ counts the series $\Lambda_w\Diag S\bmod p^\alpha=\Diag(C_wS\bmod p^\alpha)$, and each of them is a function of $C_wS\bmod p^\alpha$.
\end{proof}

\subsection{Normalisation}\label{subsec:norm}
Let $N,Q\in\Zp[z]$ with $c:=Q(0)\in\Zp^\times$, and put $R=N/Q\in\Zp[[z]]$ and $G=\Diag R$. Replacing $(N,Q)$ by $(c^{-1}N,c^{-1}Q)$ changes neither $R$ nor the supports, degrees and Newton polytopes. \emph{From now on we assume $Q(0)=1$.} We write $\Delta=\Delta(Q)$, which contains $0$, and $K=K(N)$. Both lie in the box $\prod_i[0,h_i]$.

\subsection{From sections to digit matrices}\label{subsec:matrix-primer}

What must an automaton remember after reading a word $w$?
Equation~\eqref{eq:Lambdaw} gives a direct answer: the remaining
input indexes the coefficients of the residual series $\Lambda_wG$.
By Lemma~\ref{lem:diag}, we can retain the multivariate section
$C_w(N/Q)$ instead. Its constant coefficient is the output.
We now describe what must be stored to compute the next section.

Start modulo $p$. Write $z^p=(z_1^p,\ldots,z_m^p)$. Since
$Q(z)^p\equiv Q(z^p)\pmod p$, Lemma~\ref{lem:sections}(a) gives
\[
 C_r(M/Q)\equiv C_r(MQ^{p-1})/Q\pmod p.
\]
Thus the denominator stays fixed. After the first digit, the new
numerator is supported on $J$ (Corollary~\ref{cor:mJ}).
Order the monomials $z^u$, $u\in J$, and write their coefficients
as a row vector $[M]$. Multiplication by $Q^{p-1}$ followed by
$C_r$ is linear. Its entry in source row $u$ and target column $b$
is the coefficient of $z^{pb+r\one-u}$ in $Q^{p-1}$.
An exponent outside the support gives entry zero. Multiplying $[M]$
on the right by this matrix and reducing modulo $p$ gives the next
numerator. Its constant coefficient is the output, since $Q(0)=1$.

At precision $p^\alpha$, the polynomials $\delta_t$ defined by
\eqref{eq:deltadef} supply the higher-precision terms.
In particular $\delta_1=Q^{p-1}$ and $p^{t-1}$ divides every
coefficient of $\delta_t$ (Lemma~\ref{lem:delta}). For
$1\le t\le\alpha$ and $0\le v<p^t$, form the matrix
\[
 (B_{t,v})_{u,b}=[z^{p^t b+v\one-u}]\delta_t,
 \qquad u,b\in J .
\]
Here $[z^a]$ means coefficient extraction. These are matrices over
$\Zp$ of the operators $T_{t,v}$ in Definition~\ref{def:T};
Lemma~\ref{lem:stable} proves that their values stay supported on $J$.

For the first $\alpha-1$ input digits the automaton stores the word
itself. Starting with the $\alpha$th digit, its state consists of the
last $\alpha-1$ digits $W$ and $\alpha$ coefficient rows
$H_0,\ldots,H_{\alpha-1}$, with $H_s$ stored modulo $p^{\alpha-s}$.
Definition~\ref{def:history} obtains these rows from the successive
numerators of \eqref{eq:mrec}, starting with $m_0=N$.
On appending $r$, let $v_t$ be the integer represented by the last
$t$ digits of $Wr$, read least significant first. The new first row is
\[
 H'_0=\sum_{t=1}^{\alpha} H_{t-1}B_{t,v_t}\pmod {p^\alpha}.
\]
Any lifts of the old rows give the same answer: the missing multiple
of $p^{\alpha-t+1}$ in row $H_{t-1}$ is multiplied by a matrix
divisible by $p^{t-1}$. Shift each remaining row one place,
$H'_s=H_{s-1}\bmod p^{\alpha-s}$, and keep the last $\alpha-1$
digits of $Wr$. Proposition~\ref{prop:transition} proves this update;
Definition~\ref{def:decode} recovers the section from the stored
data, and Corollary~\ref{cor:automaton} gives its output.

\subsection{The appendix example modulo two}\label{subsec:matrix-example}

Take $P=x^2y^2+(1+x+x^2)y+x^2$ from Appendix~\ref{app:zero}.
In diagonal coordinates its series is $F=\Diag(N/Q)$, where the numerator and denominator in
Section~\ref{sec:cor} reduce modulo $2$ to
\[
 N=y+xy^2+x^2y^3,\qquad Q=1+xy+x^2y+x^2y^2+x^2y^3 .
\]
Work over $\F_2$ and order the monomials after the first digit as
$(1,xy,xy^2,x^2y^2,x^2y^3)$. Applying $M\mapsto C_r(MQ)$ to
each basis monomial gives the following two matrices, with source
monomials indexing rows:
\[
 B_{1,0}=\begin{pmatrix}
 1&1&0&0&0\\0&1&0&0&0\\0&0&0&0&0\\0&1&0&1&0\\0&0&0&1&1
 \end{pmatrix},\qquad
 B_{1,1}=\begin{pmatrix}
 1&0&0&0&0\\1&1&0&0&0\\0&1&1&0&0\\0&1&0&0&0\\0&0&0&0&0
 \end{pmatrix}.
\]
For instance, of the monomials in $xy^2Q$, precisely $x^3y^3$ and
$x^3y^5$ have both exponents odd. Applying $C_1$ gives $xy+xy^2$,
the third row of $B_{1,1}$. The initial numerator $N$ includes $y$,
so compute the first digit directly: $C_0(NQ)$ and $C_1(NQ)$ have
rows $(0,1,1,1,1)$ and $(0,1,1,0,0)$, respectively. The empty word
has output zero. Every subsequent output is the first entry of its
row. For example, the digits of $2$ arrive as $0,1$.
The product $(0,1,1,1,1)B_{1,1}=(1,1,1,0,0)$ gives
$[x^2]F=1\pmod2$, which also follows by comparing the coefficients
of $x$ and $x^2$ in $P(x,F)=0$.
The two-by-two matrices in the proof of Lemma~\ref{lem:odd-history-coordinates}
describe the restriction to the numerators $u+vx$ in the central-row
coordinates of Section~\ref{subsec:odd-upper}. The five rows above also
include the other monomials needed just after the first digit.

\section{The \texorpdfstring{$\delta$}{delta}-decomposition and levels}\label{sec:delta}

\subsection{The \texorpdfstring{$\delta$}{delta}-polynomials}
Define $\delta_s\in\Zp[z]$ for $s\ge1$ recursively by
\begin{equation}\label{eq:deltadef}
Q^{p^s-1}=\sum_{t=1}^s\delta_t\,\sigma^t\bigl(Q^{p^{s-t}-1}\bigr),
\end{equation}
that is, $\delta_1=Q^{p-1}$ and $\delta_s=Q^{p^s-1}-\sum_{t=1}^{s-1}\delta_t\sigma^t(Q^{p^{s-t}-1})$.

We shall use two elementary facts, valid for every prime $p$ (including $p=2$) and for elements of $\Zp[z]$ or $\Zp[[z]]$:
\begin{enumerate}[label=(F\arabic*)]
\item\label{F1} $\sigma(A)\equiv A^p\pmod p$, and hence $\sigma^t(A)\equiv A^{p^t}\pmod p$ for $t\ge0$. Indeed $(\sum c_uz^u)^p\equiv\sum c_u^pz^{pu}$ and $c^p\equiv c\pmod p$ for $c\in\Zp$. Note that $\sigma(A)$ is again integral, so the congruence can be iterated.
\item\label{F2} If $A\equiv B\pmod{p^e}$ with $e\ge1$, then $A^p\equiv B^p\pmod{p^{e+1}}$. Indeed, writing $A=B+p^eX$, the term $\binom pjB^{p-j}p^{ej}X^j$ is divisible by $p^{e+1}$ for $j=1$ (as $p\mid\binom p1$) and for $j\ge2$ (as $ej\ge e+1$). Consequently $A\equiv B\pmod p$ implies $A^{p^e}\equiv B^{p^e}\pmod{p^{e+1}}$.
\end{enumerate}

\begin{lemma}[Vlasenko {\cite[Lemmas 3 and 4]{Vl18}}]\label{lem:delta}
For every $s\ge1$:
\begin{enumerate}[label=\textup{(\alph*)}]
\item $\supp\delta_s\subset(p^s-1)\Delta$;
\item $\delta_s\in p^{s-1}\Zp[z]$.
\end{enumerate}
\end{lemma}

\begin{proof}
(a) For $s=1$, $\delta_1=Q^{p-1}$. For $s\ge2$, $\supp Q^{p^s-1}\subset(p^s-1)\Delta$, and by induction
\[
\supp\bigl(\delta_t\sigma^t(Q^{p^{s-t}-1})\bigr)\subset(p^t-1)\Delta+p^t(p^{s-t}-1)\Delta=(p^s-1)\Delta,
\]
since $a\Delta+b\Delta=(a+b)\Delta$ for a convex set $\Delta$ and $a,b\ge0$.

(b) Induction on $s$; the case $s=1$ is trivial. Let $s\ge2$ and assume $\delta_t\in p^{t-1}\Zp[z]$ for $t<s$. For $1\le t<s$,
\[
\sigma^t(Q^{p^{s-t}-1})=\sigma^t(Q^{p^{s-1-t}-1})\cdot\sigma^t(Q^{p-1})^{p^{s-1-t}}.
\]
By \ref{F1}, $\sigma^t(Q^{p-1})\equiv Q^{(p-1)p^t}\pmod p$, and by \ref{F2}, $\sigma^t(Q^{p-1})^{p^{s-1-t}}\equiv Q^{(p-1)p^{s-1}}\pmod{p^{s-t}}$. Multiplying by $\delta_t\in p^{t-1}$ we get
\[
\delta_t\sigma^t(Q^{p^{s-t}-1})\equiv\delta_t\sigma^t(Q^{p^{s-1-t}-1})\,Q^{p^s-p^{s-1}}\pmod{p^{s-1}}.
\]
Summing over $t<s$ and using \eqref{eq:deltadef} for $s-1$,
\[
\delta_s\equiv Q^{p^s-1}-Q^{p^s-p^{s-1}}\sum_{t=1}^{s-1}\delta_t\sigma^t(Q^{p^{s-1-t}-1})=Q^{p^s-1}-Q^{p^s-p^{s-1}}Q^{p^{s-1}-1}=0\pmod{p^{s-1}}.\qedhere
\]
\end{proof}

\begin{lemma}\label{lem:inverse}
For every $s\ge1$,
\[
Q^{-1}\equiv\sum_{t=1}^s\delta_t\,\sigma^t(Q^{-1})\pmod{p^s\Zp[[z]]}.
\]
Consequently $Q^{-1}=\sum_{t\ge1}\delta_t\sigma^t(Q^{-1})$ in $\Zp[[z]]$, the series converging $p$-adically.
\end{lemma}

\begin{proof}
Since $Q(0)=1$, the series $Q$ and $\sigma^t(Q)$ are units of $\Zp[[z]]$. Put $u_t=\sigma^t(Q)\,Q^{-p^t}$. By \ref{F1}, $u_t\in1+p\Zp[[z]]$, so by \ref{F2}, $u_t^{p^e}\in1+p^{e+1}\Zp[[z]]$ for $e\ge0$. Since
\[
\sigma^t(Q^{p^{s-t}-1})=\sigma^t(Q)^{p^{s-t}}\sigma^t(Q)^{-1}=u_t^{p^{s-t}}Q^{p^s}\sigma^t(Q^{-1}),
\]
dividing \eqref{eq:deltadef} by $Q^{p^s}$ gives $Q^{-1}=\sum_{t=1}^s\delta_t\sigma^t(Q^{-1})\,u_t^{p^{s-t}}$. Therefore
\[
Q^{-1}-\sum_{t=1}^s\delta_t\sigma^t(Q^{-1})=\sum_{t=1}^s\delta_t\sigma^t(Q^{-1})\bigl(u_t^{p^{s-t}}-1\bigr)\in\sum_{t=1}^sp^{t-1}p^{s-t+1}\Zp[[z]]=p^s\Zp[[z]],
\]
using Lemma~\ref{lem:delta}(b). The terms with $t>s$ lie in $p^s\Zp[[z]]$, which gives the second statement.
\end{proof}

\subsection{Levels}

\begin{definition}\label{def:levels}
A \emph{level expansion} of $S\in\Zp[[z]]$ is a sequence $(M_j)_{j\ge0}$ of polynomials in $\Zp[z]$ with $M_j\to0$ $p$-adically and
\[
S=\sum_{j\ge0}M_j\,\sigma^j(Q^{-1}).
\]
The polynomial $M_j$ is the \emph{numerator at level $j$}.
\end{definition}

\begin{lemma}[one digit]\label{lem:onestep}
If $(M_j)_{j\ge0}$ is a level expansion of $S$ and $r\in\Sig$, then $\bigl(C_r(M_{j+1}+\delta_{j+1}M_0)\bigr)_{j\ge0}$ is a level expansion of $C_rS$.
\end{lemma}

\begin{proof}
By Lemma~\ref{lem:inverse}, $M_0Q^{-1}=\sum_{j\ge0}\delta_{j+1}M_0\,\sigma^{j+1}(Q^{-1})$. Hence $S=\sum_{j\ge0}A_j\,\sigma\bigl(\sigma^j(Q^{-1})\bigr)$ with $A_j=M_{j+1}+\delta_{j+1}M_0$, and $A_j\to0$. By continuity of $C_r$ and Lemma~\ref{lem:sections}(a) applied to $\zs^{-r}A_j\in\Lc$,
\[
C_rS=\sum_{j\ge0}C\bigl(\zs^{-r}A_j\,\sigma(\sigma^j(Q^{-1}))\bigr)=\sum_{j\ge0}C_r(A_j)\,\sigma^j(Q^{-1}).
\]
Each $C_r(A_j)$ is a polynomial by Lemma~\ref{lem:sections}(c), and $C_r(A_j)\to0$.
\end{proof}

\begin{definition}\label{def:T}
For $t\ge1$, $0\le v<p^t$ and $M\in\Zp[z]$ put
\[
T_{t,v}(M)=C^t\bigl(\zs^{-v}\,\delta_t\,M\bigr)\in\Zp[z].
\]
By Lemma~\ref{lem:delta}(b), $T_{t,v}(M)\in p^{t-1}\Zp[z]$. For a word $w=r_0\cdots r_{k-1}$ and $1\le t\le i\le k$ let
\[
v_{i,t}=v_{i,t}(w)=\sum_{q=0}^{t-1}r_{i-t+q}\,p^q\in[0,p^t),
\]
the integer whose digits (least significant first) are $r_{i-t},\dots,r_{i-1}$. The \emph{zero-level numerators} of $w$ are
\begin{equation}\label{eq:mrec}
m_0=N,\qquad m_i=\sum_{t=1}^iT_{t,v_{i,t}}(m_{i-t})\qquad(1\le i\le k).
\end{equation}
Note that $m_i$ depends only on $r_0,\dots,r_{i-1}$.
\end{definition}

\begin{proposition}[levels along a word]\label{prop:levels}
Let $w=r_0\cdots r_{k-1}$ be a word. Then $C_w(R)$ has the level expansion
\[
M^{(k)}_0=m_k,\qquad M^{(k)}_j=\sum_{t=1}^kC^t\bigl(\zs^{-v_{k,t}}\,\delta_{j+t}\,m_{k-t}\bigr)\quad(j\ge1).
\]
In particular $M^{(k)}_j\in p^j\Zp[z]$ for $j\ge1$.
\end{proposition}

\begin{proof}
For $k\ge1$ the two displayed formulas are the cases $j=0$ and $j\ge1$ of the single formula $M^{(k)}_j=\sum_{t=1}^kC^t(\zs^{-v_{k,t}}\delta_{j+t}m_{k-t})$, by \eqref{eq:mrec}. We argue by induction on $k$. For $k=0$, $R=N\cdot Q^{-1}$ has the level expansion $M_0=N$, $M_j=0$ for $j\ge1$. Suppose the claim holds for $w$ of length $k$, and let $r=r_k$. By Lemma~\ref{lem:onestep}, $C_{wr}(R)=C_r(C_wR)$ has the level expansion $M^{(k+1)}_j=C_r(M^{(k)}_{j+1}+\delta_{j+1}m_k)$. First,
\[
C_r(\delta_{j+1}m_k)=C\bigl(\zs^{-v_{k+1,1}}\delta_{j+1}m_{(k+1)-1}\bigr),
\]
since $v_{k+1,1}=r_k$. Second, by Lemma~\ref{lem:sections}(b),
\[
C_r\bigl(C^t(\zs^{-v_{k,t}}\delta_{j+1+t}m_{k-t})\bigr)=C^{t+1}\bigl(\zs^{-(v_{k,t}+r_kp^t)}\delta_{j+1+t}m_{k-t}\bigr),
\]
and $v_{k,t}+r_kp^t=v_{k+1,t+1}$. Renaming $t+1$ as $t$, we obtain $M^{(k+1)}_j=\sum_{t=1}^{k+1}C^t(\zs^{-v_{k+1,t}}\delta_{j+t}m_{k+1-t})$ for every $j\ge0$, as claimed. The divisibility follows from $\delta_{j+t}\in p^{j+t-1}\Zp[z]$, and it also shows $M^{(k)}_j\to0$.
\end{proof}

\section{A support set stable under all digit operators}\label{sec:J}

Let $J=J_p(N,Q)$ be as in Definition~\ref{def:J}, and let $L_J=\bigoplus_{u\in J}\Zp z^u$ be the polynomials supported on $J$.

\begin{lemma}\label{lem:Jbox}
\begin{enumerate}[label=\textup{(\alph*)}]
\item $J\subset\prod_{i=1}^m\{0,1,\dots,h_i\}$, so $|J|\le D=\prod_i(h_i+1)$.
\item If $\Delta\subset\Delta'$ and $K\subset K'$, then $J$ is contained in the set defined by $\Delta'$ and $K'$.
\end{enumerate}
\end{lemma}

\begin{proof}
(a) If $u+\eta\one=(1-\varepsilon)a+\varepsilon b$ with $a\in\Delta$, $b\in K$, then $0\le u_i\le u_i+\eta\le(1-\varepsilon)h_i+\varepsilon h_i=h_i$. (b) is clear.
\end{proof}

\begin{lemma}[stability]\label{lem:stable}
Let $t\ge1$ and $0\le v<p^t$. Then $T_{t,v}(L_J)\subset L_J$ and $T_{t,v}(N)\in L_J$.
\end{lemma}

\begin{proof}
Let $u\in J\cup\supp N$. If $u\in J$, choose $\eta\in[0,1)$, $\varepsilon\in(0,1/p]$, $a\in\Delta$ and $b\in K$ with $u+\eta\one=(1-\varepsilon)a+\varepsilon b$. If $u\in\supp N$, put $\eta=0$, $\varepsilon=1$, $b=u$ and let $a\in\Delta$ be arbitrary. The monomials of $\zs^{-v}\delta_tz^u$ are $z^{u+e-v\one}$ with $e\in\supp\delta_t\subset(p^t-1)\Delta$, so $e=(p^t-1)d$ with $d\in\Delta$ (Lemma~\ref{lem:delta}(a)). Such a monomial survives $C^t$ only if $u+e-v\one=p^tu'$ with $u'\in\Z^m$, and then $u'\in\Z^m_{\ge0}$ by Lemma~\ref{lem:sections}(c). Put $\eta'=(\eta+v)/p^t$. Since $0\le\eta<1$ and $0\le v\le p^t-1$, we have $\eta'\in[0,1)$. Moreover
\[
u'+\eta'\one=p^{-t}(u+\eta\one)+(1-p^{-t})d=(1-\varepsilon')a'+\varepsilon'b,
\]
where $\varepsilon'=\varepsilon p^{-t}\in(0,1/p]$ and
\[
a'=\frac{p^{-t}(1-\varepsilon)\,a+(1-p^{-t})\,d}{1-\varepsilon'}
\]
is a convex combination of $a,d\in\Delta$: its coefficients are nonnegative and add up to $(p^{-t}-\varepsilon p^{-t}+1-p^{-t})/(1-\varepsilon')=1$. Hence $a'\in\Delta$ and $u'\in J$. Since $T_{t,v}$ is $\Zp$-linear, both claims follow.
\end{proof}

\begin{corollary}\label{cor:mJ}
For every word $w$ and every $i\ge1$, the zero-level numerator $m_i$ lies in $L_J$.
\end{corollary}

\begin{proof}
By \eqref{eq:mrec}, $m_i=T_{i,v_{i,i}}(N)+\sum_{t=1}^{i-1}T_{t,v_{i,t}}(m_{i-t})$, and we apply Lemma~\ref{lem:stable} and induction on $i$.
\end{proof}

\begin{remark}\label{rem:lifts}
Only Lemma~\ref{lem:stable} is used about $J$ below. Hence Theorem~\ref{thm:main} holds with $|J|$ replaced by $|J'|$ for any finite $J'\subset\Z^m_{\ge0}$ satisfying the conclusion of Lemma~\ref{lem:stable}. Moreover $K_\alpha(G)$ depends only on $N\bmod p^\alpha$ and $Q\bmod p^\alpha$. One may therefore compute $J$ from any lifts of these reductions with unit constant term; for instance, lifts with coefficients in $\{0,1,\dots,p^\alpha-1\}$ have supports contained in those of $N$ and $Q$.
\end{remark}

\begin{example}\label{ex:multinomial}
Let $Q=1-z_1-\dots-z_m$ and $N=1$. Then $\Delta$ is the standard simplex and $K=\{0\}$. If $u+\eta\one\in(1-\varepsilon)\Delta$ with $u\in\Z_{\ge0}^m$, then $\sum_iu_i+m\eta\le1-\varepsilon<1$, so $u=0$. Thus $J=\{0\}$, while the box has $2^m$ points. Theorem~\ref{thm:main} gives
\[
K_\alpha\Bigl(\sum_{n\ge0}\frac{(mn)!}{n!^m}x^n\Bigr)\le p^{2\alpha-1}+\frac{p^\alpha-1}{p-1}
\]
for every $m\ge1$ and every prime $p$. Here all zero-level numerators are constants.
\end{example}

\section{Histories, windows and the count}\label{sec:count}

Fix $\alpha\ge1$. For $e\ge0$ write $(\Z/p^e)^J=L_J/p^eL_J$ for polynomials supported on $J$ with coefficients modulo $p^e$.

\subsection{Precision bookkeeping}

\begin{lemma}\label{lem:precision}
Let $e\ge0$, $t\ge1$, $j\ge0$, $0\le v<p^t$, and let $H\in(\Z/p^e)^J$ have lift $\widetilde H\in L_J$. Then
\[
C^t\bigl(\zs^{-v}\,\delta_{j+t}\,\widetilde H\bigr)\bmod p^{\,e+j+t-1}
\]
does not depend on the lift. In particular $T_{t,v}$ induces a well-defined map $(\Z/p^e)^J\to(\Z/p^{e+t-1})^J$, which we denote by the same symbol.
\end{lemma}

\begin{proof}
Two lifts differ by $p^eY$ with $Y\in L_J$, and $C^t(\zs^{-v}\delta_{j+t}p^eY)\in p^{e+j+t-1}\Zp[z]$ by Lemma~\ref{lem:delta}(b). The last statement is the case $j=0$, combined with Lemma~\ref{lem:stable}.
\end{proof}

\subsection{Digit values}
For a digit string $U=(u_1,\dots,u_L)\in\Sig^L$ and integers $s\ge0$, $t\ge1$ with $s+t\le L$, let
\begin{equation}\label{eq:digitvalue}
\dg Ust=\sum_{q=0}^{t-1}u_{L-s-t+1+q}\,p^q\in[0,p^t).
\end{equation}
This is the integer whose digits, least significant first, are the $t$ digits of $U$ that end $s$ places before the end of $U$. Directly from the definition, for $r\in\Sig$ and $Ur=(u_1,\dots,u_L,r)$:
\begin{enumerate}[label=(D\arabic*)]
\item\label{D1} $\dg{Ur}01=r$ and $\dg{Ur}0{t+1}=\dg U0t+rp^t$ for $1\le t\le L$;
\item\label{D2} $\dg{Ur}st=\dg U{s-1}t$ for $s\ge1$ and $s+t\le L+1$;
\item\label{D3} if $U'$ is obtained from $U$ by deleting initial digits, then $\dg{U'}st=\dg Ust$ whenever $s+t\le|U'|$.
\end{enumerate}
For a word $w=r_0\cdots r_{k-1}$, regarded as the string $(r_0,\dots,r_{k-1})$, we have
\begin{equation}\label{eq:vwindow}
v_{i,t}(w)=\dg w{k-i}t\qquad(1\le t\le i\le k).
\end{equation}

\subsection{Histories and windows}

\begin{definition}\label{def:history}
Let $w$ be a word of length $k\ge\alpha$. Its \emph{window} $W(w)\in\Sig^{\alpha-1}$ is the string of its last $\alpha-1$ digits. Its \emph{history} is
\[
H(w)=(H_0,\dots,H_{\alpha-1}),\qquad H_s=m_{k-s}\bmod p^{\alpha-s}\in(\Z/p^{\alpha-s})^J .
\]
This makes sense by Corollary~\ref{cor:mJ}, since $k-s\ge k-\alpha+1\ge1$.
\end{definition}

By \eqref{eq:vwindow} and \ref{D3},
\begin{equation}\label{eq:vk}
v_{k-s,t}(w)=\dg wst=\dg{W(w)}st\qquad\text{whenever }s+t\le\alpha-1 .
\end{equation}

\begin{definition}[decoding]\label{def:decode}
For $W\in\Sig^{\alpha-1}$ and $H=(H_0,\dots,H_{\alpha-1})\in\prod_{s=0}^{\alpha-1}(\Z/p^{\alpha-s})^J$ define $\Psi(H,W)\in(\Z/p^\alpha)[[z]]$ by
\begin{equation}\label{eq:Psi}
\Psi(H,W)=\widetilde H_0\,Q^{-1}+\sum_{j=1}^{\alpha-1}\ \sum_{t=1}^{\alpha-j}C^t\bigl(\zs^{-\dg W0t}\,\delta_{j+t}\,\widetilde H_t\bigr)\,\sigma^j(Q^{-1})\pmod{p^\alpha},
\end{equation}
where $\widetilde H_s\in L_J$ are lifts. The value does not depend on the lifts: the first term only involves $H_0$ modulo $p^\alpha$, and by Lemma~\ref{lem:precision} with $e=\alpha-t$ and $j\ge1$, the summand indexed by $(j,t)$ is determined modulo $p^{\alpha+j-1}$, hence modulo $p^\alpha$. Note that $t\le\alpha-j\le\alpha-1=|W|$.
\end{definition}

\begin{lemma}[decoding]\label{lem:decode}
For every word $w$ with $|w|\ge\alpha$,
\[
C_w(R)\equiv\Psi\bigl(H(w),W(w)\bigr)\pmod{p^\alpha}.
\]
\end{lemma}

\begin{proof}
Let $k=|w|$ and use the level expansion of Proposition~\ref{prop:levels}. The level-$0$ numerator is $m_k$, which reduces to $H_0$. For $j\ge\alpha$, $M^{(k)}_j\in p^j\Zp[z]\subset p^\alpha\Zp[z]$. For $1\le j\le\alpha-1$ consider the summand of $M^{(k)}_j$ indexed by $t$. If $t>\alpha-j$, it lies in $p^{j+t-1}\Zp[z]\subset p^\alpha\Zp[z]$. If $t\le\alpha-j$, then $t\le\alpha-1\le k-1$, and $v_{k,t}=\dg{W(w)}0t$ by \eqref{eq:vk} with $s=0$. Moreover $m_{k-t}$ is a lift of $H_t$, which is given modulo $p^{\alpha-t}$, so the summand agrees modulo $p^{\alpha+j-1}$ with the corresponding summand of \eqref{eq:Psi} by Lemma~\ref{lem:precision}.
\end{proof}

\subsection{Compatible histories}

\begin{definition}\label{def:compatible}
Let $W\in\Sig^{\alpha-1}$. A tuple $H\in\prod_{s=0}^{\alpha-1}(\Z/p^{\alpha-s})^J$ is \emph{compatible with $W$} if, for $0\le s\le\alpha-2$,
\begin{equation}\label{eq:compat}
H_s\equiv\sum_{t=1}^{\alpha-s-1}T_{t,\dg Wst}\bigl(H_{s+t}\bigr)\pmod{p^{\alpha-s-1}} .
\end{equation}
The right-hand side is well defined in $(\Z/p^{\alpha-s-1})^J$ by Lemma~\ref{lem:precision}, because $H_{s+t}$ is given modulo $p^{\alpha-s-t}$. Let $\Hc_W$ be the set of histories compatible with $W$. Since the maps $T_{t,v}$ are additive, $\Hc_W$ is a subgroup of $\prod_s(\Z/p^{\alpha-s})^J$.
\end{definition}

\begin{lemma}\label{lem:compatible}
For every word $w$ with $|w|\ge\alpha$ we have $H(w)\in\Hc_{W(w)}$.
\end{lemma}

\begin{proof}
Let $k=|w|$ and $0\le s\le\alpha-2$. By \eqref{eq:mrec}, $m_{k-s}=\sum_{t=1}^{k-s}T_{t,v_{k-s,t}}(m_{k-s-t})$. The summands with $t\ge\alpha-s$ lie in $p^{t-1}\Zp[z]\subset p^{\alpha-s-1}\Zp[z]$. For $t\le\alpha-s-1$ (note $t<k-s$) we have $v_{k-s,t}=\dg{W(w)}st$ by \eqref{eq:vk}. By Lemma~\ref{lem:precision}, $T_{t,v}(m_{k-s-t})\bmod p^{\alpha-s-1}$ depends only on $m_{k-s-t}\bmod p^{\alpha-s-t}$, which is $H_{s+t}$. Reducing modulo $p^{\alpha-s-1}$ gives \eqref{eq:compat}.
\end{proof}

\begin{lemma}[exact count]\label{lem:count}
For every $W\in\Sig^{\alpha-1}$ we have $|\Hc_W|=p^{|J|\alpha}$. More precisely, the map sending $H\in\Hc_W$ to the top $p$-adic digits of its $|J|\alpha$ coefficients is a bijection onto $\prod_{s=0}^{\alpha-1}(\Z/p)^J$.
\end{lemma}

\begin{proof}
For $0\le s\le\alpha-1$ let $\Hc^{(s)}_W$ be the set of tuples $(H_s,\dots,H_{\alpha-1})$ satisfying \eqref{eq:compat} for all $s'$ with $s\le s'\le\alpha-2$. Then $\Hc^{(\alpha-1)}_W=(\Z/p)^J$ has $p^{|J|}$ elements. Let $s\le\alpha-2$. A tuple in $\Hc^{(s)}_W$ consists of a tuple in $\Hc^{(s+1)}_W$ together with an element $H_s\in(\Z/p^{\alpha-s})^J$ whose reduction modulo $p^{\alpha-s-1}$ is prescribed by \eqref{eq:compat}. The reduction map $(\Z/p^{\alpha-s})^J\to(\Z/p^{\alpha-s-1})^J$ is surjective with fibres of size $p^{|J|}$, and an element of a fibre is determined by the top digits $\lfloor\overline H_s/p^{\alpha-s-1}\rfloor$ of its coefficients, where $\overline H_s$ has coefficients in $[0,p^{\alpha-s})$. Hence $|\Hc^{(s)}_W|=p^{|J|}|\Hc^{(s+1)}_W|$, and the top digits can be chosen freely at each level. Therefore $|\Hc_W|=|\Hc^{(0)}_W|=p^{|J|\alpha}$.
\end{proof}

\subsection{Proof of Theorem~\ref{thm:main}}\label{subsec:proof-main}
By Lemma~\ref{lem:diag}, $K_\alpha(G)\le|\{C_wR\bmod p^\alpha:\ w\in\Sig^*\}|$. Words of length less than $\alpha$ contribute at most $\sum_{k=0}^{\alpha-1}p^k=(p^\alpha-1)/(p-1)$ elements. For $|w|\ge\alpha$, Lemmas~\ref{lem:decode} and~\ref{lem:compatible} give $C_wR\bmod p^\alpha=\Psi(H(w),W(w))$ with $H(w)\in\Hc_{W(w)}$. By Lemma~\ref{lem:count} there are at most
\[
\sum_{W\in\Sig^{\alpha-1}}|\Hc_W|=p^{\alpha-1}\cdot p^{|J|\alpha}=p^{(|J|+1)\alpha-1}
\]
such values. Finally $|J|\le D$ by Lemma~\ref{lem:Jbox}. \qed

\subsection{The explicit automaton}

\begin{proposition}[transitions]\label{prop:transition}
Let $W\in\Sig^{\alpha-1}$, $H\in\prod_{s=0}^{\alpha-1}(\Z/p^{\alpha-s})^J$ and $r\in\Sig$. Let $W'$ be the string of the last $\alpha-1$ digits of $Wr$, and define $H'=(H'_0,\dots,H'_{\alpha-1})$ by
\begin{equation}\label{eq:transition}
H'_0=\sum_{t=1}^{\alpha}T_{t,\dg{Wr}0t}\bigl(H_{t-1}\bigr)\in(\Z/p^\alpha)^J,\qquad H'_s=H_{s-1}\bmod p^{\alpha-s}\quad(1\le s\le\alpha-1).
\end{equation}
Write $\tau_r(H,W)=(H',W')$. Then:
\begin{enumerate}[label=\textup{(\alph*)}]
\item $C_r\,\Psi(H,W)=\Psi(\tau_r(H,W))$ in $(\Z/p^\alpha)[[z]]$, for every $H$, compatible or not;
\item if $H\in\Hc_W$, then $H'\in\Hc_{W'}$;
\item for every word $w$ with $|w|\ge\alpha$, $\tau_r(H(w),W(w))=(H(wr),W(wr))$.
\end{enumerate}
\end{proposition}

The summand of $H'_0$ indexed by $t$ is well defined in $(\Z/p^\alpha)^J$ by Lemma~\ref{lem:precision}, since $H_{t-1}$ is given modulo $p^{\alpha-t+1}$.

\begin{proof}
(a) Choose lifts $\widetilde H_s\in L_J$ and write $\Psi(H,W)$ as the reduction of $\sum_{j\ge0}M_j\sigma^j(Q^{-1})$, where $M_0=\widetilde H_0$, $M_j=\sum_{t=1}^{\alpha-j}C^t(\zs^{-\dg W0t}\delta_{j+t}\widetilde H_t)$ for $1\le j\le\alpha-1$, and $M_j=0$ for $j\ge\alpha$. By Lemma~\ref{lem:onestep}, $C_r$ of this series has the level expansion $N_j=C_r(M_{j+1}+\delta_{j+1}\widetilde H_0)$. For $j\ge\alpha$, $N_j=C_r(\delta_{j+1}\widetilde H_0)\in p^j\Zp[z]$ vanishes modulo $p^\alpha$. For $0\le j\le\alpha-1$ we compute, using $\dg{Wr}01=r$, Lemma~\ref{lem:sections}(b) and \ref{D1},
\[
N_j=C^1\bigl(\zs^{-\dg{Wr}01}\delta_{j+1}\widetilde H_0\bigr)+\sum_{t=1}^{\alpha-j-1}C^{t+1}\bigl(\zs^{-\dg{Wr}0{t+1}}\delta_{j+1+t}\widetilde H_t\bigr)
=\sum_{t=1}^{\alpha-j}C^t\bigl(\zs^{-\dg{Wr}0t}\delta_{j+t}\widetilde H_{t-1}\bigr).
\]
For $j=0$ this is $\sum_{t=1}^\alpha T_{t,\dg{Wr}0t}(\widetilde H_{t-1})$, a lift of $H'_0$. Now let $1\le j\le\alpha-1$. The level-$j$ numerator of $\Psi(H',W')$ is $\sum_{t=1}^{\alpha-j}C^t(\zs^{-\dg{W'}0t}\delta_{j+t}\widetilde H'_t)$. Here $\dg{W'}0t=\dg{Wr}0t$ by \ref{D3}, as $t\le\alpha-1$. Also $\widetilde H'_t\equiv\widetilde H_{t-1}\pmod{p^{\alpha-t}}$, so by Lemma~\ref{lem:precision} the two summands agree modulo $p^{\alpha+j-1}$. Hence $N_j\equiv M'_j\pmod{p^\alpha}$ for all $j$, and $C_r\Psi(H,W)=\Psi(H',W')$.

(b) First let $s=0$. Reduce \eqref{eq:transition} modulo $p^{\alpha-1}$. The term $t=\alpha$ lies in $p^{\alpha-1}$ and disappears. For $t\le\alpha-1$, the term $T_{t,\dg{Wr}0t}(H_{t-1})\bmod p^{\alpha-1}$ depends only on $H_{t-1}\bmod p^{\alpha-t}=H'_t$ by Lemma~\ref{lem:precision}, and $\dg{Wr}0t=\dg{W'}0t$ by \ref{D3}. This gives \eqref{eq:compat} for $(H',W')$ at $s=0$.

Next let $1\le s\le\alpha-2$. Relation \eqref{eq:compat} for $(H,W)$ at index $s-1$ reads
\[
H_{s-1}\equiv\sum_{t=1}^{\alpha-s}T_{t,\dg W{s-1}t}(H_{s-1+t})\pmod{p^{\alpha-s}} .
\]
Reduce it modulo $p^{\alpha-s-1}$. The term $t=\alpha-s$ disappears. For $t\le\alpha-s-1$, the term depends only on $H_{s+t-1}\bmod p^{\alpha-s-t}=H'_{s+t}$. Moreover $\dg W{s-1}t=\dg{Wr}st=\dg{W'}st$ by \ref{D2} and \ref{D3}. Since $H'_s=H_{s-1}\bmod p^{\alpha-s}$, we obtain \eqref{eq:compat} for $(H',W')$ at index $s$.

(c) Let $k=|w|$. By \eqref{eq:mrec}, $m_{k+1}=\sum_{t=1}^{k+1}T_{t,v_{k+1,t}(wr)}(m_{k+1-t})$. Modulo $p^\alpha$ the terms with $t\ge\alpha+1$ vanish. For $t\le\alpha$, $v_{k+1,t}(wr)=\dg{wr}0t=\dg{Wr}0t$ by \eqref{eq:vwindow} and \ref{D3}, because $Wr$ is the suffix of $wr$ of length $\alpha$. The term depends only on $m_{k+1-t}\bmod p^{\alpha-t+1}=H_{t-1}$. Hence $H(wr)_0=H'_0$. For $s\ge1$, $H(wr)_s=m_{k+1-s}\bmod p^{\alpha-s}=H_{s-1}\bmod p^{\alpha-s}$. Finally $W(wr)=W'$.
\end{proof}

\begin{corollary}[explicit automaton]\label{cor:automaton}
Consider the automaton with input alphabet $\Sig$ whose states are the words of length less than $\alpha$ together with the pairs $(H,W)$ with $W\in\Sig^{\alpha-1}$ and $H\in\Hc_W$. Its initial state is the empty word. Its transitions are $w\xrightarrow{r}wr$ for $|w|<\alpha-1$, $w\xrightarrow{r}(H(wr),W(wr))$ for $|w|=\alpha-1$, and $(H,W)\xrightarrow{r}\tau_r(H,W)$. Its outputs are $g(v(w))\bmod p^\alpha$ at a word $w$, and
\[
\omega(H,W)=[z^0]\Psi(H,W)=H_0(0)+\sum_{j=1}^{\alpha-1}\sum_{t=1}^{\alpha-j}[z^{\dg W0t\one}]\bigl(\delta_{j+t}\widetilde H_t\bigr)\bmod p^\alpha
\]
at $(H,W)$. This automaton generates $(g(n)\bmod p^\alpha)_{n\ge0}$, reading digits least significant first, and it has at most $p^{(|J|+1)\alpha-1}+(p^\alpha-1)/(p-1)$ states.
\end{corollary}

\begin{proof}
By Proposition~\ref{prop:transition}(b),(c) the transitions are well defined and the state reached after reading $w$ with $|w|\ge\alpha$ is $(H(w),W(w))$. Its output is $[z^0]\Psi(H(w),W(w))=[z^0]C_wR=[x^0]\Lambda_wG=g(v(w))$ modulo $p^\alpha$, by Lemmas~\ref{lem:decode} and~\ref{lem:diag} and \eqref{eq:Lambdaw}. The displayed formula for $\omega$ holds because $[z^0]\sigma^j(Q^{-1})=1$ and $[z^0]C^t(\zs^{-v}Y)=[z^{v\one}]Y$. The bound on the number of states is Lemma~\ref{lem:count}.
\end{proof}

\begin{remark}\label{rem:automaton}
(i) The decoding map $\Psi$ need not be injective, and different pairs may represent the same kernel element; the minimal automaton is obtained from this one by restricting to the reachable states and applying Moore minimisation. (ii) The transitions are given by the reductions modulo $p^\alpha$ of the $|J|\times|J|$ matrices of $T_{t,v}$ on $L_J$ for $1\le t\le\alpha$ and $0\le v<p^t$; there are about $p^{\alpha+1}/(p-1)$ of them. Computing them requires the polynomials $\delta_t$, whose degrees are $O(p^t)$ for fixed $Q$. Afterwards each state is stored as $\alpha$ vectors of length $|J|$ together with $\alpha-1$ digits, and each transition costs $O(\alpha|J|^2)$ operations modulo $p^\alpha$ (Section~\ref{sec:computations}). (iii) For $\alpha=1$, after the first digit, the state is a single $H_0\in(\Z/p)^J$ with transition $H_0\mapsto T_{1,r}(H_0)=C(\zs^{-r}Q^{p-1}H_0)$. This is the classical construction modulo $p$; Theorem~\ref{thm:main} extends it to higher precision at the price of only a linear growth of the exponent.
\end{remark}

The numerator update in Remark~\ref{rem:automaton}(iii) at precision one is the specialization of the Cartier transition in the proof of \cite[Theorem~2.1 and Algorithm~1]{RY15}. For partition refinement of states with outputs, see the proof of \cite[Theorem~6]{Mo56}.

\begin{remark}[words of length $\alpha-1$]\label{rem:short-histories}
Suppose that the initial numerator $N$ lies in $L_J$. Then Definition~\ref{def:history} and Lemmas~\ref{lem:decode} and~\ref{lem:compatible} also hold for words $w$ of length $k=\alpha-1$, with the same proofs. The oldest entry of the history is then $H_{\alpha-1}=m_0=N$, which is supported on $J$, so Lemma~\ref{lem:precision} applies to it. In \eqref{eq:compat} the summand with $t=\alpha-1-s$ is the direct term $T_{k-s,v_{k-s,k-s}}(m_0)$ of \eqref{eq:mrec}, and $v_{k-s,k-s}=\dg{W(w)}st$ because $s+t=\alpha-1$. In the decoding formula every index satisfies $t\le\alpha-j\le k$. Section~\ref{sec:odd-exponent} uses this with $N$ replaced by an element of $W_0$.
\end{remark}

\section{Furstenberg series and a quadratic family}\label{sec:cor}

\subsection{Proof of Corollary~\ref{cor:furst}}\label{subsec:proof-furst}
Write $P=\sum_{i,j}c_{ij}x^iy^j$ with $c_{00}=0$ and $c_{01}=\partial_yP(0,0)\in\Zp^\times$. Then
\[
P(xy,y)=\sum c_{ij}x^iy^{i+j},\qquad y\,\partial_yP(xy,y)=\sum j\,c_{ij}x^iy^{i+j}.
\]
Every monomial with $c_{ij}\ne0$ has $i+j\ge1$, so $Q_F=P(xy,y)/y=\sum c_{ij}x^iy^{i+j-1}$ and $N_F=y\partial_yP(xy,y)$ are polynomials. The only monomial of $Q_F$ with exponent $(0,0)$ comes from $(i,j)=(0,1)$, so $Q_F(0,0)=c_{01}$ is a unit. By \cite[Lemma~6.3]{DL87} (see \cite[Theorem~4]{RY26}), $F=\Diag(N_F/Q_F)$. Every exponent $(u_1,u_2)$ of $Q_F$ has $0\le u_1\le h$ and $-1\le u_2-u_1\le d-1$, and every exponent of $N_F$ has $0\le u_1\le h$ and $1\le u_2-u_1\le d$. Hence every point of $(1-\varepsilon)\Delta(Q_F)+\varepsilon K(N_F)$ has first coordinate in $[0,h]$ and coordinate difference in $[-1+2\varepsilon,d-1+\varepsilon]$. If $u\in J$, then $u+\eta\one$ is such a point for some $\eta\in[0,1)$ and $\varepsilon\in(0,1/p]$, and adding $\eta\one$ does not change $u_2-u_1$. Since $0<\varepsilon\le1/2$, the integer $u_2-u_1$ lies in $\{0,\dots,d-1\}$, and $u_1\in\{0,\dots,h\}$ by Lemma~\ref{lem:Jbox}(a). Thus $|J|\le(h+1)d$, and Theorem~\ref{thm:main} gives the bound. In the central-row coordinates of Section~\ref{sec:odd-exponent} and of \cite{RSY25}, obtained by the shear $x^iy^j\mapsto x^iy^{j-i}$, this says that $J$ lies in the box $0\le i\le h$, $0\le j\le d-1$, which contains every state of the automaton of \cite{RSY25} except the initial one \cite[Proposition~9]{RSY25}. \qed

\subsection{The quadratic family}\label{subsec:proof-family}
Let $a,b,c,d\in\Z_2$ and $P=bx^2y^2+(1+ax+dx^2)y+cx^2$, so that $\partial_yP(0,0)=1$. Then
\[
Q_F=1+a\,xy+c\,x^2y+d\,x^2y^2+b\,x^2y^3,\qquad N_F=y+a\,xy^2+d\,x^2y^3+2b\,x^2y^4 .
\]

\begin{proof}[Proof of Corollary~\ref{cor:family}]
By Lemma~\ref{lem:Jbox}(b) we may enlarge the polytopes to
\[
\Delta'=\conv\{(0,0),(2,1),(2,3)\}=\{0\le X\le2,\ X/2\le Y\le3X/2\},
\]
\[
K'=\conv\{(0,1),(2,3),(2,4)\}=\{0\le X\le2,\ 1+X\le Y\le1+3X/2\},
\]
which contain $\supp Q_F$ and $\supp N_F$ for all parameter values. Let $0<\varepsilon\le1/2$ and $(X,Y)=(1-\varepsilon)(X_1,Y_1)+\varepsilon(X_2,Y_2)$ with $(X_1,Y_1)\in\Delta'$ and $(X_2,Y_2)\in K'$. Then $0\le X\le2$, and
\[
Y\ge(1-\varepsilon)\tfrac{X_1}2+\varepsilon(1+X_2)\ge\tfrac X2+\varepsilon,\qquad Y\le(1-\varepsilon)\tfrac{3X_1}2+\varepsilon\bigl(1+\tfrac{3X_2}2\bigr)=\tfrac{3X}2+\varepsilon .
\]
If $X=2$, then $X_1=X_2=2$ and $Y\ge(1-\varepsilon)+3\varepsilon=1+2\varepsilon$. Now let $u=(u_1,u_2)\in\Z^2_{\ge0}$ and $\eta\in[0,1)$ with $(u_1+\eta,u_2+\eta)$ of this form. Then:
\begin{itemize}
\item $u_1+\eta\le2$, so $u_1\le2$, and $\eta=0$ if $u_1=2$;
\item if $u_1=0$: $u_2+\eta\le\tfrac32\eta+\varepsilon$ gives $u_2\le\eta/2+\varepsilon<1$, so $u=(0,0)$;
\item if $u_1=1$: $u_2+\eta\ge\frac{1+\eta}2+\varepsilon$ gives $u_2\ge\frac{1-\eta}2+\varepsilon>0$, and $u_2+\eta\le\frac{3(1+\eta)}2+\varepsilon$ gives $u_2\le\frac32+\frac\eta2+\varepsilon<\frac52$; so $u_2\in\{1,2\}$;
\item if $u_1=2$: $1+2\varepsilon\le u_2\le3+\varepsilon$, so $u_2\in\{2,3\}$.
\end{itemize}
Hence $J_2(N_F,Q_F)\subset\{(0,0),(1,1),(1,2),(2,2),(2,3)\}$, so $|J|\le5$, and Theorem~\ref{thm:main} gives $K_\alpha(F)\le2^{6\alpha-1}+2^\alpha-1$. The lower bound for odd parameters is Theorem~\ref{thm:zero}(b) in the appendix.
\end{proof}

\begin{remark}\label{rem:numbers}
For odd $a,b,c,d$, as in Section~\ref{sec:odd-exponent}, the defining polynomial has bidegree $(2,2)$, so we take $h=d=2$ in \cite[Theorem~2]{RY26} and in Corollary~\ref{cor:furst}. The main term of the explicit bound of \cite[Theorem~2]{RY26}, which also gives the exponent of the asymptotic bound of \cite[Theorem~1]{RY26}, is $p^{N}$ with $N=\frac16\alpha(\alpha+1)(7\alpha+5)$; this equals $4,19,52,110,200$ for $\alpha=1,\dots,5$. The exponent of Corollary~\ref{cor:family} is $6\alpha-1$, that is $5,11,17,23,29$. The bound of Corollary~\ref{cor:furst} has exponent $7\alpha-1$, that is $6,13,20,27,34$, and becomes smaller than the cubic exponent from $\alpha=2$ on. This illustrates that all comparisons are asymptotic in $\alpha$ for fixed input.
\end{remark}

\begin{remark}
In the variables of \cite{RY26}, which use the shear $x\mapsto xy^{-1}$ and the central row $\CT_y$ instead of the diagonal, the five-point set containing $J$ above becomes $\{(0,0),(1,0),(1,1),(2,0),(2,1)\}$.
\end{remark}

\begin{remark}\label{rem:J-equality}
For odd $a,b,c,d$ the inclusion in the proof of Corollary~\ref{cor:family} is an equality. Here $\Delta(Q_F)$ is the triangle with vertices $(0,0),(2,1),(2,3)$ and $K(N_F)$ is the triangle with vertices $(0,1),(2,3),(2,4)$. With $(\varepsilon,\eta)=(\tfrac14,\tfrac12),(\tfrac14,0),(\tfrac12,0),(\tfrac14,0),(\tfrac12,0)$ respectively,
\[
(\tfrac12,\tfrac12)=\tfrac34(\tfrac23,\tfrac13)+\tfrac14(0,1),\quad
(1,1)=\tfrac34(1,\tfrac23)+\tfrac14(1,2),\quad
(1,2)=\tfrac12(2,3)+\tfrac12(0,1),
\]
\[
(2,2)=\tfrac34(2,\tfrac53)+\tfrac14(2,3),\quad
(2,3)=\tfrac12(2,3)+\tfrac12(2,3),
\]
where in each sum the first point lies in $\Delta(Q_F)$ and the second in $K(N_F)$. Thus all five points lie in $J$.
\end{remark}

\section{The exact exponent for the odd quadratic family}\label{sec:odd-exponent}

The lower bound in Corollary~\ref{cor:family} can be strengthened to a matching exponential rate. Throughout this section the parameters $a,b,c,d$ are odd elements of $\Z_2$, and $F(0)=0$ is the root of
\[
 bx^2F^2+(1+ax+dx^2)F+cx^2=0.
\]

\begin{theorem}\label{thm:odd-exponent}
For every $\alpha\ge1$,
\begin{equation}\label{eq:odd-exponent-bounds}
 \frac{9\,8^\alpha}{64\alpha-32}\le K_\alpha(F)
 \le 2^{3\alpha-1}
   +3(2^{\alpha+1}+\alpha+1)(2^\alpha-1).
\end{equation}
Consequently
\[
 \lim_{\alpha\to\infty}\frac{\log_2K_\alpha(F)}{\alpha}=3.
\]
\end{theorem}

The proof uses the histories of Section~\ref{sec:count}, first for numerators of the form $u+vx$ with $u,v\in\Z_2$ and then in normalized series coordinates. The normalized full-window frames are independent of the odd parameters; the entire kernel cardinality need not be. For example, $K_2=53$ for $(a,b,c,d)=(1,1,1,1)$ (Table~\ref{tab:main}), while $K_2=61$ for each of the four families in Table~\ref{tab:other}, whose values of $K_3$ are $264$, $266$ and $268$.

\subsection{Numerators of the form \texorpdfstring{$u+vx$}{u+vx} and three exceptional rays}
\label{subsec:odd-upper}

We use the central-row coordinates of Appendix~\ref{app:zero}:
\[
 P_1=1+ax+dx^2,\qquad
 f=P/y=1+ax+x^2(by+d+cy^{-1}),\qquad
 N=P_1y+2bx^2y^2,\qquad R=N/f.
\]
Thus $F=\CT_yR$. Write $\mathcal C_r=\Lambda_{r,0}$,
$\mathcal C=\mathcal C_0$ and $\mathcal C_w$ for the corresponding
word operator; then $\CT_y\mathcal C_wR=\Lambda_wF$.
Under the shear $x^iy^j\mapsto x^iy^{i+j}$, these operators and
$\delta_t(f)$ become the operators and $\delta$-polynomials of
Sections~\ref{sec:delta}--\ref{sec:count}. In particular their history,
compatibility and transition formulas hold with
\[
 T_{t,v}(M)=\mathcal C^t\bigl(x^{-v}\delta_t(f)M\bigr),\qquad
 J_{\rm sh}=\{(0,0),(1,0),(1,1),(2,0),(2,1)\}.
\]
The support of $\delta_t(f)$ is contained in $(2^t-1)\Delta$, where
\[
 \Delta=\{(i,j):i\le2,\ 2j\le i,\ -2j\le i\}.
\]
Put
\[
 W_0=\langle1,x\rangle_{\Z_2},\qquad
 V_0=\langle x^2,x^2y\rangle_{\Z_2},\qquad
 V_1=\langle xy\rangle_{\Z_2}.
\]
A history supported on $W_0$ has every entry in this rank-two module at its respective precision. We call such a history \emph{ordinary} in this section. An ordinary window fibre consists of all compatible histories supported on $W_0$ for a specified window.

\begin{lemma}\label{lem:odd-exceptional}
Every $T_{t,v}$ preserves $W_0$. The zero-level numerator $m_k(w)$
belongs to $W_0$ unless
\[
 w\in\mathcal E:=\{\emptyset\}\cup
 \{0^k,1^k,01^{k-1}:k\ge1\}.
\]
The language $\mathcal E$ is closed under taking prefixes.
\end{lemma}

\begin{proof}
A source monomial $x^iy^j$ can contribute to the target
$x^Iy^J$ only through the exponent
$(2^tI+v-i,2^tJ-j)$ of $\delta_t(f)$. To reach either target $x^2$ or
$x^2y$, its first coordinate must be at most $2^{t+1}-2$; hence
$i\ge v+2$. For a source in $J_{\rm sh}$ this is possible only if
$v=0$ and the source belongs to $V_0$. To reach $xy$, the upper
facet of $(2^t-1)\Delta$ gives
\[
 2(2^t-j)\le 2^t+v-i,\qquad v\ge 2^t+i-2j.
\]
For the sources $1,x,xy,x^2,x^2y$ the thresholds are respectively
$2^t,2^t+1,2^t-1,2^t+2,2^t$. Thus only $v=2^t-1$ and a source in $V_1$ can
contribute. Since $L_{J_{\rm sh}}=W_0\oplus V_0\oplus V_1$ is
stable, these two observations prove the assertion about $W_0$ and
show that on the quotient modulo $W_0$ only the $V_0$ block survives
when $v=0$, only the $V_1$ block survives when $v=2^t-1$, and all
other blocks vanish.

For the four possible source monomials of $N$, with exponents
$(0,1),(1,1),(2,1),(2,2)$, the first-coordinate test permits $V_0$
output only when $v=0$. The four thresholds for $V_1$ output are
$2^t-2,2^t-1,2^t,2^t-2$, so only $v=2^t-2$ or $v=2^t-1$ can contribute.
Project the recurrence~\eqref{eq:mrec} to the two summands of
$L_{J_{\rm sh}}/W_0$. Inductively, a nonzero $V_0$ component requires
an all-zero suffix and an all-zero preceding prefix, or the direct
source term with $v=0$; thus its word is $0^k$. A nonzero $V_1$
component requires an all-one suffix and a preceding prefix of the
form $1^j$ or $01^{j-1}$, or the direct source term with
$v=2^k-1$ or $2^k-2$. These are precisely $1^k$ and $01^{k-1}$.
This proves the claim, and prefix closure follows from the displayed
list of words.
\end{proof}

\begin{lemma}\label{lem:odd-rays}
Modulo $2^\alpha$, each of the three rays
$\mathcal C_{0^k}R$ ($k\ge0$), $\mathcal C_{1^k}R$ ($k\ge1$) and
$\mathcal C_{01^{k-1}}R$ ($k\ge1$) takes at most
$L_\alpha=2^{\alpha+1}+\alpha+1$ distinct values. Hence the states $\mathcal C_uR$ with $u\in\mathcal E$ take at most $3L_\alpha$ values.
\end{lemma}

\begin{proof}
We extend the annihilation argument of Appendix~\ref{app:zero} to
$M_j=x^jN$, $j=0,1,2$. For the all-zero word the support induction
in Lemma~\ref{lem:stable} gives
\[
 \supp m_k(M_j)\subset
 (1-2^{-k})\Delta+2^{-k}\operatorname{conv}(\supp M_j).
\]
On the second polytope the maxima of $X$ and $2Y-X$ are $j+2$ and
$2-j$, and the minimum of $X+2Y$ is $j+2>0$. Consequently points
in the displayed interpolation satisfy
\[
 X\le2+j2^{-k},\qquad 2Y-X\le(2-j)2^{-k},\qquad
 X+2Y\ge(j+2)2^{-k}>0.
\]
The remaining facet $-2Y\le X$ already holds. Taking integer points
shows
\[
 m_k(M_j)\in L_{J_0},\qquad
 J_0=\{(1,0),(2,0),(2,1)\},\qquad k\ge h_j,
 \quad (h_0,h_1,h_2)=(2,1,2).
\]

More generally, suppose an integral numerator $M$ has
$m_k(M)\in L_{J_0}$ for $k\ge h$, where $h\ge1$. Let
$\Theta_s$ and $\chi(t)=\sum_{r=0}^3\chi_rt^r=(t^2-1)(t-\mu)$,
$\mu$ odd, be as in Appendix~\ref{app:zero}. We claim that
\begin{equation}\label{eq:odd-extended-annihilation}
 \chi(\mathcal C)\mathcal C^k(M/f)
       \in2^{k-h+1}\widehat A_+\qquad(k\ge h-1).
\end{equation}
The notation $\widehat A_+$ denotes the series ring supported on
$i\ge0$, $i+j\ge0$, as defined in that appendix; all series here
belong to it. For $k\ge h$, split the zero-level recurrence as
\[
 m_k=\sum_{t=1}^{k-h}E_t(m_{k-t})+s_k,\qquad
 s_k=\sum_{t=k-h+1}^{k}E_t(m_{k-t}),
 \quad E_t(M)=\mathcal C^t(\delta_t(f)M).
\]
The first sum and $m_k$ lie in $L_{J_0}$, so $s_k$ does too, and
$s_k\in2^{k-h}L_{J_0}$. The right convolution
identity~\eqref{eq:rightconv} gives by induction
\[
 [m_k]=\sum_{l=h}^{k}[s_l]\Theta_{k-l}.
\]
In the residual $\sum_r\chi_rm_{k+r}$ the terms with $l\le k$
are divisible by
$2^{l-h}2^{k-l+1}=2^{k-h+1}$ by Lemma~\ref{lem:ratio}; terms
with $l>k$ are divisible by $2^{l-h}$, which is at least the same
power. This proves the required divisibility at level zero.
At level $j\ge1$, use~\eqref{eq:zerolevels}. For $t\le k-h$,
the inner residual at time $k-t$ is divisible by $2^{k-t-h+1}$;
multiplication by $\delta_{j+t}(f)$ supplies $2^{j+t-1}$, giving
$2^{k-h+j}$. For $t\ge k-h+1$, the factor
$\delta_{j+t}(f)$ alone is divisible by $2^{k-h+j}$. These
powers are at least $2^{k-h+1}$. The convergent level expansion
preserves the divisibility. The case $k=h-1$ has exponent zero
and is automatic. This proves~\eqref{eq:odd-extended-annihilation}.

At $k=\alpha+h-1$, the state $S=\mathcal C^k(M/f)$ is killed
by $\chi(\mathcal C)$ modulo $2^\alpha$. Put
$g_m(t)=\sum_{i=0}^{m-1}t^{2i}$. The polynomial
\[
 t^{2m}-1-g_m(\mu)(t^2-1)
   =(t^2-1)(g_m(t)-g_m(\mu))
\]
is divisible by $\chi(t)$. For $m=2^\alpha$, the factorisation
$g_m(\mu)=\prod_{i=0}^{\alpha-1}(1+\mu^{2^{i+1}})$ is divisible
by $2^\alpha$, since every factor is even. Thus
$(\mathcal C^{2^{\alpha+1}}-1)S=0$ modulo $2^\alpha$.
The orbit of $M/f$ therefore has preperiod at most
$\alpha+h-1$ and period dividing $2^{\alpha+1}$.

For a word of length $k$ and digit value $v$, its sheared operator
is $\mathcal C^k(x^{-v}\,\cdot\,)$. Coefficient extraction gives
\[
 \mathcal C_{0^k}R=\mathcal C^kR,\qquad
 \mathcal C_{1^k}R=x^{-1}\mathcal C^k(xR),\qquad
 \mathcal C_{01^{k-1}}R=x^{-1}\mathcal C^k(x^2R)\quad(k\ge1).
\]
The division by $x$ in the last two identities preserves
congruences; the selected numerators are divisible by $x$.
Apply the orbit estimate with $h=h_j\le2$. Each ray has at most
$2^{\alpha+1}+\alpha+1$ values, as required.
\end{proof}

\begin{proposition}\label{prop:odd-upper}
For every $\alpha\ge1$,
\[
 K_\alpha(F)\le2^{3\alpha-1}
    +3(2^{\alpha+1}+\alpha+1)(2^\alpha-1).
\]
\end{proposition}

\begin{proof}
Let a word first leave $\mathcal E$ at its prefix of length $e$.
Every later zero-level numerator belongs to $W_0$, by
Lemma~\ref{lem:odd-exceptional}. If its length $k$ satisfies
$k-e\ge\alpha-1$, its entire precision-$\alpha$ history is
supported on $W_0$. Since all $T_{t,v}$ preserve $W_0$, the proof
of Lemma~\ref{lem:count} gives exactly $2^{2\alpha}$ compatible
$W_0$-histories per window. There are $2^{\alpha-1}$ windows,
so these states have at most $2^{3\alpha-1}$ distinct decoded
values.

Otherwise the word is $urv$, where $u\in\mathcal E$, $r$ is an
exit digit and $|v|\le\alpha-2$. By Lemma~\ref{lem:odd-rays}
there are at most $3L_\alpha$ states after $u$, at most two
choices of $r$, and $2^{\alpha-1}-1$ suffixes of length at most
$\alpha-2$. Including the exceptional states themselves gives
at most
\[
 3L_\alpha\bigl(1+2(2^{\alpha-1}-1)\bigr)
   =3L_\alpha(2^\alpha-1)
\]
additional states. For $\alpha=1$ the suffix sum is empty; the
same formula remains valid. Every word is covered by one of
these cases, including words too short to fill a history.
Taking the central row cannot increase the count.
\end{proof}

\subsection{Window planes and normalized frames}\label{subsec:odd-frames}
Set $R_\alpha=\Z/2^\alpha\Z$, $P_\theta=1+ax+dx^2$, and
\[
 B_\theta=\CT_y\frac1f=(P_\theta^2-4bcx^4)^{-1/2},\qquad B_\theta(0)=1,
 \qquad W_0=\langle1,x\rangle_{\Z_2}.
\]
Here the constant-term identity follows by expanding $1/f$ in powers of $x^2(by+cy^{-1})/P_\theta$: only the even powers contribute, and their constant terms are $\binom{2j}{j}(bc)^j$. Thus
\[
 B_\theta=P_\theta^{-1}\sum_{j\ge0}\binom{2j}{j}
      (bcx^4/P_\theta^2)^j\in\Z_2[[x]].
\]
For a binary word $w=r_0\cdots r_{k-1}$ read in that order, put $v(w)=\sum_{i=0}^{k-1}r_i2^i$ and $\Lambda_w=\Lambda_{r_{k-1}}\allowbreak\cdots\Lambda_{r_0}=\Lambda_{2^k,v(w)}$, where $\Lambda_{2^k,v}g=\sum_ng(2^kn+v)x^n$. Write
\[
 \widehat E_w=(\Lambda_wB_\theta,\Lambda_w(xB_\theta)),\qquad
 L(h)=\binom{[x^0]h}{[x^1]h},\qquad M_w=L(\widehat E_w).
\]
The map $L$ is applied to each column of a row of series. If $B_\theta=\sum b_nx^n$, its row of index $n$ is $(b_{2^kn+v},b_{2^kn+v-1})$, with $b_j=0$ for $j<0$.
Modulo $2$, $B_\theta=1/(1+x+x^2)$, and sections on the numerator space $\langle1,x\rangle$ have row matrices
\[
 \begin{pmatrix}1&1\\0&1\end{pmatrix},\qquad
 \begin{pmatrix}1&0\\1&1\end{pmatrix}.
\]
Both are invertible, as is extraction of the first two coefficients of $(u+vx)/(1+x+x^2)$. Hence $M_w\in GL_2(\Z_2)$. Define
\begin{equation}\label{eq:odd-normalized-frame}
 E_w=\widehat E_wM_w^{-1},\qquad
 \mathcal G_{\alpha,w}=E_wR_\alpha^2\subset R_\alpha[[x]].
\end{equation}
The first two rows of $E_w$ form the identity. The plane $\mathcal G_{\alpha,w}$ therefore has exactly $4^\alpha$ elements, each uniquely specified by its first two coefficients.

\begin{lemma}\label{lem:odd-frame-history}
For $|w|=\alpha-1$, the plane $\mathcal G_{\alpha,w}$ is the central-row image of the compatible histories supported on $W_0$ with window $w$. For $1\le s\le\alpha$ and $|w|\ge\alpha-1$,
\begin{equation}\label{eq:odd-coherence}
 E_w\equiv E_{w_s}\pmod{2^s},
\end{equation}
where $w_s$ consists of the last $s-1$ digits of $w$.
\end{lemma}

\begin{proof}
The support calculation in Lemma~\ref{lem:odd-exceptional} shows that every $T_{t,v}$ preserves $W_0$. The triangular proof of Lemma~\ref{lem:count} consequently gives $4^\alpha$ compatible $W_0$-histories for a fixed window. Start with $M/f$, where $M\in W_0$. After $\alpha-1$ digits all retained zero-level numerators, including the oldest numerator $m_0=M$, lie in $W_0$. By Remark~\ref{rem:short-histories}, applied with $N$ replaced by $M\in W_0\subset L_{J_{\rm sh}}$, the compatibility and decoding arguments of Section~\ref{sec:count} apply at this length. The map $M\mapsto m_{\alpha-1}$ is linear; modulo $2$ its matrix is a product of the two displayed invertible matrices. Varying $M$ modulo $2^\alpha$ therefore supplies $4^\alpha$ different compatible histories, hence all of them. Their central rows are exactly $\widehat E_wR_\alpha^2$ and are distinct because $M_w$ is invertible.

For a longer word, any initial segment maps $W_0$ to itself. Modulo $2$ the induced map on the two coordinates is again invertible, so varying the initial $M$ still supplies the entire fixed-window history group. Now reduce a history from precision $\alpha$ to precision $s$, discard its older entries, and retain its terminal $s-1$ window digits. In a compatibility relation the discarded terms vanish because $T_{t,v}$ is divisible by $2^{t-1}$. In the decoding formula they vanish by Lemma~\ref{lem:precision}. Thus the reduced plane is contained in $\mathcal G_{s,w_s}$. It still has $4^s$ elements, since its first two coefficient map is invertible, so the inclusion is equality. Equal planes with first-two-row normalization have equal frames: the column with initial coefficients $(1,0)$, and likewise that with $(0,1)$, is unique. This proves \eqref{eq:odd-coherence}.
\end{proof}

\subsection{Primitive histories are reachable}

Write $W_0=\langle1,x\rangle_{\Z_2}$. For a window
$w\in\{0,1\}^{\alpha-1}$, let $\mathcal H_{\alpha,w}$ denote the
compatible histories of Definition~\ref{def:compatible} whose entries
belong to $W_0$, using the sheared operators. We call $H_0=u+vx$
\emph{primitive} if $u$ or $v$ is odd.

\begin{lemma}\label{lem:odd-history-coordinates}
For every window $w$, the map $H\mapsto H_0$ identifies
$\mathcal H_{\alpha,w}$ with $(\Z/2^\alpha)^2$. Reducing the precision
from $\alpha$ to $\beta\le\alpha$ reduces $H_s$ modulo $2^{\beta-s}$
for $s<\beta$ and keeps the last $\beta-1$ digits of the window.
These reductions commute with digit transitions. Each transition
between two ordinary window fibres is an invertible linear map on
$H_0$.
\end{lemma}

\begin{proof}
The exact matrices of $T_{1,r}$ on $W_0$, in row coordinates, are
\[
 G_{1,0}=\begin{pmatrix}1&d\\0&a\end{pmatrix},\qquad
 G_{1,1}=\begin{pmatrix}a&0\\1&d\end{pmatrix}.
\]
They are invertible modulo $2$. Suppose first that a compatible
history has $H_0=0$. Reducing its compatibility equations successively
modulo $2$ shows that every $H_s$ is even: at each step only
$T_{1,r}(H_{s+1})$ survives, because $T_{t,v}$ is divisible by
$2^{t-1}$. Divide the entries by $2$, omit the last entry, and retain
the last $\alpha-2$ window digits. This gives a compatible history at
precision $\alpha-1$ with first entry zero. Induction, starting with
$\alpha=1$, proves that the original history is zero. The triangular
count of Lemma~\ref{lem:count}, applied to the invariant module
$W_0$, gives $|\mathcal H_{\alpha,w}|=2^{2\alpha}$, so the injective
map to $H_0\in(\Z/2^\alpha)^2$ is bijective.

In a compatibility equation at position $s$, the terms discarded by
reduction to precision $\beta$ have $t\ge\beta-s$ and hence are
zero modulo $2^{\beta-s-1}$. The retained digit values depend only
on the retained suffix. Thus the asserted reduction preserves
compatibility and, by uniqueness, sends each canonical history with
given $H_0$ to the canonical reduced history. In the transition
formula~\eqref{eq:transition}, the terms $t>\beta$ vanish modulo
$2^\beta$; changing $H_{t-1}$ by a multiple of
$2^{\beta-t+1}$ also makes no change modulo $2^\beta$. The history
shift and window suffix update commute with reduction as well.
Finally the transition on $H_0$ is linear, and modulo $2$ is
$G_{1,r}$. Its determinant is odd, proving invertibility at every
precision.
\end{proof}

\begin{lemma}\label{lem:odd-reachability}
Every primitive history in $\mathcal H_{\alpha,w}$ is reached from
$R=N/f$ by a digit word. Within each window its decoded central-row
series are distinct. Consequently each window contributes exactly
$3\cdot4^{\alpha-1}$ distinct reachable primitive series.
\end{lemma}

\begin{proof}
At the zero window let $A_\alpha$ be the digit-$0$ matrix on $H_0$.
Lemma~\ref{lem:odd-history-coordinates} shows that these matrices are
coherent, and hence define $A\in\operatorname{GL}_2(\Z_2)$.
Modulo $2$, $A$ is $\left(\begin{smallmatrix}1&1\\0&1\end{smallmatrix}\right)$.
Therefore $A^2=I+2D_0$ for an integral matrix $D_0$, and repeated
squaring gives
\[
 A_\alpha^{2^\alpha}=I\pmod{2^\alpha}.
\]
Let $B_\alpha$ be the matrix of
the return word $10^{2^{\alpha+1}}$ starting at the zero window. It is
invertible. These matrices are coherent: after reduction to
precision $\beta$, the first $2^{\beta+1}$ zeros return the window to
zero, and the remaining zeros multiply by
$A_\beta^{2^{\alpha+1}-2^{\beta+1}}=I$, since
$2^\beta\mid 2^{\alpha+1}-2^{\beta+1}$. They define
$B\in\operatorname{GL}_2(\Z_2)$. Every positive word in
$A_\alpha,B_\alpha$ is an actual return word. Inverses can also be
realized by positive words, since each invertible matrix over the
finite ring $\Z/2^\alpha$ has finite order.

We now compute the return matrices algebraically, uniformly in all four odd parameters. Set
\[
 a=1+2r,\quad d=1+2s,\quad b=1+2u,\quad c=1+2v,\quad t=u+v.
\]
All matrices in the next calculation are modulo $8$. For
$i,k\in\{0,1\}$ put
\[
 (G_{j,e})_{i,k}=[x^{2^jk+e-i}y^0]\delta_j(f).
\]
The identities
\[
 \delta_1=f,\qquad \delta_2=f^3-f\sigma(f),\qquad
 \delta_3=f^7-f\sigma(f^3)-\delta_2\sigma^2(f)
\]
give the following finite coefficient table. Write
\[
 g=2+6r+4t+2s+4rs,\quad h=6+6r+4s,\quad j_*=4t+6s+4rs.
\]
Then
\[
\begin{gathered}
 G_{1,0}=\begin{pmatrix}1&1+2s\\0&1+2r\end{pmatrix},\qquad
 G_{1,1}=\begin{pmatrix}1+2r&0\\1&1+2s\end{pmatrix},\\
 G_{2,0}=\begin{pmatrix}0&g\\0&h\end{pmatrix},\quad
 G_{2,1}=\begin{pmatrix}2+4r&j_*\\0&g\end{pmatrix},\quad
 G_{2,2}=\begin{pmatrix}4+6r+4s&6+4t+6s\\2+4r&j_*\end{pmatrix},\\
 G_{3,0}=4\begin{pmatrix}0&1+r+t+rs\\0&r+s\end{pmatrix},\quad
 G_{3,1}=4\begin{pmatrix}1&t+rs\\0&1+r+t+rs\end{pmatrix},\\
 G_{3,2}=4\begin{pmatrix}1+r&r+s\\1&t+rs\end{pmatrix},\quad
 G_{3,4}=4\begin{pmatrix}r+s&s+rs\\1+r&t+s+rs\end{pmatrix}.
\end{gathered}
\]
For completeness these entries can be obtained directly without
parameter specialization: expand the three displayed polynomials,
discard terms of degree at least three in $2r,2s,2u,2v$, and use
$4z^2=4z\pmod8$. The latter reduction is valid for every
$z\in\Z_2$, since $z(z-1)$ is even. For example the coefficients
of $y^0$ in $\delta_2$ at $x$-degrees $0,\ldots,6$ are
\[
\begin{gathered}
0,\quad 2a,\quad 3a^2+2d-a,\quad a^3+6ad-a^2,\\
3a^2d+3d^2+6bc-ad-d,\quad a(3d^2+6bc-d),\quad
 d^3+6bcd-d^2,
\end{gathered}
\]
which give the second line of the table by the indicated coefficient
extraction. The third-level entries follow by substituting
$\delta_2=f^3-f\sigma(f)$ into the displayed formula for
$\delta_3$ and the same degree and square reductions.

Let
$P_0=\left(\begin{smallmatrix}1&1\\0&-1\end{smallmatrix}\right)$
and $J_0'=\left(\begin{smallmatrix}1&1\\0&1\end{smallmatrix}\right)$.
At precision $3$, the canonical zero-window history with
$H_0=z=(z_0,z_1)$ has
\[
 H_1=zP_0\pmod4,\qquad H_2=z\pmod2.
\]
Indeed compatibility says
$H_1=H_2G_{1,0}\pmod2$ and
$H_0=H_1G_{1,0}+H_2G_{2,0}\pmod4$.
Eliminating $H_2$ makes the latter matrix
$G_{1,0}+G_{1,0}^{-1}G_{2,0}=P_0\pmod4$, whose inverse is itself.
The transition formula and the table now give the following four
explicit matrix identities:
\[
\begin{aligned}
 A&=G_{1,0}+P_0G_{2,0}+G_{3,0}
       =\begin{pmatrix}1&5\\0&3\end{pmatrix},\\
 C_1&=G_{1,1}+P_0G_{2,2}+G_{3,4}
       =\begin{pmatrix}7&6\\3&1\end{pmatrix},\\
 C_2&=C_1G_{1,0}+G_{2,1}+J_0'G_{3,2}
       =\begin{pmatrix}1&5\\7&6\end{pmatrix},\\
 C_3&=C_2G_{1,0}+C_1G_{2,0}+G_{3,1}
       =\begin{pmatrix}5&0\\7&5\end{pmatrix}.
\end{aligned}
\]
Here $C_1,C_2,C_3$ are the successive matrices after the digits $1,0,0$. In each, the three summands are the three terms of the transition formula~\eqref{eq:transition}, whose digit values $(\dg{Wr}01,\dg{Wr}02,\dg{Wr}03)$ are $(1,2,4)$, $(0,1,2)$ and $(0,0,1)$. Every parameter term cancels in
the four displayed products. Since $A^2=
\left(\begin{smallmatrix}1&4\\0&1\end{smallmatrix}\right)$ and
$A^4=I\pmod8$, the remaining fourteen zeros in $10^{16}$ give
\begin{equation}\label{eq:odd-return-mod8}
 A=\begin{pmatrix}1&5\\0&3\end{pmatrix},\qquad
 B=C_3A^{14}=\begin{pmatrix}5&4\\7&1\end{pmatrix}\pmod8.
\end{equation}

Modulo $4$ the orbit of $(1,0)$ contains all twelve primitive
vectors. The following table gives an actual producing word for
each vector; $e$ denotes the empty word.
\[
\begin{array}{c|c@{\qquad}c|c@{\qquad}c|c}
(1,0)&e &(1,1)&A &(0,1)&AB\\
(0,3)&ABA &(3,1)&ABB &(1,3)&ABAB\\
(3,2)&ABBA &(2,1)&ABBB &(1,2)&ABABA\\
(2,3)&ABABB &(3,3)&ABABBB &(3,0)&ABABBBA
\end{array}
\]
Put $X=A^2$ and $Y=AB^4A$. Multiplication in
\eqref{eq:odd-return-mod8} yields
\[
 X=I+4E\pmod8,\qquad Y=I+4D\pmod8,\qquad
 E=\begin{pmatrix}0&1\\0&0\end{pmatrix},\quad
 D=\begin{pmatrix}1&0\\1&1\end{pmatrix}.
\]
Repeated squaring gives, for $m\ge2$,
\[
 X^{2^{m-2}}=I+2^mE\pmod{2^{m+1}},\qquad
 Y^{2^{m-2}}=I+2^mD\pmod{2^{m+1}}.
\]
Indeed a perturbation divisible by $4$ acquires exactly one extra
factor $2$ upon squaring at the precision in question; its square
has higher valuation. For a row $z=(1,0)\pmod2$, the two rows
$zE=(0,1)$ and $zD=(1,0)$ modulo $2$ are independent. Thus these
powers correct either coordinate bit at level $m$ without changing
lower bits. Induction shows that the orbit of $(1,0)$ contains every
vector congruent to $(1,0)$ modulo $4$. For any primitive target
$z^*$, choose $g$ from the table with
$(1,0)g=z^*\pmod4$. The vector $z^*g^{-1}$ is congruent to
$(1,0)$ modulo $4$, so it is reachable by the lifting argument;
applying $g$ reaches $z^*$. This proves transitivity on primitive
rows modulo every $2^\alpha$ (and modulo $2$ by reduction).

There is an actual primitive starting history at the zero window.
The word $10$ leaves the exceptional language, so all subsequent
zero-level numerators lie in $W_0$. Modulo $2$, the first digit gives
$x+xy$, and digit $0$ sends this to $x$: it kills $xy$ and fixes
$x$. Each further zero digit fixes $x$ modulo $2$. Therefore
$10^{\alpha+1}$ has a wholly ordinary history at the zero window
with primitive $H_0$. Applying actual return words realizes the
whole primitive set there by transitivity. To reach another
prescribed window, append its $\alpha-1$ digits. Each transition is
invertible on $H_0$ by Lemma~\ref{lem:odd-history-coordinates}, so
the image of the primitive set is exactly the primitive set in the
new window. Uniqueness of the compatible history with a given
$H_0$ proves the claimed reachability.

Finally decoding is linear on each window. Modulo $2$ it is
$\CT_y((u+vx)/f)$, whose constant and $x$ coefficients are
$u$ and $v-au$. Their matrix has determinant one. Hence these same
two coefficient functionals at precision $\alpha$ give an invertible
matrix over $\Z/2^\alpha$, and the decoded series are distinct.
There are $4^\alpha-4^{\alpha-1}=3\cdot4^{\alpha-1}$ primitive
coordinate rows.
\end{proof}

\subsection{Universality of the normalized window frames}\label{subsec:odd-universality}
The following argument proves universality for arbitrary odd $2$-adic parameters, rather than by checking finitely many residue classes.

\begin{proposition}\label{prop:odd-universality}
Let $E_w^{(\theta)}$ denote \eqref{eq:odd-normalized-frame} for an odd quadruple $\theta=(a,b,c,d)$, and let $E_w^{(1)}$ denote it for $(1,1,1,1)$. If $|w|=k$, then
\[
 E_w^{(\theta)}\equiv E_w^{(1)}\pmod{2^{k+1}}.
\]
In particular, the full-window frames at precision $\alpha$, and their planes, are independent of $\theta$.
\end{proposition}

\begin{proof}
Put $q=1+x+x^2$, $D=q^2-4x^4$, $B=D^{-1/2}$ and
\[
 G=((1+x)B,xB),\qquad G_\theta=((1+ax)B_\theta,xB_\theta).
\]
These rows have the same section planes as $(B,xB)$ and $(B_\theta,xB_\theta)$, respectively. We work in the $2$-adic completion of $\Z_2[x,q^{-1}]$; all rational denominators have constant term one. For $m\ge0$ let
\[
 \mathcal R_m=\{N/q^m:\deg N\le2m+1\},\qquad
 \mathcal R_m^0=\{N/q^m:\deg N\le2m\},
\]
with integral polynomial numerators. They are nested in $m$, and
$\mathcal R_m^0\mathcal R_n\subset\mathcal R_{m+n}$. Congruences in this completed rational ring imply coefficientwise series congruences.

Write $P_\theta=q+2H$, $bc=1+2v_0$, and $Z=qH+H^2-bcx^4$. Thus $D_\theta=P_\theta^2-4bcx^4=q^2+4Z$. Set
\[
 T=\frac{x^4(v_0-2H/q-2H^2/q^2)}{D_\theta}.
\]
A direct calculation gives
\begin{equation}\label{eq:odd-deformation-root}
 (B_\theta/B)^2=(q/P_\theta)^2(1+8T),\qquad
 B_\theta/B=(q/P_\theta)\sum_{r\ge0}\binom{1/2}{r}8^rT^r.
\end{equation}
The second identity uses the square root with constant term one. For $r\ge1$,
\[
 v_2\left(8^r\binom{1/2}{r}\right)=r+s_2(r),
\]
where $s_2(r)$ counts the binary ones of $r$. Indeed
$\binom{1/2}{r}=(-1)^{r-1}\binom{2r}{r}/(4^r(2r-1))$, and counting multiples of powers of two gives $v_2(n!)=n-s_2(n)$. This also proves convergence of \eqref{eq:odd-deformation-root}.

The expansion
\[
 D_\theta^{-1}=q^{-2}\sum_{h\ge0}(-4Z/q^2)^h
\]
shows that $T\in\sum_{j\ge0}2^j\mathcal R_{2j+2}^0$. To check the indices, the three terms in its numerator have (precision, pole) pairs $(0,2)$, $(1,3)$ and $(1,4)$, while the $h$th geometric term adds $(2h,2h)$. Thus
$T^r\in\sum_{j\ge0}2^j\mathcal R_{2j+2r}^0$.
A nonconstant square-root term of total precision $n$ in \eqref{eq:odd-deformation-root} consequently has pole at most $2n-2s_2(r)\le2n-2$. Combining this with
$q/P_\theta=\sum_{h\ge0}(-2H/q)^h$ gives
\[
 B_\theta/B-1\in2\mathcal R_1^0+\sum_{n\ge2}2^n\mathcal R_{2n-2}^0,
 \qquad
 (G_\theta-G)/B\in2\mathcal R_1+\sum_{n\ge2}2^n\mathcal R_{2n-2}.
\]
For the second inclusion multiply by $(1+x,x)$ and add $((a-1)xB_\theta/B,0)$; the extra factor $a-1$ is even and obeys the same bounds.

Let $H_j=x^j\partial_x^{[j]}$, where $\partial_x^{[j]}x^n=\binom njx^{n-j}$ is the $j$th Hasse derivative, and put $V_j=(H_jG)/B$. Differentiation shows
\[
 V_j=N_j/D^j,\qquad \deg N_j\le4j+1,
\]
with integral polynomial numerators. For integrality, the series $H_jG$, $B^{-1}$ and $D^j$ are integral, and their product is the polynomial $N_j$; no division by $j!$ in $\Z_2$ is used.
Modulo $2$, the columns of $V_j$ have pole at most $j$ at $q$. At a root $\rho$ of $q$ in $\F_4$, their leading principal parts for $j>0$ are a common nonzero multiple of $(1+\rho,\rho)$. The latter two entries are independent over $\F_2$: their sum is $1$, and neither is in $\F_2$. Hasse differentiation of a simple pole has a nonzero leading term at every order. At $j=0$ the polynomials $1+x,x$ are a basis. It follows successively by pole order that the $2m+2$ columns of $V_0,\ldots,V_m$ are a basis of $\mathcal R_m/2\mathcal R_m$.

For $j>0$ choose any lift $v_j\in\mathcal R_j$ of $V_j\bmod2$. Since
\[
 V_j=\frac{N_j}{q^{2j}}
  \sum_{h\ge0}\binom{j+h-1}{h}(4x^4/q^2)^h,
\]
the leading difference has even numerator of pole at most $2j$, and the remaining terms have precision at least $2h$. Therefore
\begin{equation}\label{eq:odd-hasse-tail}
 V_j-v_j\in\sum_{t\ge1}2^t\mathcal R_{2j+2\lfloor t/2\rfloor}.
\end{equation}
For $j=0$ take $v_0=V_0$, so that its tail is zero.

Define filtered tails
\[
 \mathcal D_{\ge n}=\sum_{m\ge n}2^m\mathcal R_{2^{m-1}}\qquad(n\ge1).
\]
The initial deformation lies in $\mathcal D_{\ge1}$, because $2m-2\le2^{m-1}$ for $m\ge2$. At stage $n$, expand its coefficient of $2^n$ modulo $2$ in the basis given by $V_0,\ldots,V_{2^{n-1}}$, and subtract the corresponding $2^n$ multiples with constant binary matrix coefficients. The residual lies in $\mathcal D_{\ge n+1}$: the leading difference is even, and \eqref{eq:odd-hasse-tail} fits the later bounds since, for $j\le2^{n-1}$ and $t\ge1$,
\[
 2j+2\lfloor t/2\rfloor\le2^n+2\lfloor t/2\rfloor
      \le2^{n+t-1}.
\]
For $t=1$ this is equality at the largest $j$; for $t\ge2$ it follows from $2\lfloor t/2\rfloor\le2^n(2^{t-1}-1)$.
Iteration yields convergent constant $2\times2$ matrices $A_j$ such that
\begin{equation}\label{eq:odd-hasse-expansion}
 G_\theta=\sum_{j\ge0}(H_jG)A_j,\quad
 A_0\equiv I\pmod2,\quad
 v_2(A_j)\ge1+\lceil\log_2j\rceil\quad(j\ge1).
\end{equation}
The matrices tend to zero $2$-adically and the $V_j$ are integral, so $\sum_jV_jA_j=G_\theta/B$ converges in the rational completion. Multiplication by $B$ gives the displayed identity in $\Z_2[[x]]$, with convergence in the $2$-adic topology.

Coefficient extraction gives the exact operator identity
\[
 \Lambda_{2^k,v}H_j=\binom{2^kx\partial_x+v}{j}\Lambda_{2^k,v}.
\]
Vandermonde's identity gives, for $n\ge0$,
\[
 \binom{2^kn+v}{j}-\binom vj
 =\sum_{i=1}^j\binom{2^kn}{i}\binom v{j-i}.
\]
If $\lfloor\log_2j\rfloor\le k$, each summand is divisible by
$2^{k-\lfloor\log_2j\rfloor}$, since
$\binom{2^kn}{i}=(2^kn/i)\binom{2^kn-1}{i-1}$; for $n=0$ it is zero. If $\lfloor\log_2j\rfloor>k$, integrality suffices. After multiplication by $A_j$ in \eqref{eq:odd-hasse-expansion}, the nonconstant-in-$n$ part thus vanishes modulo $2^{k+1}$. We obtain
\[
 \Lambda_{2^k,v}G_\theta\equiv(\Lambda_{2^k,v}G)C_v\pmod{2^{k+1}},
 \qquad C_v=\sum_{j=0}^v\binom vjA_j\equiv I\pmod2.
\]
The constant matrix $C_v$ is invertible and disappears on first-two-coefficient normalization. The proposition follows.
\end{proof}

\subsection{Weak contact and the energy bound}\label{subsec:odd-energy}
For two windows of length $\alpha-1$, define their clipped contact by
\[
 r_\alpha(w,v)=\max\{0\le s\le\alpha:E_w\equiv E_v\pmod{2^s}\}.
\]
For distinct windows let $\ell(w,v)$ be the length of their longest common terminal suffix; this counts the last applied digits.

\begin{lemma}\label{lem:odd-weak-contact}
For distinct full windows,
\[
 r_\alpha(w,v)\le\min\{\alpha,\ell(w,v)+3\}.
\]
For every $s\ge1$, an equivalence class of length-$(s-1)$ windows with equal normalized frames modulo $2^s$ has at most four members.
\end{lemma}

\begin{proof}
By Proposition~\ref{prop:odd-universality} it suffices to use the all-one parameters. Proposition~\ref{prop:odd-rational-contact} in Appendix~\ref{app:odd-certificate} proves the asserted contact bound, by a normal rational lattice and a ten-dimensional quotient on which the normalized difference cannot vanish.
For the second statement let $M_s$ be the largest equivalence-class size. If $s\le3$ there are at most $2^{s-1}\le4$ windows in total. If $s\ge4$, any two different members of a class have contact $s$, so their common terminal suffix has length at least $s-3$. Once a member is fixed, only its two older digits can vary, leaving at most four members.
\end{proof}

\begin{proof}[Proof of Theorem~\ref{thm:odd-exponent}]
The upper bound is Proposition~\ref{prop:odd-upper}. For the lower bound, put $n=2^{\alpha-1}$, and let $\mathcal P_w$ be the primitive part of $\mathcal G_{\alpha,w}$, namely $E_wz$ for coordinate columns $z\in R_\alpha^2$ with at least one odd entry. By Lemma~\ref{lem:odd-reachability}, all these states belong to the kernel of $F\bmod2^\alpha$, and
\[
 |\mathcal P_w|=4^\alpha-4^{\alpha-1}=3\,4^{\alpha-1}.
\]
If a state belongs to both planes, its coordinate column in both is the same, since the first two frame rows are the identity. Hence $(E_w-E_v)z=0$. If $r=r_\alpha(w,v)<\alpha$, a row attaining the smallest valuation gives one congruence
$c_0z_0+c_1z_1=0\pmod{2^\alpha}$ with $\min(v_2(c_0),v_2(c_1))=r$. Dividing by $2^r$ leaves a unit coefficient modulo $2^{\alpha-r}$, so there are exactly $2^{\alpha+r}$ solutions in $R_\alpha^2$. If $r=\alpha$, the same bound is the full-plane size $4^\alpha$. Consequently
\begin{equation}\label{eq:odd-intersection}
 |\mathcal P_w\cap\mathcal P_v|\le2^{\alpha+r_\alpha(w,v)}.
\end{equation}

At precision $s$, let $m_{s,1},\ldots,m_{s,t_s}$ be the sizes of the frame equivalence classes among length-$(s-1)$ windows; their sum is $2^{s-1}$ and their maximum is $M_s$. Coherence \eqref{eq:odd-coherence} gives the exact ordered-pair count
\[
 N_s:=\#\{(w,v):r_\alpha(w,v)\ge s\}
   =2^{2(\alpha-s)}\sum_jm_{s,j}^2
   \le M_s2^{2\alpha-s-1}.
\]
Since $2^r=1+\sum_{s=1}^\alpha2^{s-1}\one_{r\ge s}$,
\begin{equation}\label{eq:odd-contact-energy}
 \sum_{w,v}2^{r_\alpha(w,v)}
 \le n^2\left(1+\sum_{s=1}^\alpha M_s\right).
\end{equation}
Let $h(g)=\#\{w:g\in\mathcal P_w\}$. Equations \eqref{eq:odd-intersection}--\eqref{eq:odd-contact-energy} imply
\[
 \sum_gh(g)=n\,3\,4^{\alpha-1},\qquad
 \sum_gh(g)^2\le2^\alpha n^2
       \left(1+\sum_{s=1}^\alpha M_s\right).
\]
Cauchy--Schwarz therefore gives
\[
 K_\alpha(F)\ge\left|\bigcup_w\mathcal P_w\right|
 \ge\frac{9\,2^{3\alpha-4}}{1+\sum_{s=1}^\alpha M_s}.
\]
There is one empty window at $s=1$, so $M_1=1$; Lemma~\ref{lem:odd-weak-contact} gives $M_s\le4$. Thus the denominator is at most $4\alpha-2$, yielding the stated lower bound, including $\alpha=1$. Taking logarithms of the two bounds gives liminf at least $3$ and limsup at most $3$, since the logarithmic loss $\log_2\alpha$ is $o(\alpha)$.
\end{proof}

\section{Rank versus cardinality}\label{sec:rank}

For a sequence $g$ we write $\kerp(g)$ for its $p$-kernel. A series over $\QQ$ is called D-finite if it satisfies a nonzero homogeneous linear differential equation with polynomial coefficients. For the two examples in Proposition~\ref{prop:rank}, see \cite[Theorem~2.1]{Stan80} for algebraic series and \cite[Theorem~3.1]{BLS17} for rational diagonals. The two Mahler results used below are Becker's implication from regularity to a functional equation \cite[Theorem~1]{Bec94} and B\'ezivin's rationality theorem \cite{Bez94}, in the form of \cite[Theorem~1]{BCR13}. For the algebraic case see \cite[Corollary~8]{BCR13}, which identifies \cite[Theorem~5.1.7]{Nis96} as its source.

\begin{proposition}\label{prop:rank}
Let $F\in\Z[[x]]$ be D-finite but not rational; for example, $F$ may be an irrational algebraic series, or an irrational diagonal of a rational function with integer coefficients \cite{Lip88}. Let $p$ be a prime. For $\alpha\ge1$ let $\rho_\alpha$ be the least $r$ such that $\kerp(F\bmod p^\alpha)$ is contained in a $\Z/p^\alpha$-submodule of $(\Z/p^\alpha)^{\Z_{\ge0}}$ generated by $r$ elements. Set $\rho_\alpha=+\infty$ if there is no such finite $r$. Then $\rho_\alpha\to\infty$. In particular, the minimal number of states of a $p$-linear scheme generating $F\bmod p^\alpha$ tends to infinity with $\alpha$.
\end{proposition}

\begin{proof}
If $\kerp(F\bmod p^{\alpha+1})$ lies in a submodule generated by $r$ elements, then reducing modulo $p^\alpha$ shows that $\kerp(F\bmod p^\alpha)$ lies in a submodule generated by the $r$ reductions. Hence $\rho_\alpha\le\rho_{\alpha+1}$. Suppose $\rho_\alpha\le r$ for all $\alpha$. Let $w_0,\dots,w_r$ be words and $n_0,\dots,n_r$ indices, and let $\mu=\det\bigl([x^{n_l}]\Lambda_{w_i}F\bigr)_{i,l}\in\Z$. Modulo $p^\alpha$ the rows lie in a module generated by $r$ sequences $e_1,\dots,e_r$, so the matrix factors as $XY$ with $X$ of size $(r+1)\times r$ and $Y$ of size $r\times(r+1)$ over $\Z/p^\alpha$. Expand the determinant multilinearly in the rows of $XY$, each a combination of the $r$ rows of $Y$. Every choice of $r+1$ rows repeats one, so alternation gives $\det(XY)=0$. Hence $p^\alpha\mid\mu$ for every $\alpha$, so $\mu=0$. All $(r+1)$-minors of the matrix with rows $\Lambda_wF$ therefore vanish. If $r+1$ such sequences were linearly independent, their coordinate functionals would span the dual of their $(r+1)$-dimensional span, since a sequence killed by every coordinate is zero. Choosing a basis of these functionals would give a nonzero finite minor, a contradiction. Thus the $\QQ$-span of $\{\Lambda_wF\}$ has dimension at most $r$. That is, $F$ is $p$-regular over $\QQ$.

A $p$-regular power series is $p$-Mahler \cite{Bec94}; we recall the argument. Choose a $\QQ$-basis $H=(H_1=F,\dots,H_d)^T$ of the Cartier span. The identity $H_i(x)=\sum_{r=0}^{p-1}x^r(\Lambda_rH_i)(x^p)$ gives $H(x)=A(x)H(x^p)$ for a polynomial matrix $A$. Iteration expresses each $F(x^{p^i})$, $0\le i\le d$, as a polynomial row vector times $H(x^{p^d})$. Linear dependence of these $d+1$ row vectors over $\QQ(x)$, followed by clearing denominators, gives $\sum_{i=j}^d a_i(x)F(x^{p^i})=0$ with $a_j\ne0$. Write $a_i(x)=\sum_{v=0}^{p^j-1}x^v b_{i,v}(x^{p^j})$. Separating the residue classes of exponents gives $\sum_{i=j}^d b_{i,v}(x)F(x^{p^{i-j}})=0$ for every $v$. Choose $v$ with $b_{j,v}\ne0$ and omit any zero final coefficients. Since $F\ne0$, this is a $p$-Mahler equation with both endpoint coefficients nonzero.

By a theorem of B\'ezivin \cite{Bez94} (see \cite[Theorem~1]{BCR13}), a D-finite $p$-Mahler power series in $\mathbb C[[x]]$ is a rational function. (For algebraic $F$ one may instead use the rational--transcendental dichotomy for Mahler functions \cite[Theorem~5.1.7]{Nis96}, \cite{BCR13}.)

A series with rational coefficients that is rational over $\mathbb C$ is rational over $\QQ$: normalize a complex denominator to have constant term $1$. Its finitely many unknown coefficients solve rational linear equations expressing the vanishing of the product beyond the numerator degree. A finite set of rows spans this system (including its constant terms), and Gaussian elimination gives a rational solution. The resulting numerator is also rational. This contradicts the hypotheses on $F$. The last statement holds because the kernel of a sequence generated by a scheme with $s$ states lies in a module generated by $s$ elements.
\end{proof}

For example, the series of Corollary~\ref{cor:family} with $a=b=c=d=1$ is irrational: its discriminant $(1+x+x^2)^2-4x^4=(1+x-x^2)(1+x+3x^2)$ is not a square: both factors are squarefree, with discriminants $5$ and $-11$, and their difference is $4x^2$, so they have no common root. Their nonconstant product is squarefree and cannot be a square in $\QQ(x)$. Therefore no description of its kernels by modules of bounded rank exists.

By contrast, for each window $W$ the set $\Hc_W$ is a group of order exactly $p^{|J|\alpha}$ (Lemma~\ref{lem:count}), and $\Psi(\cdot,W)$ is additive. The multivariate sections $C_wR$ for $|w|\ge\alpha$ thus lie in a union of $p^{\alpha-1}$ subgroups of $(\Z/p^\alpha)[[z]]$, each of order at most $p^{|J|\alpha}$. Applying $\Diag$ gives the corresponding statement for the univariate kernel elements in $(\Z/p^\alpha)[[x]]$. We make no claim about the ranks of these groups. Even in the trivial case $Q=N=1$ (so $J=\{0\}$, $\delta_1=1$ and $\delta_t=0$ for $t\ge2$), the group $\Hc_W$ for $W=0^{\alpha-1}$ is cyclic of order $p^\alpha$, so its rank is $|J|$ rather than $|J|\alpha$. The bound of Theorem~\ref{thm:main} is a statement about cardinality.

\emph{Rational special cases.} Precision-linear exponents are already elementary for sequences with a fixed finite matrix representation, so the one-variable rational specialization is not claimed as new. If $a(n)=\lambda A^nv$ with $A\in M_r(\Zp)$, then $a(p^kn+s)=(\lambda A^s)(A^{p^k})^nv$. Modulo $p^\alpha$ there are at most $p^{(r+r^2)\alpha}$ pairs $(\lambda A^s,A^{p^k})$. A rational series has such a representation: normalize $Q=1+q_1x+\cdots+q_dx^d$ and put $r=\max(d,\deg N+1,1)$. The identity $Q\sum_na(n)x^n=N$ gives $a(n)=-\sum_{i=1}^d q_i a(n-i)$ for $n\ge r$. Padding the coefficients to order $r$ gives a fixed matrix on the vector $(a(n),\dots,a(n+r-1))^T$ for every $n\ge0$, with output its first coordinate. The new general conclusion concerns arbitrary rational diagonals, for which a fixed-rank representation need not exist.

\section{Computations}\label{sec:computations}

The computations in this section are exact finite checks and do not establish asymptotic laws. The finite proof certificate for Theorem~\ref{thm:odd-exponent} is given separately in Appendix~\ref{app:odd-lattice}.

\subsection{The family with \texorpdfstring{$a=b=c=d=1$}{a=b=c=d=1}}
Let $F$ be the Furstenberg series of $P=x^2y^2+(1+x+x^2)y+x^2$ and $p=2$. Table~\ref{tab:main} reports exact computations of its minimal kernel size. All state identities, polynomial identities and partition refinements in these computations use exact arithmetic and full equality checks.

\begin{table}[ht]
\centering
\small
\begin{tabular}{rrrrrr}
\toprule
$\alpha$ & $K_\alpha(F)$ & history states & rational states & $K_\alpha/8^\alpha$ & $\log_2K_\alpha$ \\
\midrule
1 & 14 & 14 & 14 & 1.750 & 3.81\\
2 & 53 & 71 & 53 & 0.828 & 5.73\\
3 & 277 & 466 & 303 & 0.541 & 8.11\\
4 & 1\,357 & 2\,875 & 1\,519 & 0.331 & 10.41\\
5 & 6\,753 & 19\,597 & 8\,676 & 0.206 & 12.72\\
6 & 38\,741 & 143\,777 & 56\,335 & 0.148 & 15.24\\
7 & 353\,479 & 1\,098\,769 & 399\,129 & 0.169 & 18.43\\
8 & 2\,747\,455 & 8\,587\,651 & 2\,991\,449 & 0.164 & 21.39\\
9 & 21\,151\,568 & --- & 23\,136\,226 & 0.158 & 24.33\\
\bottomrule
\end{tabular}
\caption{Exact finite kernel computations for the quadratic family at $p=2$. ``Rational states'' counts the states reached in the independent construction below, before minimization. The history construction was run through $\alpha=8$; the rational construction through $\alpha=9$. Theorem~\ref{thm:odd-exponent} determines the asymptotic exponent.}\label{tab:main}
\end{table}

\emph{Two constructions.} The first is the history automaton of Corollary~\ref{cor:automaton}. Our implementation builds its transition matrices by exact integer polynomial multiplication, uses full equality to identify states, and minimizes the complete graph. It gives the history-state counts and $K_\alpha$ through $\alpha=8$.

The independent rational construction does not use histories or the $\delta$-decomposition. Put $f=1+x+x^2(y+1+y^{-1})$, $N=y+xy+x^2y+2x^2y^2$, $E=(f^2-\sigma(f))/2$, and $e=\lceil\alpha/2\rceil$. Represent a state by $V/f^\alpha$. The binomial expansion, modulo $2^\alpha$, gives
\[
C_r(V/f^\alpha)=\sum_{j=0}^{\alpha-1}c_j\,
 \frac{C_r(Vf^{\alpha\bmod2}E^j)}{f^{e+j}},
\qquad c_j=(-2)^j\binom{e+j-1}{j},
\]
where here $C_r=C\circ x^{-r}$ in the sheared variables of the appendix. Indeed, write $f^\alpha=f^{2e}/f^{\alpha\bmod2}$, expand $(\sigma(f)+2E)^{-e}$, discard terms divisible by $2^\alpha$, and use the projection formula for $C$. To return to denominator $f^\alpha$, form
\[
W=\sum_{j=0}^{\alpha-1}c_j\,C_r(Vf^{\alpha\bmod2}E^j)f^{\alpha-1-j}
\]
and compute a candidate polynomial $V'$ with $W=f^{e-1}V'$ modulo $2^\alpha$. The complete polynomial residual, including coefficients outside the stored frame, is checked on every reached state; thus no truncated-series equality is assumed. The initial numerator is $Nf^{\alpha-1}$, and its two outgoing transitions are checked separately.

Both $f$ and $E$ are invariant under $y\mapsto y^{-1}$. The Cartier sections commute with this involution, and the central-row output kills every $y^j-y^{-j}$. Therefore one may quotient by their span, storing the sum of the coefficients at $j$ and $-j$ for $j>0$, with the coefficient at $j=0$ unchanged. No division by $2$ is used. Multiplication is always by a symmetric polynomial. A finite frame $0\le i\le2\alpha$, $0\le j\le\lceil i/2\rceil$ supplies the candidate numerator; its validity follows from the full residual checks. All residual constraints are obtained by exact polynomial operations; invertible row operations reduce them before the graph is built.

At $\alpha=8$, the coefficients $c_5,c_6,c_7$ vanish modulo $256$. The transition is consequently the polynomial
\[
V'=\sum_{j=0}^4c_j\,C_r(VE^j)f^{4-j}\pmod{256},
\]
and every full residual vanishes identically. A separate implementation independently recomputed all $2\cdot89^2$ matrix entries and both initial transition vectors. The resulting complete graph has $2\,991\,449$ reachable states. Independent sorting-based partition refinement gives $65\,536$, $2\,369\,812$, $2\,747\,455$, $2\,747\,455$ classes, confirming $K_8$ without using a history graph.

The complete synchronous product with the independently rebuilt history graph at $\alpha=8$ has $8\,587\,651$ reachable pairs, and all paired outputs agree. This checks equality for every input word, beyond the agreement of the two minimum counts.

At $\alpha=9$ the graph has $23\,136\,226$ states, and all $1\,503\,854\,625$ reduced residual constraints passed; the minimum has $21\,151\,568$ states. A separate implementation reconstructed both transition matrices, the complete residual matrices and initial vectors by direct polynomial division, and independently checked graph reachability and the minimum by sorting-based refinement. This is one mathematical construction with separately checked matrices and minimization, rather than agreement of two constructions at that precision. Building and minimizing the graphs took $6.70$ seconds at $\alpha=8$ and $85.39$ seconds at $\alpha=9$ on a single core; matrix generation is logged separately. A third, separately written minimizer reproduces both minima and the final partitions. Direct coefficient recurrence checks agree for $n<4096$ at every $\alpha\le9$, and for $n<2^{18}$ at $\alpha=8$ and $\alpha=9$. Source code, complete logs, matrix certificates, graph checksums and regeneration instructions are in \texttt{anc/}.

\subsection{Other series}
Table~\ref{tab:other} collects further minimal kernel sizes. All $57$ entries were reproduced by new exact code: the non-multinomial entries by an independent implementation of the history construction, and the multinomial row by a separate fixed-denominator polynomial Cartier construction. The latter construction also independently reproduces every Catalan and Example~13 entry, all listed entries at $p=3,5$, and the first three entries in each other odd-parameter row. Support-set membership is certified by rational primal/dual or infeasibility certificates; floating-point linear programming is used only to propose candidates. Exact random tests in two and three variables at $p=2,3,5$ are logged separately.

\begin{table}[ht]
\centering
\small
\begin{tabular}{llll}
\toprule
series & $p$ & $|J|$ & $K_1,K_2,\dots$\\
\midrule
family, $(a,b,c,d)=(3,5,7,-1)$ & 2 & 5 & 14, 61, 266, 1345, 6736, 38693, 353099\\
family, $(-1,-3,5,9)$ & 2 & 5 & 14, 61, 268, 1331, 6762, 38572, 353435\\
family, $(5,1,-7,3)$ & 2 & 5 & 14, 61, 268, 1323, 6750, 38769, 353640\\
family, $(1,1,1,3)$ & 2 & 5 & 14, 61, 264, 1334, 6748, 38397, 353101\\
family, $(1,1,1,1)$ & 3 & 5 & 15, 134, 2058\\
family, $(1,1,1,1)$ & 5 & 5 & 11, 454\\
$C(x)-1$ (Catalan), $P=x(1+y)^2-y$ & 2 & 4 & 4, 6, 13, 28, 58, 136, 309\\
$C(x)-1$ & 3 & 4 & 6, 23, 113, 629\\
\cite[Ex.~13]{RY26}, $P=xy^2+(x+1)y+x$ & 2 & 4 & 7, 15, 23, 79, 220, 707\\
\cite[Ex.~13]{RY26} & 3 & 4 & 6, 24, 116\\
$(3n)!/n!^3$, $Q=1-z_1-z_2-z_3$ & 2 & 1 & 2, 3, 8, 14\\
\bottomrule
\end{tabular}
\caption{Minimal kernel sizes for further series; $|J|=|J_p(N,Q)|$ is computed from the integral polynomials $N,Q$ of the representation used in each row. Every listed entry was independently reproduced; the multinomial row is confirmed through $\alpha=4$.}\label{tab:other}
\end{table}

\subsection{Finite ratios and the proved exponent}
For the family with $a=b=c=d=1$, the ratio $K_\alpha/8^\alpha$ equals $0.148$, $0.169$, $0.164$ and $0.158$ for $\alpha=6,7,8,9$, and the ratios of successive history-state counts are $7.3$, $7.6$ and $7.8$. Theorem~\ref{thm:odd-exponent} proves the exponent $3$ for the entire odd family independently of these data.

\section{Open questions}\label{sec:questions}

\begin{enumerate}[label=\arabic*.]
\item \emph{The general exponent.} Determine whether $\lim_\alpha\alpha^{-1}\log_pK_\alpha(G)$ exists for a general rational diagonal, and compute it in terms of $N$, $Q$ and $p$. For the odd quadratic family the limit is $3$. In all examples of Section~\ref{sec:computations} the observed growth is well below the slope $|J|+1$; for the odd quadratic family the slope is $3$, while $|J|+1=6$. We know no example in which the slope $|J|+1$ is attained. Is the general slope governed by $|J|$, or by a unit-root part of the operators $T_{t,v}$?
\item \emph{The leading constant.} Does $K_\alpha(F)/8^\alpha$ converge for the odd quadratic family? Determine its limit, if it exists, and its dependence on the parameters. Can the factor $\alpha$ in the lower bound of Theorem~\ref{thm:odd-exponent} be removed?
\item \emph{Lower bounds beyond the odd family.} Find sharp exponential lower bounds for other natural classes of rational diagonals, and identify conditions under which their upper and lower exponents agree.
\item \emph{Beyond a chosen unit-denominator representation.} Treat algebraic series whose given annihilating polynomial has nonunit derivative at the origin, distinguishing failure of the Furstenberg conversion from nonexistence of another unit-denominator diagonal representation. For the central trinomial coefficients, the equation $(x+1)(3x-1)y^2+1=0$ has $y$-derivative zero modulo $2$. Rowland and Yassawi discuss this example in \cite[\S10, p.~37, before Lemma~49]{RY26}. A different representation gives
\[
\sum_{n\ge0}\CT_y(y^{-1}+1+y)^n x^n
=\Diag\frac1{1-x(1+y+y^2)}.
\]
Indeed, the diagonal coefficient is $[y^n](1+y+y^2)^n$. For this representation $J=\{0\}$: membership implies $u_1+\eta\le1-\varepsilon<1$, hence $u_1=0$, and $u_2+\eta\le2\eta$, hence $u_2=0$; conversely $0\in J$. Theorem~\ref{thm:main} therefore gives $K_\alpha\le p^{2\alpha-1}+(p^\alpha-1)/(p-1)$.
\item \emph{Sharper support sets.} Find stable sets smaller than $J$, for instance by examining the ranks and invariant subspaces of the digit operators modulo $p$, and decide whether such refinements reach the true slope.
\item \emph{Uniformity.} Prove bounds that are good simultaneously in $\alpha$, in the degrees and in $p$. Can the two-variable bounds of \cite[Theorems~1 and~56]{RY26} and the count in Theorem~\ref{thm:main} be combined into a single estimate that captures their respective advantages?
\item \emph{Structure.} Describe the groups $\Hc_W$ (their Smith forms and ranks) and the fibres of $\Psi$. How much of the gap between Table~\ref{tab:main} and the bound comes from non-injectivity of $\Psi$, and how much from unreachable histories?
\item \emph{The zero orbit.} Determine the exact preperiod of the observed sequences $\Lambda_0^kF$ in Theorem~\ref{thm:zero}. Which other polynomials admit an exact formula for this orbit, including cases where the first matrix on the proposed invariant subspace is singular modulo the prime?
\item \emph{Computations.} Recompute $K_9$ by a second mathematical construction and extend the computations to $\alpha=10$.
\end{enumerate}

\section*{Use of AI tools}
The author directed the preparation and verification of this manuscript and selected the results included in it. Claude (Anthropic) assisted in verifying the zero-orbit proof, finding and proving the linear kernel bound for the two-variable family, independently checking the proof of Theorem~\ref{thm:odd-exponent}, searching the literature, performing computations, and drafting the text. ChatGPT (OpenAI), including Codex, was used to independently verify the proofs, extend the bound to diagonals in any number of variables, audit the literature, perform independent exact computations, and prepare the submission files. Under the author's direction, ChatGPT was also used to derive the proof of Theorem~\ref{thm:odd-exponent} and to draft Section~\ref{sec:odd-exponent} and Appendix~\ref{app:odd-certificate}; a separate ChatGPT session audited this proof independently. For the second version, ChatGPT revised the text and Claude checked the revision against the sources; an independent Claude session found the sharper form of Corollary~\ref{cor:furst} and the missing step in the proof of Lemma~\ref{lem:odd-finite}. Before submission, two further independent sessions, one of each system, audited the text; the ChatGPT session reran the ancillary code and recomputed Tables~\ref{tab:main} and~\ref{tab:other}. The author checked all results and takes full responsibility for the content of this paper.

\appendix

\section{The exact zero orbit for the quadratic family}\label{app:zero}

In this appendix $p=2$, $a,b,c,d\in\Z_2^\times$, and
\[
P=bx^2y^2+P_1(x)\,y+cx^2,\qquad P_1(x)=1+ax+dx^2,\qquad f=P/y=1+ax+x^2(by+d+cy^{-1}),
\]
\[
N=y\,\partial_yP=P_1(x)\,y+2bx^2y^2,\qquad R=N/f .
\]
We use the coordinate change $x\mapsto xy^{-1}$ of \cite[\S2, p.~7, before (2)]{RY26}, which sends a diagonal coefficient to a coefficient of $y^0$. Let $A=\Z_2[y,y^{-1}][[x]]$ be the ring of power series in $x$ whose coefficients are Laurent polynomials in $y$; $f$ is a unit of $A$. Let $\sigma(x)=x^2$, $\sigma(y)=y^2$, let $\Lambda_{r,0}$ act by $\sum c_{ij}x^iy^j\mapsto\sum c_{2i+r,2j}x^iy^j$, put $C=\Lambda_{0,0}$, and let $\CT_y$ take the coefficient of $y^0$. Lemma~\ref{lem:sections}(a) holds in $A$ with the same proof, and $\CT_y\circ C=\Lambda_0\circ\CT_y$.

Let $\widehat A_+$ be the set of formal sums $\sum c_{ij}x^iy^j$ ($c_{ij}\in\Z_2$) supported on $\{i\ge0,\ i+j\ge0\}$. Each coefficient of a product of two such sums is a finite sum, so $\widehat A_+$ is a ring. It contains $A_+:=A\cap\widehat A_+$ as a subring, and $f,N,R\in A_+$. The substitution $\varphi(g)(x,y)=g(xy,y)$, that is $x^iy^j\mapsto x^iy^{i+j}$, is a ring isomorphism $\widehat A_+\to\Z_2[[x,y]]$, and
\[
\varphi\circ\sigma=\sigma\circ\varphi,\qquad\varphi\circ\Lambda_{r,0}=C_r\circ\varphi,\qquad\Diag\circ\varphi=\CT_y,\qquad\varphi(f)=Q_F,\qquad\varphi(N)=N_F .
\]
Transporting Sections~\ref{sec:setting}--\ref{sec:delta} by $\varphi$, Lemmas~\ref{lem:delta}, \ref{lem:inverse}, \ref{lem:onestep} and Proposition~\ref{prop:levels} hold in $\widehat A_+$ (with its $2$-adic topology), with $f$ in place of $Q$ and $\Lambda_{r,0}$ in place of $C_r$. In particular, for a Laurent polynomial $M\in A_+$ and the all-zero word,
\begin{equation}\label{eq:zerolevels}
C^k(M/f)=\sum_{j\ge0}M^{(k)}_j\sigma^j(1/f),\qquad m_k=\sum_{t=1}^kE_t(m_{k-t}),\qquad M^{(k)}_j=\sum_{t=1}^kC^t(\delta_{j+t}m_{k-t})\ (j\ge1),
\end{equation}
where $E_t(M)=C^t(\delta_tM)$, $\delta_t=\delta_t(f)$ and $m_0=M$. Also $\CT_yR=F$ is the Furstenberg series.

\begin{theorem}\label{thm:zero}
For every $\alpha\ge1$:
\begin{enumerate}[label=\textup{(\alph*)}]
\item the sequence $(C^kR\bmod2^\alpha)_{k\ge0}$ has exact preperiod $\alpha+1$ and exact period $2^{\alpha+1}$; it takes exactly $2^{\alpha+1}+\alpha+1$ distinct values;
\item the sequence $(\Lambda_0^kF\bmod2^\alpha)_{k\ge0}$ has exact eventual period $2^{\alpha+1}$; in particular $K_\alpha(F)\ge2^{\alpha+1}$;
\item let $\widetilde f\in\Z_2[x,y,y^{-1}]$ satisfy $\widetilde f\equiv f\pmod2$ and have the same monomial support as $f$, let $q=2^{\alpha-1}$, $S_0=Nf^{q-1}\bmod2^\alpha$ and $\lambda_{0,0}(S)=\Lambda_{0,0}(S\widetilde f^{\,q})\bmod2^\alpha$ (the initial state and the digit-$0$ transition of \cite[\S2]{RY26} with the lift $\widetilde f$). Then the $\lambda_{0,0}$-orbit of $S_0$ has exactly $2^{\alpha+1}+\alpha+1$ elements.
\end{enumerate}
For $P=x^2y^2+(x^2+x+1)y+x^2$, statement (c) is the orbit-size formula conjectured in the appendix of \cite{RY26}, where it was verified for $\alpha\le14$.
\end{theorem}

The initial state and transition in part (c) are given in \cite[\S2, (2) and the notation on p.~8]{RY26}. The conjecture and the computed range appear in the final paragraph of their appendix before the tables \cite[Appendix, p.~47]{RY26}. The proof occupies the rest of this appendix. The main step is the annihilation estimate of Proposition~\ref{prop:annihilation}.

\subsection{A rank-three subspace with invertible first matrix}\label{subsec:rank-three}
The Newton polygon of $f$ is the triangle $\Delta$ with vertices $(0,0)$, $(2,-1)$ and $(2,1)$. Its lattice points are $(0,0)$, $(1,0)$, $(2,-1)$, $(2,0)$ and $(2,1)$. The linear form $\ell(i,j)=i+2j$ is nonnegative on $\Delta$ and vanishes exactly on the edge $\mathcal E=[(0,0),(2,-1)]$. Put
\[
J_0=(\Delta\cap\Z^2)\setminus\mathcal E=\{(1,0),(2,0),(2,1)\},
\]
so that $L_{J_0}$ is spanned by $x$, $x^2$ and $x^2y$. For $u\in J_0$ and $e'\in(2^s-1)\Delta$, a lattice point $(u+e')/2^s$ lies in $\Delta$ and satisfies $\ell\ge2^{-s}\ell(u)>0$, hence lies in $J_0$. Therefore, by Lemma~\ref{lem:delta}(a), the maps $E_t$ and $M\mapsto C^s(f^{2^s-1}M)$ preserve $L_{J_0}$. Let $\Gamma_t$ and $\Theta_s$ be their matrices with respect to the ordered basis $(x,x^2,x^2y)$, acting on row vectors: $M$ with row vector $[M]$ is sent to the polynomial with row vector $[M]\Gamma_t$, respectively $[M]\Theta_s$. Then $\Gamma_t\equiv0\pmod{2^{t-1}}$, and $(\Theta_s)_{u,w}$ is the coefficient of $x^{2^sw_1-u_1}y^{2^sw_2-u_2}$ in $f^{2^s-1}$.

By \eqref{eq:deltadef} and Lemma~\ref{lem:sections}(a),
\[
C^s(f^{2^s-1}M)=\sum_{t=1}^sC^{s-t}\bigl(f^{2^{s-t}-1}E_t(M)\bigr),\qquad\text{that is,}\qquad\Theta_s=\sum_{t=1}^s\Gamma_t\Theta_{s-t}\quad(\Theta_0=I).
\]
In terms of generating series in a formal variable $Z$, $\bigl(I-\sum_t\Gamma_tZ^t\bigr)\bigl(\sum_s\Theta_sZ^s\bigr)=I$. Writing $H(Z)=\sum_t\Gamma_tZ^t$, the series $\sum_{n\ge0}H(Z)^n$ is defined coefficientwise, since $H(0)=0$, and is a two-sided inverse of $I-H(Z)$. Thus also
\begin{equation}\label{eq:rightconv}
\Theta_s=\sum_{t=1}^s\Theta_{s-t}\Gamma_t\qquad(s\ge1).
\end{equation}
With \eqref{eq:zerolevels} and induction, this gives $[m_k]=[M]\Theta_k$ for $M\in L_{J_0}$.

The matrix
\[
\Theta_1=\Gamma_1=\begin{pmatrix}a&0&0\\1&d&0\\0&c&b\end{pmatrix}
\]
is invertible modulo $2$, and $\Theta_s\equiv\Theta_1^s\pmod2$ because $\Gamma_t\equiv0\pmod 2$ for $t\ge2$. Hence $\Theta_s\in GL_3(\Z_2)$. Put $\Phi_s=\Theta_{s+1}\Theta_s^{-1}$.

\begin{lemma}\label{lem:ratio}
$\Phi_s-\Phi_{s-1}\in2^sM_3(\Z_2)$ for $s\ge1$. Hence $\Phi=\lim\Phi_s$ exists and $\Phi-\Phi_s\in2^{s+1}M_3(\Z_2)$. If $\chi(t)=\det(tI-\Phi)=\sum_{r=0}^3\chi_rt^r$, then $\sum_{r=0}^3\chi_r\Theta_{n+r}\in2^{n+1}M_3(\Z_2)$ for every $n\ge0$.
\end{lemma}

\begin{proof}
By \eqref{eq:rightconv}, $\Theta_{s+1}-\Phi_{s-1}\Theta_s=\sum_{t=1}^s(\Theta_{s+1-t}-\Phi_{s-1}\Theta_{s-t})\Gamma_t+\Gamma_{s+1}$. Using $\Theta_{s+1-t}=\Phi_{s-t}\Theta_{s-t}$, this equals $\sum_{t=2}^s(\Phi_{s-t}-\Phi_{s-1})\Theta_{s-t}\Gamma_t+\Gamma_{s+1}$. By induction on $s$, $\Phi_{s-t}-\Phi_{s-1}=-\sum_{i=s-t+1}^{s-1}(\Phi_i-\Phi_{i-1})\in2^{s-t+1}$, while $\Gamma_t\in2^{t-1}$. So the right-hand side lies in $2^s$, and multiplying by $\Theta_s^{-1}$ gives the first claim. Then $\Theta_{n+r}=\Phi_{n+r-1}\cdots\Phi_n\Theta_n\equiv\Phi^r\Theta_n\pmod{2^{n+1}}$, and $\chi(\Phi)=0$ gives $\sum_r\chi_r\Theta_{n+r}\equiv\chi(\Phi)\Theta_n=0\pmod{2^{n+1}}$. To see the matrix identity directly, write $\operatorname{adj}(tI-\Phi)=t^2I+tB_1+B_0$ in the adjugate identity. Coefficient comparison gives $B_1=\Phi+\chi_2I$, $B_0=\Phi^2+\chi_2\Phi+\chi_1I$, and $-\Phi B_0=\chi_0I$, hence $\chi(\Phi)=0$.
\end{proof}

Vlasenko proves a congruence for successive matrix quotients indexed by the interior lattice points \cite[Theorem~1(ii)]{Vl18} and discusses the full set of lattice points in \cite[Remark~6]{Vl18}. Lemma~\ref{lem:ratio} uses the index set $J_0$, with the identity as Frobenius on $\Z_2$. In the face topology of \cite[Proposition~3.4]{BV21}, deleting the closed edge $\mathcal E$ gives the open subset $\Delta\setminus\mathcal E$. Its integral points are precisely $J_0$, the indexing convention of \cite[Definition~3.5]{BV21}.

\subsection{The eigenvalues}
For $u=(1,0)$ and $w\in\{(2,0),(2,1)\}$, the exponent $2^sw-u$ has $x$-degree $2^{s+1}-1>\deg_xf^{2^s-1}$. For $u=(2,0)$ and $w=(2,1)$, it is $(2^{s+1}-2,2^s)$, whose $y$-degree exceeds that of the top $x$-degree part $(x^2(by+d+cy^{-1}))^{2^s-1}$. Hence $\Theta_s$ is lower triangular, with diagonal entries
\[
\beta_s=[x^{2^s-1}y^0]f^{2^s-1},\qquad\vartheta_s=\CT_y(by+d+cy^{-1})^{2^s-1},\qquad b^{2^s-1}.
\]
So $\Phi$ is lower triangular, with diagonal entries $\mu_1=\lim\beta_{s+1}/\beta_s$, $\mu_2=\lim\vartheta_{s+1}/\vartheta_s$ and $\mu_3=\lim b^{2^s}=1$. They are odd, since $\Phi\equiv\Theta_1\pmod2$.

\begin{lemma}\label{lem:lambda2}
$\mu_2=-1$.
\end{lemma}

\begin{proof}
Let $\kappa=1+(d/c)y+(b/c)y^2\in\Z_2[y]$, so that $1/(by+d+cy^{-1})=(y/c)/\kappa=\sum_{n\ge1}u_ny^n$. We apply the one-variable case ($m=1$) of the constructions above to $\kappa$, with the digit-$0$ operators $E_t(M)=C_y^t(\delta_t(\kappa)M)$. The Newton polytope of $\kappa$ is $[0,2]$, and $\{1\}$ is stable under these operators, since $(1+e')/2^t\in[2^{-t},2-2^{-t}]$ for $e'\in[0,2(2^t-1)]$. The corresponding $1\times1$ matrices are $\theta_s=[y^{2^s-1}]\kappa^{2^s-1}=c^{1-2^s}\vartheta_s$, so their ratio limit is $\mu_2\cdot\lim c^{-2^s}=\mu_2$. Lemma~\ref{lem:ratio} and the proof of Proposition~\ref{prop:annihilation}(i) below, which uses only Lemma~\ref{lem:ratio} and \eqref{eq:zerolevels}, apply verbatim with $M=y/c$. They give $(C_y-\mu_2)C_y^k(y/(c\kappa))\in2^{k+1}\Z_2[[y]]$. The coefficient of $y$ in this series is $u_{2^{k+1}}-\mu_2u_{2^k}$.

Modulo $2$, $\sum u_ny^n\equiv y/(1+y+y^2)=(y+y^2)/(1+y^3)$, so $u_{2^k}$ is odd. Hence $\mu_2=\lim u_{2^{k+1}}/u_{2^k}$. Work in the complete ring $\mathcal O=\Z_2[T]/(T^2+T+1)$, with $\zeta$ the class of $T$ and residue field $\F_4$. The polynomial $bY^2+dY+c$ has the two simple residue roots $\zeta,\zeta^2$. Each lifts uniquely: a root $r_n$ modulo $2^n$ is corrected to $r_n+2^n\epsilon$ modulo $2^{n+1}$ by choosing $\epsilon$ in the residue field, since the derivative $2br_n+d$ is a unit. Completeness supplies roots $\rho_1=\zeta\xi_1$, $\rho_2=\zeta^2\xi_2$ with $\xi_1,\xi_2\in1+2\mathcal O$. By (F2), $\xi_i^{2^k}\equiv1\pmod{2^{k+1}}$. Partial fractions give $u_n=(\rho_2^{-n}-\rho_1^{-n})/(b(\rho_1-\rho_2))$. As $k\to\infty$ we have $\xi_1^{2^k},\xi_2^{2^k}\to1$, so $\rho_2^{-2^k}-\rho_1^{-2^k}$ is asymptotic to $\zeta^{-2\cdot2^k}-\zeta^{-2^k}$. This equals $\zeta-\zeta^2$ if $2^k\equiv1\pmod3$ and $-(\zeta-\zeta^2)$ if $2^k\equiv2\pmod3$. Since $(\zeta-\zeta^2)^2=-3$ is a unit and $2^{k+1}\equiv-2^k\pmod3$, the ratio $u_{2^{k+1}}/u_{2^k}$ tends to $-1$.
\end{proof}

Thus $\chi(t)=(t^2-1)(t-\mu)$ with $\mu=\mu_1\in\Z_2$ odd. We do not need the value of $\mu$. For the five odd integer parameter sets in Table~\ref{tab:main} and Table~\ref{tab:other}, exact computations give $\mu\equiv91\pmod{128}$ and $\mu^2-\mu+2\equiv0\pmod{128}$.

\subsection{Annihilation}

\begin{proposition}\label{prop:annihilation}
\begin{enumerate}[label=\textup{(\roman*)}]
\item For $M\in L_{J_0}$ and $k\ge0$: $\chi(C)C^k(M/f)\in2^{k+1}\widehat A_+$.
\item For $k\ge1$: $\chi(C)C^kR\in2^{k-1}\widehat A_+$.
\end{enumerate}
\end{proposition}

\begin{proof}
By \eqref{eq:zerolevels} and linearity, the level-$j$ numerator of $\chi(C)C^k(M/f)$ is $\sum_r\chi_rM^{(k+r)}_j$.

(i) At level $0$ it equals $[M]\sum_r\chi_r\Theta_{k+r}\in2^{k+1}$ by Lemma~\ref{lem:ratio}. At level $j\ge1$ it is
\[
\sum_{t\ge1}C^t\Bigl(\delta_{j+t}\sum_{r:\,k+r\ge t}\chi_rm_{k+r-t}\Bigr).
\]
For $t\le k$ the inner sum is the level-$0$ residual at $k-t$, which lies in $2^{k-t+1}$; multiplied by $\delta_{j+t}\in2^{j+t-1}$ it lies in $2^{k+j}\subset2^{k+1}$. For $t>k$ we have $\delta_{j+t}\in2^{j+k}\subset2^{k+1}$.

(ii) Now $M=N$ and $m_k=m_k(N)$. The induction in the proof of Lemma~\ref{lem:stable}, with $v=0$ and $\eta=0$, shows that $\supp m_k\subset(1-2^{-k})\Delta+2^{-k}K_N$, where $K_N=\conv\supp N\subset\conv\{(0,1),(1,1),(2,1),(2,2)\}$. The facets of $\Delta$ are $i\le2$, $2j\le i$ and $-2j\le i$. On $K_N$ the first and third inequalities hold, and $2j-i\le2$. So a lattice point of $(1-2^{-k})\Delta+2^{-k}K_N$ satisfies $2j-i\le2^{1-k}$, hence $2j-i\le0$ for $k\ge2$; it also satisfies $\ell\ge2^{-k}\min_{K_N}\ell>0$. Hence $m_k\in L_{J_0}$ for $k\ge2$.

For $k\ge2$ write $m_k=\sum_{t=1}^{k-2}E_t(m_{k-t})+s_k$ with $s_k=E_{k-1}(m_1)+E_k(m_0)$. Then $s_k\in2^{k-2}$, and $s_k\in L_{J_0}$ because all other terms are. Induction with \eqref{eq:rightconv} gives $[m_k]=\sum_{l=2}^k[s_l]\Theta_{k-l}$ for $k\ge2$. For $k\ge2$, the level-$0$ residual is therefore
\[
\sum_r\chi_r[m_{k+r}]=\sum_{l=2}^k[s_l]\sum_r\chi_r\Theta_{k-l+r}+\sum_r\chi_r\sum_{l=k+1}^{k+r}[s_l]\Theta_{k+r-l}.
\]
In the first sum each term lies in $2^{l-2}\cdot2^{k-l+1}=2^{k-1}$, and in the second each term lies in $2^{l-2}\subset2^{k-1}$. At a level $j\ge1$ we argue as in (i). For $t\le k-2$ the inner sum lies in $2^{k-t-1}$, and the product with $\delta_{j+t}$ lies in $2^{k+j-2}\subset2^{k-1}$. For $t\ge k-1$ we have $\delta_{j+t}\in2^{j+k-2}\subset2^{k-1}$. For $k=1$ the claim is trivial.
\end{proof}

\subsection{Reduction modulo 2}
Modulo $2$ all parameters may be replaced by $1$, and $f^2\equiv\sigma(f)$. Hence $C(M/f)\equiv C(Mf)/f\pmod2$, and $C^kR\equiv M_k/f$ with $M_0=N$ and $M_{k+1}=C(fM_k)$. A direct computation gives, modulo $2$,
\[
M_0=(1+x+x^2)y,\qquad M_1=x+xy+x^2+x^2y,\qquad M_2=x^2y,
\]
\[
M_3=x^2(1+y),\qquad M_4=x+x^2y,\qquad M_5=x+x^2(1+y),
\]
and $M_6=M_2$. Hence $(C^2-1)C^kR\equiv x/f\pmod2$ for $k\ge2$, and the sequence $(C^kR\bmod2)_{k\ge2}$ has exact period $4$. The coefficient of $x^1y^0$ in $x/f$ is $1$. Since this monomial lies in the central row, $(\Lambda_0^2-1)\Lambda_0^kF$ also has odd coefficient of $x$ for $k\ge2$, and $(\Lambda_0^kF\bmod2)_{k\ge2}$ has exact period $4$.

\subsection{The preperiod witness}
Let $\widehat F=c/(bF)$ be the second root of $P$ in $y$, so that $P=bx^2(y-F)(y-\widehat F)$ in the ring $\Z_2[y,y^{-1}]((x))$, where $x$ is invertible. In this ring
\[
R=\frac{y^2\,\partial_yP}{P}=y^2\Bigl(\frac1{y-F}+\frac1{y-\widehat F}\Bigr)=\sum_{n\ge0}F^ny^{1-n}-\sum_{n\ge0}y^{n+2}\widehat F^{-n-1},
\]
both series converging $x$-adically, since $F$ and $1/\widehat F=bF/c$ have positive $x$-order; both lie in $A$. Hence $[y^j]R=-(bF/c)^{j-1}$ for $j\ge2$. Write $F=-cx^2V$. Then $P_1V=1+bcx^4V^2$, so $V=1-ax+(a^2-d)x^2+O(x^3)$. For $k\ge1$ and $n=2^k$,
\[
B_k:=[x^2y]C^kR=[x^{2n}y^n]R=b^{n-1}[x^2]V^{n-1}=b^{n-1}(n-1)\bigl(\tfrac n2a^2-d\bigr).
\]
For $\ell\ge1$ and $n'=2^\ell n$, put $U(n)=(n-1)(\frac n2a^2-d)$. Then
\[
U(n')-U(n)=2^{k-1}(2^\ell-1)\bigl(a^22^k(2^\ell+1)-(2d+a^2)\bigr)
\]
has $2$-adic valuation exactly $k-1$, while $b^{n'}-b^n\in2^{k+2}$ because $b^{2^k}\equiv1\pmod{2^{k+2}}$: an odd $b$ has $b^2\equiv1\pmod8$, and each further squaring raises this precision by at least one, by (F2). Since $B_{k+\ell}-B_k=b^{-1}\bigl(b^{n'}(U(n')-U(n))+(b^{n'}-b^n)U(n)\bigr)$, we get
\[
v_2(B_{k+\ell}-B_k)=k-1\qquad(k,\ell\ge1).
\]

\subsection{Proof of Theorem~\ref{thm:zero}}\label{subsec:proof-zero}
(a) Let $S^*=C^{\alpha+1}R\bmod2^\alpha$. By Proposition~\ref{prop:annihilation}(ii), $\chi(C)S^*=0$. For $m\ge1$ let $g_m(t)=\sum_{j=0}^{m-1}t^{2j}$. Then
\[
t^{2m}-1-g_m(\mu)(t^2-1)=(t^2-1)\bigl(g_m(t)-g_m(\mu)\bigr)
\]
is divisible by $\chi(t)$, so $(C^{2m}-1)S^*=g_m(\mu)(C^2-1)S^*$. Write $m=2^ru'$ with $u'$ odd. Then $g_m(\mu)=g_{u'}(\mu)\prod_{i<r}(1+\mu^{2u'2^i})$, where $g_{u'}(\mu)$ is a sum of $u'$ odd terms and each factor is $\equiv2\pmod4$. Hence $v_2(g_m(\mu))=v_2(m)$.

Taking $m=2^\alpha$ shows that $S^*$ is periodic with period dividing $2^{\alpha+1}$, so the preperiod is at most $\alpha+1$. Let $\varpi$ be the exact period. If $\varpi$ were odd, then reduction modulo $2$ would give an odd period of $(C^kR\bmod2)_{k\ge2}$, contradicting the exact period $4$. So $\varpi=2m$, and $0=(C^{2m}-1)S^*=g_m(\mu)(C^2-1)S^*$. The coefficient of $x$ in $(C^2-1)S^*$ is a unit, so $2^\alpha\mid g_m(\mu)$, hence $2^\alpha\mid m$ and $\varpi=2^{\alpha+1}$.

The witness shows that $C^\alpha R\not\equiv C^{\alpha+\ell}R\pmod{2^\alpha}$ for every $\ell\ge1$, so $C^\alpha R$ is not periodic and the preperiod is exactly $\alpha+1$. The orbit of $R$ under $C$ consists of the non-periodic elements $C^kR$ with $k\le\alpha$, followed by a cycle of length $2^{\alpha+1}$. By minimality of the preperiod and of the period, these $2^{\alpha+1}+\alpha+1$ elements are pairwise distinct.

(b) Put $e_k=\Lambda_0^kF\bmod2^\alpha=\CT_y(C^kR)\bmod2^\alpha$. Applying $\CT_y$, the argument of (a) uses three facts:
\begin{enumerate}[label=(\roman*)]
\item $\chi(\Lambda_0)e_{\alpha+1}=0$, by Proposition~\ref{prop:annihilation}(ii) and $\CT_y\circ C=\Lambda_0\circ\CT_y$;
\item for $k\ge2$ the coefficient of $x$ in $(\Lambda_0^2-1)e_k$ is odd;
\item $(e_k\bmod2)_{k\ge2}$ has exact period $4$.
\end{enumerate}
Facts (ii) and (iii) are proved in the subsection on reduction modulo $2$. With them, the argument of (a) shows that the eventual period of $(e_k)$ is exactly $2^{\alpha+1}$. Hence the kernel of $F\bmod2^\alpha$ contains at least $2^{\alpha+1}$ distinct elements $e_k$. The preperiod of $(e_k)$ can be smaller than $\alpha+1$, because the witness monomial $x^2y$ does not lie in the central row.

(c) Since $\widetilde f\equiv f\pmod2$, facts \ref{F1} and \ref{F2} give $\widetilde f^{\,q}\equiv f^q\pmod{2^\alpha}$ and $\widetilde f^{\,2q}\equiv\sigma(\widetilde f)^q\pmod{2^\alpha}$. Hence, by Lemma~\ref{lem:sections}(a) in $A$,
\[
C(S/\widetilde f^{\,q})=C(S\widetilde f^{\,q}/\widetilde f^{\,2q})\equiv C(S\widetilde f^{\,q})/\widetilde f^{\,q}\pmod{2^\alpha}.
\]
Division by the unit $\widetilde f^{\,q}$ is injective. It maps $S_0$ to $N/f=R$ modulo $2^\alpha$ and intertwines $\lambda_{0,0}$ with $C$. So the $\lambda_{0,0}$-orbit of $S_0$ is in bijection with the orbit of $R$, and (a) applies. \qed

\section{The rational lattice and the contact certificate}\label{app:odd-certificate}\label{app:odd-lattice}

This appendix proves the contact estimate used in
Section~\ref{sec:odd-exponent}. Throughout the appendix the parameters
are $a=b=c=d=1$. Put
\[
 q=1+x+x^2,\qquad D=q^2-4x^4,\qquad B=D^{-1/2},\qquad U=(1,x).
\]
These are the $q$, $D$ and $B$ of the proof of Proposition~\ref{prop:odd-universality}. The operators $C_r,T_r$, the frames $P_w,Q_w$, the matrix $H$ and the starred polynomials $C_*,H_*,J_*,T_*,h_*$ below are local to this appendix; in particular $C_r$ and $T_r$ are not the operators $C$ and $T_{t,v}$ of Sections~\ref{sec:delta}--\ref{sec:count}, and $P_w,Q_w$ are not the polynomials $P,Q$.
The square root has constant coefficient $1$. For a scalar series $h$,
let $L(h)=([x^0]h,[x^1]h)^{\mathsf t}$, and apply $L$ columnwise
to a two-column frame. Write
\[
 \pi_Uh=h-U L(h),\qquad
 \mathcal N(P)=P(LP)^{-1}
\]
when $LP$ is invertible. All valuations below are the minimum of the
$2$-adic valuations of the coefficients, including both columns for a
frame. We call this the physical valuation, to distinguish it from classes in the quotient $\mathcal Q$ below. A binary word is written in order of application of its sections.
Define
\[
 A=B(x)/B(x^2),\quad C_r(h)=\Lambda_r(Ah),\quad
 P_\varnothing=U,\quad P_{wr}=\mathcal N(C_rP_w).
\]
These frames are the rational normalizations of the physical frames
$E_w$ of Section~\ref{sec:odd-exponent}. More explicitly, with
\[
 H=L(BU)=\begin{pmatrix}1&0\\-1&1\end{pmatrix},
 \qquad E_w=B P_wH^{-1},
\]
the first two rows of $E_w$ are the identity. Indeed
$\Lambda_r(Bh)=B C_r(h)$, and this identity intertwines the two
normalizations. Multiplication by $B$ and $B^{-1}$, and right
multiplication by a constant integral unit matrix, preserve valuation.
Consequently $v_2(E_w-E_v)=v_2(P_w-P_v)$.

\subsection{An all-precision rational multiplier}

\begin{lemma}\label{lem:odd-ratio}
In the $2$-adic completion of $\Z_2[x,q^{-1}]$ one has
\begin{equation}\label{eq:odd-ratio-expansion}
 A=q+\frac{2T_*}{q}+\sum_{n\ge2}\frac{2^nA_n}{q^{4n-5}},
 \qquad T_*=-x-x^2-x^3+x^4,
\end{equation}
where $A_n\in\Z_2[x]$, $\deg A_n\le8n-8$. In particular every
summand has growth at infinity at most $2$. The identity
$\Lambda_r(Bh)=B C_r(h)$ holds for every integral formal series $h$.
\end{lemma}

\begin{proof}
Set
\[
 C_*=1+x^2+3x^4,\qquad
 H_*=x^5+x^6+2x^7+4x^9+2x^{10}+3x^{11}-3x^{12}.
\]
Polynomial multiplication gives
$C_*^2D-q^2D(x^2)=8H_*$. Therefore, on the branches with constant
coefficient $1$,
\[
 A=\frac{C_*}{q}\left(1-\frac{8H_*}{DC_*^2}\right)^{1/2}.
\]
For $k\ge1$,
$v_2(8^k\binom{1/2}{k})=k+s_2(k)$, where $s_2(k)$ is the number
of nonzero binary digits of $k$. To check this, use
$\binom{1/2}{k}=(-1)^{k-1}\binom{2k}{k}/(4^k(2k-1))$ and
$v_2(n!)=n-s_2(n)$. This also proves convergence in the indicated
completion. It agrees coefficientwise with the formal-series
expansion, since $H_*$ has positive order.

Put $h_*=x^5(1+x)^3$ and $J_*=(H_*-q^2h_*)/2$. The latter is integral,
since $H_*\equiv q^2h_*\pmod2$. As $C_*=q^2+2T_*$,
\[
 \frac{H_*}{DC_*^2}=
 \left(\frac{h_*}{q^4}+\frac{2J_*}{q^6}\right)
 (1-4x^4/q^2)^{-1}(1+2T_*/q^2)^{-2}.
\]
In the $k$th square-root term choose $j$ factors $2J_*/q^6$ and
indices $l,m$ in the last two geometric expansions. Before
multiplication by $C_*/q$ the pole exponent and coefficient valuation
satisfy
\[
 e=4k+2j+2l+2m,\qquad
 n\ge k+s_2(k)+j+2l+m,
\]
hence $e\le4n-4$. Multiplication by the $q$ part of $C_*/q$ leaves
pole at most $4n-5$; its $2T_*/q$ part, with one more unit of valuation,
leaves pole at most $4n-7$ at the new weight. The $k=0$ term is
$q+2T_*/q$. Every displayed factor has growth at infinity at most zero,
except $C_*/q$, which has growth at most two. Collecting at each
valuation therefore gives \eqref{eq:odd-ratio-expansion}, including
its degree bound. At any fixed valuation only finitely many
$k,j,l,m$ contribute, so this collection controls the entire tail.
Finally $Bh=B(x^2)Ah$ and
$\Lambda_r(S(x^2)T(x))=S(x)\Lambda_r(T(x))$ prove the operator
identity.
\end{proof}

\subsection{The normal lattice and the difference cocycle}

For $e\ge0$ use the pole spaces of the proof of Proposition~\ref{prop:odd-universality},
\[
 \mathcal R_e=\{N/q^e:N\in\Z_2[x],\ \deg N\le2e+1\}.
\]
Thus $\mathcal R_0=\langle1,x\rangle$. Define closed normal modules
\begin{align*}
 \mathcal L&=\ker L\cap\left\{\sum_{j\ge0}2^jh_j:
                           h_j\in \mathcal R_{1+2j}\right\},\\
 \mathcal H_3&=\ker L\cap\left\{\sum_{j\ge3}2^jh_j:
                           h_j\in \mathcal R_{1+2j}\right\},\qquad
 \mathcal Q=\mathcal L/(2\mathcal L+\mathcal H_3).
\end{align*}
All sums converge $2$-adically coefficientwise. Each summand in a
normal element may be made normal by applying $\pi_U$: the removed
polynomial pieces sum to $UL(h)=0$.

\begin{lemma}\label{lem:odd-lattice}
The divided operators $T_r=\pi_U C_r/2$ preserve $\mathcal L$ and
send $\mathcal H_3$ into $2\mathcal L$. The frames satisfy
$P_w=U+2Q_w$ with $Q_w\in\mathcal L$. If
$Q_r=(\mathcal N(C_rU)-U)/2$, then
\begin{equation}\label{eq:odd-forcing-congruence}
 Q_{wr}\equiv Q_r\pmod{2\mathcal L}.
\end{equation}
The quotient $\mathcal Q$ is ten-dimensional over $\F_2$, with
ordered basis
\begin{equation}\label{eq:odd-quotient-basis}
 \frac{x^2}{q},\frac{x^3}{q};\qquad
 \frac{2x^i}{q^3}\ (4\le i\le7);\qquad
 \frac{4x^i}{q^5}\ (8\le i\le11).
\end{equation}
If two normalized frames differ by $2^sZ$, where $s\ge1$ and
$Z$ has columns in $\mathcal L$, then after a common digit $r$
the difference is $2^{s+1}Z'$, with
\begin{equation}\label{eq:odd-finite-cocycle}
 [Z']=[T_rZ-Q_rLC_rZ]M_r^{-1},\qquad
 M_0=\begin{pmatrix}1&0\\1&1\end{pmatrix},\quad
 M_1=\begin{pmatrix}1&1\\0&1\end{pmatrix}
 \quad\text{over }\F_2.
\end{equation}
The expression in brackets depends only on $[Z]\in\mathcal Q^2$.
\end{lemma}

\begin{proof}
An exact rational-section formula controls every pole.
For $e\ge1$, in the term indexed by $i$ write
$e-i=2j+\varepsilon$, with $\varepsilon\in\{0,1\}$. The identities
$q(x)q(-x)=q(x^2)$, $q(-x)=q-2x$ and
$q^2=q(x^2)+2xq$ give
\begin{equation}\label{eq:odd-rational-section}
 \Lambda_r\left(\frac N{q^e}\right)=
 \sum_{i=0}^e\sum_{l=0}^j
 \binom ei(-2)^i\binom jl2^l
 \frac{\Lambda_r(Nx^{i+l}q^{\varepsilon+l})}{q^{e-j+l}}.
\end{equation}
If $\deg N\le2e+3$, the degree of each numerator on the right
is at most $2(e-j+l)+1$. Indeed the upper bound before section
is $2e+3+i+3l+2\varepsilon$; after section its degree is at most
the floor of half that number minus $r/2$. The required inequality
follows from $i+l+r\ge1$. In the remaining case $i=l=r=0$ the
upper degree before section is odd, so taking the floor proves the
same inequality. Since
$e-j+l=\lceil(e+i)/2\rceil+l\le\lceil e/2\rceil+i+l$, the
section belongs to
$\sum_{t\ge0}2^t\mathcal R_{\lceil e/2\rceil+t}$.
For $e=0$ a polynomial of degree at most three sections to degree
at most one.

Multiply an input $2^j\mathcal R_{1+2j}$ by the three parts of
\eqref{eq:odd-ratio-expansion}, section, and divide by two. The
resulting weight $g$ and pole bound $p$ are
\[
\begin{array}{c|cc}
 \text{multiplier part}&g&p\\ \hline
 q&j+t-1&j+t\\
 2T_*/q&j+t&j+t+1\\
 2^nA_n/q^{4n-5},\ n\ge2&j+n+t-1&j+2n+t-2.
\end{array}
\]
Every product before section has growth at infinity at most three,
so the degree estimate just proved applies. In each row
$p\le1+2g$; the sole exception $j=t=0$ in the first row is a
polynomial of degree at most one, annihilated by $\pi_U$. This
proves $T_r\mathcal L\subseteq\mathcal L$.
For $j\ge3$ the stronger inequality $p\le2g-1$ holds in all
three rows, proving $T_r\mathcal H_3\subseteq2\mathcal L$.
Projection by $\pi_U$ preserves each pole and growth bound.

The same calculation starting with the polynomial frame $U$ gives
$\mathcal N(C_rU)\in U+2\mathcal L$. Its normalization is a unit:
modulo two its first rows are precisely $M_r$ above. If
$P=U+2Q$ and $P'=\mathcal N(C_rP)$, comparison with
$\mathcal N(C_rU)=U+2Q_r$ gives a difference in $4\mathcal L$.
For, the normal part of $C_rQ$ is in $2\mathcal L$, and replacing
$\pi_U$ by $\pi_{U+2Q_r}$ subtracts $2Q_rLC_rQ$, also in
$2\mathcal L$. Right multiplication by normalization units
preserves the module. This proves the frame assertion and
\eqref{eq:odd-forcing-congruence}.

For the quotient, normal $\mathcal R_e$ has rank $2e$, with numerator
basis $x^2,\ldots,x^{2e+1}$. At weights zero, one and two, the
new ranks after quotienting the previous pole space are
$2$, $6-2=4$ and $10-6=4$. All higher weights are killed by
$\mathcal H_3$, and even coefficients in these ten generators are
killed by $2\mathcal L$. This proves both spanning and
independence in \eqref{eq:odd-quotient-basis}. Equivalently,
multiply by $q^5$ modulo eight and eliminate successively the
pivots $x^2,x^3$ with coefficient one, $x^4,\ldots,x^7$ with
coefficient two, and $x^8,\ldots,x^{11}$ with coefficient four.

For the cocycle put $P'=\mathcal N(C_rP)$. Exact matrix algebra gives
\[
 \mathcal N(C_r(P+2^sZ))-P'
 =2^s(C_rZ-P'LC_rZ)[LC_r(P+2^sZ)]^{-1}.
\]
Substitute $P'=U+2Q'$ and divide by $2^{s+1}$. The result has
columns in $\mathcal L$ by the closure just proved. In the
quotient replace $Q'$ by $Q_r$ using
\eqref{eq:odd-forcing-congruence}, and the right inverse by its
reduction $M_r^{-1}$. This yields
\eqref{eq:odd-finite-cocycle}. The scalar map kills
$2\mathcal L$ in the quotient. It also kills $\mathcal H_3$:
its $T_r$ part lies in $2\mathcal L$, and $LC_rZ$ is divisible
by eight for such an input. This proves that the finite map is
well-defined.
\end{proof}

We call $Q_0$ and $Q_1$ the \emph{forcing frames}: by \eqref{eq:odd-forcing-congruence}, modulo $2\mathcal L$ the frame $Q_w$ of a nonempty word depends only on its last digit.

\begin{lemma}\label{lem:odd-seeds}
One has $Q_0-Q_1\in2\mathcal L$. Consequently
$(P_w-P_v)/4$ has columns in $\mathcal L$ for every pair of
nonempty words. Also, for words $w,v$,
\begin{equation}\label{eq:odd-history-erasure}
 P_{wv}-P_v\in2^{|v|+1}\mathcal L^2.
\end{equation}
\end{lemma}

\begin{proof}
The forcing frames have the sharper pole bounds
\[
 Q_r\in \mathcal R_1+\sum_{k\ge1}2^k\mathcal R_{2k}
 \quad\text{columnwise}.
\]
Indeed the $2T_*/q$ term stays in $\mathcal R_1$ after section and division
by two. A term $n\ge2$ of
\eqref{eq:odd-ratio-expansion}, with additional section weight
$t$, has weight $k=n+t-1$ and pole at most $2n-2+t=2k-t$.
More explicitly,
$Q_r=(\pi_U C_rU/2)(LC_rU)^{-1}$, and the positive-weight
terms of the constant inverse multiply only lower-weight rational
terms. They preserve the displayed bounds. Modulo $2\mathcal L$, terms with $t\ge1$
already have the one-pole margin required at positive weight.
Thus at weight $k\ge1$ only the term $n=k+1,t=0$ may retain
pole $2k$. Its incoming pole is $4k-1$, which is odd. At weight
zero the only border is likewise an odd incoming pole, namely one.

To compare these borders, take an arbitrary incoming numerator $N$
at pole $2m-1$ and work over $\F_2$. Write
$Nq=a(x^2)+xb(x^2)$. The two raw sections at denominator $q^m$
have numerator rows $(a,xb)$ and $(b,a)$ for the input row
$(N,xN)$. Multiplication by $M_0^{-1}$ and $M_1^{-1}$ gives
\[
 (a+xb,xb),\qquad (b,a+b).
\]
At a root $\zeta$ of $q$ in $\F_4$, evaluation of the splitting
identity at $\zeta^2$ gives
$a(\zeta)+\zeta^2b(\zeta)=0$. Since
$\zeta^2+\zeta=1$, the two displayed rows agree at $\zeta$,
so their difference is divisible by $q$. The top pole of the
forcing difference therefore drops from $2k$ to at most $2k-1$.
Scalar normalization corrections from weight $k-1$ have pole at
most
\[
 \max\{1,2k-2\}\le2k-1\qquad(k\ge1),
\]
including the boundary $k=1$. At weight zero the pole-free normal
remainder is a polynomial of degree at most one with its first two
coefficients zero, and hence vanishes.

This cancellation lifts integrally: divide a borderline numerator
by the monic polynomial $q$, writing it as $qM+2N_1$. The first
term has the required reduced pole at the current weight; the
second carries the old pole $2k$ to weight $k+1$, where the
allowed bound in $2\mathcal L$ is $2k+1$. At weight zero its
carry is in $2\mathcal R_1$, also allowed. Iterating at every weight proves
$Q_0-Q_1\in2\mathcal L$. Formula~\eqref{eq:odd-forcing-congruence} now shows
that all nonempty $P_w$ are congruent modulo $4\mathcal L^2$.
Finally $P_w-U\in2\mathcal L^2$; each common appended digit
increases this factor by two by the exact difference calculation
in Lemma~\ref{lem:odd-lattice}. Induction proves
\eqref{eq:odd-history-erasure}.
\end{proof}

We call \eqref{eq:odd-history-erasure} \emph{history erasure}: the older digits $w$ change the frame of $v$ only by a multiple of $2^{|v|+1}$.

\subsection{The finite certificate}

Encode a vector in the ordered basis
\eqref{eq:odd-quotient-basis} by the integer whose ten binary
bits are its coordinates, with the first basis vector in bit zero.
This encoding is used only for finite linear algebra; addition is
bitwise exclusive-or. Let $S_r$ be the scalar map
$[z]\mapsto[T_rz-Q_rLC_rz]$ in
\eqref{eq:odd-finite-cocycle}.

\begin{lemma}\label{lem:odd-finite}
The columns of the two $10\times10$ matrices $S_r$, in the stated
encoding, are
\begin{align*}
 S_0&=(1,3,29,0,58,0,28,0,56,0),\\
 S_1&=(2,1,0,29,0,58,0,28,0,56).
\end{align*}
Both forcing frames $Q_r$ have column classes $(477,954)$.
For every $0\le u<4$ and $4\le v<8$, interpreting $u,v$ as
three-digit words with least significant digit applied first,
both columns of $(P_u-P_v)/4$ have class $39$.
\end{lemma}

\begin{proof}
We give bounded exact arithmetic for all entries. With $h_*,J_*$
as in Lemma~\ref{lem:odd-ratio}, the expansion there gives
\begin{equation}\label{eq:odd-A16}
 A\equiv
 \frac{q^8+2T_*q^6-4h_*q^4-8(J_*+T_*h_*)q^2-8h_*^2}{q^7}
 \pmod {16}.
\end{equation}
Here it suffices to retain the terms of the square-root expansion
whose valuation is less than four; expanding its two geometric
factors gives the displayed expression. Its numerator, by
increasing degree from zero to sixteen, is
\[
 (1,6,6,8,2,14,2,6,1,6,10,14,0,0,12,14,15)\pmod {16}.
\]
For any numerator $f$ at pole $e$, an alternative exact form of
\eqref{eq:odd-rational-section} is
\[
 \Lambda_r(f/q^e)
 =\frac{\Lambda_r(fq(-x)^e)}{q^e}.
\]
Apply this to \eqref{eq:odd-A16} times each of the ten basis
vectors and to its two polynomial columns $U$. The first-row
matrices of $C_rU$ modulo sixteen are
\[
 LC_0U=\begin{pmatrix}1&0\\1&15\end{pmatrix},\qquad
 LC_1U=\begin{pmatrix}15&1\\0&1\end{pmatrix}.
\]
Subtract the first two coefficients, divide by two, and reduce to
\eqref{eq:odd-quotient-basis}. This yields the following complete
column calculation. The columns headed $T_r$ are encoded in the
same basis; the only nonzero columns of $LC_r$ modulo two are
$LC_0(e_1)=(0,1)^{\mathsf t}$ and
$LC_1(e_2)=(0,1)^{\mathsf t}$.
\[
\begin{array}{c|rrrrrrrrrr}
 &e_1&e_2&e_3&e_4&e_5&e_6&e_7&e_8&e_9&e_{10}\\ \hline
 T_0&955&3&29&0&58&0&28&0&56&0\\
 T_1&2&955&0&29&0&58&0&28&0&56\\
 Q_0\text{ and }Q_1&477&954&&&&&&&&
\end{array}
\]
The last row lists the two forcing columns rather than ten map
columns. Thus subtracting the forcing contribution changes
$955$ to $955\mathbin{\mathrm{xor}}954=1$, and leaves all
other entries unchanged. This proves both matrices and the forcing
classes.

For completeness the seed calculation can be performed without
rational square roots. Write $B=\sum b_nx^n$, set $b_0=1$ and
$b_n=0$ for $n<0$, and use $2DB'+D'B=0$ to compute
\begin{equation}\label{eq:odd-radical-recurrence}
 2n b_n=-\sum_{j=1}^4 d_j(2n-j)b_{n-j},\qquad
 (d_1,d_2,d_3,d_4)=(2,3,2,-3).
\end{equation}
This is an exact integer recurrence; integrality also follows
from $B=q^{-1}\sum_{k\ge0}\binom{2k}{k}(x^4/q^2)^k$.
For each $0\le v<8$ form the first twelve coefficients of
\[
 \left(\sum_{n\ge0}b_{8n+v}x^n,
       \sum_{n\ge0}b_{8n+v-1}x^n\right).
\]
Multiply by $B^{-1}$ and normalize its first two rows modulo
thirty-two. Subtract any column for $u<4$ from that for $v\ge4$,
divide by four, and eliminate the pivots in
\eqref{eq:odd-quotient-basis}. The sixteen resulting pairs of
classes are
\[
\begin{array}{c|cccc}
 u\backslash v&4&5&6&7\\ \hline
0&(39,39)&(39,39)&(39,39)&(39,39)\\
1&(39,39)&(39,39)&(39,39)&(39,39)\\
2&(39,39)&(39,39)&(39,39)&(39,39)\\
3&(39,39)&(39,39)&(39,39)&(39,39).
\end{array}
\]
This is a finite certificate: \eqref{eq:odd-radical-recurrence}
requires only $b_0,\ldots,b_{95}$, and inversion and normalization
use only twelve coefficients. Its sufficiency follows from the
previous analytic lemmas. Every normal lattice element modulo
eight has denominator $q^5$ and numerator degree at most eleven.
Since $q^5(0)=1$, its first twelve coefficients determine that
numerator uniquely modulo eight.
A normal lattice element $z$ that vanishes modulo eight lies in $\mathcal H_3$. Indeed, normalizing summands, write $z=h_0+2h_1+4h_2+y$ with normal $h_j=N_j/q^{1+2j}\in\mathcal R_{1+2j}$ and $y\in\mathcal H_3$. Since $y$ and $z$ vanish modulo eight, so does the numerator $q^4N_0+2q^2N_1+4N_2$ of $z-y$. Reducing successively modulo two gives $N_0=2N_0'$, $N_1=-q^2N_0'+2N_1'$ and $N_2=-q^2N_1'+2N_2'$ with integral $N_0',N_1',N_2'$ of degrees at most $3$, $7$ and $11$. Substituting, $z-y=8N_2'/q^5=8q^2N_2'/q^7$, which is normal and lies in $8\mathcal R_7$, hence in $\mathcal H_3$. Therefore the class of a lattice element in $\mathcal Q$ depends only on the element modulo eight.
Thus no possible further
coefficient can alter the quotient calculation. Division by two
for the maps needs precision sixteen; division by four for seeds
needs precision thirty-two, exactly as used above. These formulas,
all matrix entries and the sixteen seeds are reproduced by the
self-contained script in \texttt{anc/exponent\_three/}.
\end{proof}

\begin{remark}\label{rem:odd-seed-structure}
Half of the matrix columns and the sixteen equal seed classes have a structural explanation. First, $\Lambda_1(xg)=\Lambda_0(g)$ gives $C_1(xh)=C_0(h)$. Since $Q_0-Q_1\in2\mathcal L$ (Lemma~\ref{lem:odd-seeds}), it follows that $S_1[xz]=S_0[z]$ whenever $z$ and $xz$ lie in $\mathcal L$. In the basis \eqref{eq:odd-quotient-basis} multiplication by $x$ sends $e_{2i-1}$ to $e_{2i}$, which explains $S_1(e_{2i})=S_0(e_{2i-1})$. Second, for a nonempty word $w$, \eqref{eq:odd-forcing-congruence} gives $P_w=U+2Q$ with $[Q]=(477,954)$. Applying \eqref{eq:odd-finite-cocycle} with $s=1$ to the frames $P_w$ and $U$ gives $P_{wr}=U+2Q_r+4Z_r$ with $[Z_r]=S_r(477,954)M_r^{-1}$, and the columns in Lemma~\ref{lem:odd-finite} give $(3,1)$ for both $r$. Hence for every seed pair $(P_u-P_v)/4$ has the class of $(Q_0-Q_1)/2$, so the sixteen entries of the table are one class, $39$ in each column, computed sixteen times. We keep the full table as a check.
\end{remark}

\begin{lemma}\label{lem:odd-unit}
Every product of the scalar maps $S_0,S_1$ sends the seed $39$
to a nonzero vector. Every actual lattice element whose class is
one of these images has physical valuation at most one.
\end{lemma}

\begin{proof}
Let $\mathsf W=\langle1,2,28,36\rangle$ and
$\mathsf N=\langle1,2\rangle$ in $\mathcal Q$. The four
generators of $\mathsf W$ are independent, and direct addition
of the certified columns gives
\[
\begin{array}{c|rrrr}
 z&1&2&28&36\\ \hline
 S_0z&1&3&39&29\\
 S_1z&2&1&29&58.
\end{array}
\]
Thus both spaces are invariant. In the ordered basis
$[28],[36]$ of $\mathsf W/\mathsf N$ the induced matrices are
\[
 \overline S_0=\begin{pmatrix}0&1\\1&0\end{pmatrix},
 \qquad
 \overline S_1=\begin{pmatrix}1&1\\0&1\end{pmatrix}.
\]
They are invertible. Since $39=36+1+2$ projects to the nonzero
second basis vector, no product annihilates it.

Every vector of $\mathsf W$ has zero coordinates in the final
four basis positions. If at least one of its first two coordinates is nonzero,
any actual representative has a coefficient odd modulo two. If
these coordinates are zero and its class is nonzero, its next
four coordinates represent a nonzero element of
$2\mathcal R_3/(2\mathcal R_1+4\mathcal R_3)$. Terms from $2\mathcal L$ contribute only
$2\mathcal R_1$ modulo four; terms from $\mathcal H_3$ are divisible by
eight. They cannot cancel this nonzero middle-layer class.
The representative therefore has a coefficient nonzero modulo
four. This proves the physical valuation bound. For a
two-column frame the unit matrix on the right of
\eqref{eq:odd-finite-cocycle} preserves the column span, so the
same nonannihilation conclusion applies to its common
one-dimensional seed span.
\end{proof}

\begin{proposition}\label{prop:odd-rational-contact}
For distinct full windows $w,v$ of length $\alpha-1$ with
longest common terminal suffix of length $\ell$, their clipped
physical contact satisfies
\[
 \min\{\alpha,v_2(E_w-E_v)\}
 \le\min\{\alpha,\ell+3\}.
\]
\end{proposition}

\begin{proof}
Only $\ell+3<\alpha$ needs proof. Retain the last $\ell+3$
digits in each word. Their three-digit prefixes have different
last digits and hence belong, in either order, to one of the
sixteen seed pairs in Lemma~\ref{lem:odd-finite}. By
Lemma~\ref{lem:odd-seeds} their difference is $4Z_0$, with
columns in $\mathcal L$ and common seed class $39$.
Append the common $\ell$-digit suffix. Repeated use of
\eqref{eq:odd-finite-cocycle} makes their difference
$2^{\ell+2}Z_\ell$. Lemma~\ref{lem:odd-unit} supplies a
nonzero column class with physical valuation at most one,
so the retained windows have contact at most $\ell+3$.
History erasure~\eqref{eq:odd-history-erasure}, applied to the
retained word of length $\ell+3$, changes either frame by an
element divisible by $2^{\ell+4}$. It therefore preserves
this witness of smaller valuation. The same contact bound
holds for the original windows, and the fixed physical
normalization preserves it. The cases with
$\ell+3\ge\alpha$, including the shorter windows, are
immediate from clipping.
\end{proof}

\section*{Data and code}
The ancillary directory \texttt{anc/} contains the independent rational-Cartier construction, the exact history and polynomial constructions, the support certificates, scripts checking the finite identities used in the proofs, the rational-lattice certificate for Appendix~\ref{app:odd-lattice}, and exact results and run logs. Its \texttt{README.md} gives reproduction commands and explains the graph formats; \texttt{SHA256SUMS} identifies every included file. Large complete graphs are deterministically regenerated by the supplied code, and their SHA-256 values are recorded with the computation results. No probabilistic state identification or approximate polynomial arithmetic is used in these computations.

\end{document}